\documentclass[11pt]{article}
\usepackage{amsmath}
\usepackage{amssymb,amsmath}
\usepackage{mathrsfs}
\usepackage{graphicx}
\usepackage{color}
\usepackage{mathrsfs}
\usepackage{cite}
\usepackage{indentfirst}
\usepackage[titletoc]{appendix}
\usepackage [latin1]{inputenc}
\let\f=\frac
\let\p=\psi

\let\wt=\widetilde

\def\cC{{\cal C}}

\def\cP{{\cal P}}
\def\cQ{{\cal Q}}

\def\cS{{\cal S}}

\def\cZ{{\cal Z}}

\def\p{\partial}

\def\wt{\widetilde}

\def\dv{\mbox{div}}

\def\eqdefa{\buildrel\hbox{\footnotesize def}\over =}

\def\C{\mathop{\bf C\kern 0pt}\nolimits}
\def\DD{\mathop{\bf D\kern 0pt}\nolimits}
\def\K{\mathop{\bf K\kern 0pt}\nolimits}
\def\N{\mathop{\bf N\kern 0pt}\nolimits}
\def\Q{\mathop{\bf Q\kern 0pt}\nolimits}
\def\R{\mathop{\bf R\kern 0pt}\nolimits}
\def\endproof{\hphantom{MM}\hfill\llap{$\square$}\goodbreak}

\def\ddq{\dot \Delta_q}

\renewcommand{\div}{\mbox{\rm div}\;\!}

\newcommand{\Dv}{{\rm div}}

\newcommand{\beq}{\begin{equation}}
\newcommand{\eeq}{\end{equation}}
\newcommand{\ben}{\begin{eqnarray}}
\newcommand{\een}{\end{eqnarray}}
\newcommand{\beno}{\begin{eqnarray*}}
\newcommand{\eeno}{\end{eqnarray*}}

\newtheorem{Theorem}{Theorem}[section]
\newtheorem{Definition}[Theorem]{Definition}
\newtheorem{Proposition}[Theorem]{Proposition}
\newtheorem{Lemma}[Theorem]{Lemma}
\newtheorem{Corollary}[Theorem]{Corollary}
\newtheorem{Remark}[Theorem]{Remark}

\numberwithin{equation}{section}

\allowdisplaybreaks \numberwithin{equation} {section}

\begin{document}
\title{The maximal regularity and its application to a  multi-dimensional  non-conservative viscous compressible two-fluid
model with capillarity effects  in $L^{ p}$-type framework \thanks {Research supported by the
National Natural Science Foundation of China (11501332,11771043,51976112), the  Natural Science Foundation of Shandong Province (ZR2021MA017,ZR2015AL007),
 and Young Scholars Research Fund of Shandong University of Technology.}
}
\author{ Fuyi Xu$^\dag$    \\[2mm]
 { \small  School of Mathematics and  Statistics, Shandong University of Technology,}\\
  { \small Zibo,    255049,  Shandong,    China}}
         \date{}
         \maketitle
\noindent{\bf Abstract}\ The present paper is the continuation of work \cite{XC},  devoted to extending it to a critical functional  framework which is not related to the energy space. Employing  the special dissipative structure of the non-conservative viscous compressible two-fluid
model with capillarity effects, we first exploit the maximal regularity estimates for the corresponding linearized  system in all frequencies which behaves like the heat equation.
Then we  construct the global well-posedness for the multi-dimensional model when the initial data are  close to a stable equilibrium state in the sense of suitable $L^{ p}$-type Besov
norms. As a consequence, this
allows us to work in the framework of Besov space with negative regularity indices and this fact is
particularly important when the initial data are large highly oscillating in physical dimensions $N= 2, 3$. Furthermore, based  on a refined time weighted inequalities in the Fourier spaces,   we also  establish  optimal  time decay rates  for  the  constructed global solutions under a mild additional decay assumption involving only the low frequencies of the initial  data. 
\vskip   0.2cm \footnotetext[1]{$^\dag$Corresponding author.}
\vskip   0.2cm \footnotetext[2]{E-mail addresses: zbxufuyi@163.com(F.Xu).} \setlength{\baselineskip}{20pt}

\section{Introduction and Main Results}
\setcounter{section}{1}\setcounter{equation}{0} \ \ \ \ \ 
 The models of multi-phase flows have a very broad applications of
hydrodynamics in nature and industry,  where the fluids under
investigation contain more than one component. In nature, there is a variety of different multi-phase
flow phenomena, such as sediment transport, geysers, volcanic eruptions, clouds, and
rain (see \cite{GN}). On the other hand, it has been estimated
that over half of anything which is produced in a modern industrial society depends,
to some extent, on a multi-phase flow process for their optimum design and safe operations.  In addition,  multi-phase flows also naturally appear in many contexts
within biology, ranging from tumor biology and anticancer therapies to developmental
biology and plant physiology \cite{MT}. 
Recently, due to the physical importance and mathematical challenges, the study of  mathematical properties for the models becomes  a  significant and difficult topic.

 In the present paper, we are going to study the following  multiphase  flows  model,  namely a non-conservative viscous
compressible two-fluid system with capillarity effects in $\mathbb{R}^{N}(N\geq2)$:
\begin{align}\label{equ:CTFS}
\left\{
\begin{aligned}
&\alpha^{+}+\alpha^{-}=1,\\
&\p_t(\alpha^{\pm}\rho^{\pm})+\textrm{div}(\alpha^{\pm}\rho^{\pm}u^{\pm})=0, \\
&\p_t(\alpha^{\pm}\rho^{\pm}u^{\pm})+\textrm{div}(\alpha^{\pm}\rho^{\pm}
u^{\pm}\otimes u^{\pm})+\alpha^{\pm}\nabla P^{\pm}(\rho^{\pm})-\sigma^{\pm}\alpha^{\pm}\rho^{\pm}\nabla \Delta(\alpha^{\pm}\rho^{\pm})
=\textrm{div}(\alpha^{\pm}\tau^{\pm}), \\
&P^{+}(\rho^{+})=A^{+}(\rho^{+})^{\bar{\gamma}^{+}}=P^{-}(\rho^{-})=A^{-}(\rho^{-})^{\bar{\gamma}^{-}},
\end{aligned}
\right.
\end{align}
where the variable $0\leq\alpha^{+}(x,t)\leq1$ is the volume fraction of fluid $+$ in one of the two gases, and $0\leq\alpha^{-}(x,t)\leq1$  is the  volume fraction of the other fluid$-$. Moreover, $\rho^{\pm}(x,t)\geq0$, $u^{\pm}(x,t)$ and $P^{\pm}(\rho^{\pm})$ are, respectively, the
densities, the velocities, and the two pressure functions of the fluids. Here, for the sake of technical simplicity, we restricted ourselves to the study of the
barotropic case, that is, $P^{\pm}(\rho^{\pm})=A^{\pm}(\rho^{\pm})^{\overline{\gamma}^{\pm}}$  where  $\overline{\gamma}^{\pm}>1,A^{\pm}>0$. In what follows, we set $A^{+}=A^{-}=1$ without loss of any generality. $\sigma^{\pm}>0$ are the capillary coefficients. Also, $\tau^{\pm}$ are the viscous stress tensors
\begin{equation}\label{1.2}
\tau^{\pm}:=2\mu^{\pm}D(u^{\pm})+\lambda^{\pm}\textrm{div}u^{\pm}\textrm{IN},
\end{equation}
where $D(u^{\pm})\eqdefa\frac{\nabla u^{\pm}+\nabla^{t}u^{\pm}}{2}$ stand for the deformation tensor,  the constants  $\mu^{\pm}$ and $\lambda^{\pm}$ are the (given) shear and bulk viscosity coefficients satisfying $\mu^{\pm}>0$ and $\lambda^{\pm}+2\mu^{\pm}>0$.   System \eqref{equ:CTFS} is known as a two-fluid flows model with algebraic closure,  which is widely used in industrial applications, such as nuclear, power, oil-and-gas, micro-technology and so on,  and  we refer readers to Refs \cite{GN,MT,NK1,NK2} for more discussions about this model and related models.

From a mathematical point of view, model \eqref{equ:CTFS} is a highly nonlinear  partial differential   equation  with
 the mixed hyperbolic-parabolic property. As a
matter of fact, there is no diffusion on the mass conservation equations, whereas velocity evolves according to the parabolic equations due to the  viscosity phenomena. In particular, the non-conservative pressure terms $\alpha^{\pm}\nabla P^{\pm}(\rho^{\pm})$ typically prevent one from applying
arguments used for compressible Navier-Stokes equations. Therefore, there are more challenges
associated with this type of model when it comes to mathematical analysis
(well-posedness and stability).
In the last decade, many researchers have been devoted to studying system \eqref{equ:CTFS} and have made more progress. For example, Bresch et al. \cite{BDG} first established the existence of global weak solutions to  3D   model \eqref{equ:CTFS}. Later, Bresch-Huang-Li \cite{BHL} extended the result in \cite{BDG} and proved the existence of global weak solutions to system \eqref{equ:CTFS} in one space dimension without capillarity terms.  In 2016, Evje-Wang-Wen \cite{EWW}   proved the global existence
of strong solutions to system \eqref{equ:CTFS} without capillary
terms    by the standard energy method under the condition that the initial data are
close to the constant equilibrium state in $H^{2}(\mathbb{R}^{3})$,  and then constructed the  optimal time decay rates of the global  strong solutions  if the initial data belong to $L^{1}(\mathbb{R}^{3})$ additionally.
Lai-Wen-Yao \cite{LWY}  investigated the vanishing capillarity limit of  smooth solutions  to   system  \eqref{equ:CTFS}   with  unequal pressure functions  if   $\|(\alpha^{+}\rho^{+}_0-1, \alpha^{-}\rho^{-}_0-1)\|_{H^{4}(\mathbb{R}^{3})} +\|(u^{+}_0, u^{-}_0)\|_{H^{3}(\mathbb{R}^{3})}$ are small enough.
Recently, when $\mu^{+}=\mu^{-}$ and $\lambda^{+}=\lambda^{-},$  based on the complicated  spectral analysis  of  Green's function to the linearized system and the elaborate energy estimate to the nonlinear system,  for system  \eqref{equ:CTFS},   authors \cite{CWYZ} showed  global solvability of  smooth solutions  close to an equilibrium state  in $H^s(\mathbb{R}^{3})(s\geq3)$ and further  got the   time decay rates when the initial
perturbation is bounded in $L^{1}(\mathbb{R}^{3})$. More recently, Li et al. \cite{ LWWZ} constructed  the global existence and optimal decay rates of  system  \eqref{equ:CTFS} with general constant viscosities and capillary coefficients  when initial data $\|(\alpha^{+}\rho^{+}_0-1, \alpha^{-}\rho^{-}_0-1)\|_{H^{l+1}(\mathbb{R}^{3})} +\|(u^{+}_0, u^{-}_0)\|_{H^{l}(\mathbb{R}^{3})}$ for an integer $l\geq3$ are small enough. In $L^{2}$-type  critical  Besov spaces, Xu et al. \cite{XC2} constructed the well-posedness and decay rates of
strong solutions to a  multi-dimensional  system \eqref{equ:CTFS} without  capillarity effects. Later,  Xu and Chi \cite{XC} further established the corresponding conclusions of system  \eqref{equ:CTFS}.


Here,  it should be pointed out that the existing  results mentioned above, including the global existence and time decay rates of  strong solutions to system \eqref{equ:CTFS},   are mainly based on $L^{2}$-framework, especially in critical Besov spaces. However, to our knowledge, so far there is no result on  the  global existence  and time decay estimates to system \eqref{equ:CTFS} in $L^{p}$-framework. The main motivation of this paper is to give a positive answer to this question. In particular, we prove the global well-posedness to system \eqref{equ:CTFS}  when  the initial data are  close to a stable equilibrium state in the sense of suitable $L^{ p}$-type Besov
norms, and establish  optimal  time decay rates  for  the  constructed global solutions under a mild additional decay assumption involving only the low frequencies of the initial  data. Let us emphasize  that this framework allows us to construct global solutions of system \eqref{equ:CTFS} with highly oscillating initial
velocities in larger spaces in physical dimensions $N= 2, 3$.  
However, due to the mixed hyperbolic-parabolic property of system \eqref{equ:CTFS},   as in  \cite{Dan1}, the system has to be handled differently for the low and high frequencies.
Roughly speaking, the first order terms predominate in low frequencies,
so that  system \eqref{equ:CTFS} has to be treated by means of hyperbolic energy methods, which implies that we
must treat the low frequencies regime only in spaces constructed on $L^{2}$, as it is classical that
hyperbolic systems are ill-posed in general $L^{p}$ spaces.  In contrast, in the  high frequencies,  a $L^{ p}$ approach may be used.

For the convenience of the reader, as in \cite{BDG},  we also show some  derivations for
another expression of the pressure gradient in terms of the gradients of $\alpha^{+}\rho^{+}$ and $\alpha^{-}\rho^{-}$
 by using the pressure equilibrium assumption.  Here, we only focus on the
case that $\inf\rho^{\pm}>0, 0<\alpha^{\pm}<1$  in our framework.
 The relation between the
pressures of  system \eqref{equ:CTFS} implies the following differential identities
\begin{equation}\label{1.3}\textrm{d}P^{+}=s_{+}^{2}\textrm{d}\rho^{+},\quad
\textrm{d}P^{-}=s_{-}^{2}\textrm{d}\rho^{-},\quad \hbox{where}\quad
s_{\pm}:=\sqrt{\frac{\textrm{d}P^{\pm}}{\textrm{d}\rho^{\pm}}(\rho^{\pm})}
=\sqrt{\overline{\gamma}^{\pm}\frac{P^{\pm}(\rho^{\pm})}{\rho^{\pm}}},
\end{equation}
where $s_{\pm}$ denote the sound speed of each phase respectively.

Let
\begin{equation}\label{1.4}
R^{\pm}=\alpha^{\pm}\rho^{\pm}.
\end{equation}
Resorting to $\eqref{equ:CTFS}_{1}$,  we have
\begin{equation}\label{1.5}
\textrm{d}\rho^{+}=\frac{1}{\alpha_{+}}(\textrm{d}R^{+}-
\rho^{+}\textrm{d}\alpha^{+}),
~~\textrm{d}\rho^{-}=\frac{1}{\alpha_{-}}(\textrm{d}R^{-}+
\rho^{-}\textrm{d}\alpha^{+}).
\end{equation}
Combining with   \eqref{1.3} and \eqref{1.5}, we  conclude that
\begin{equation*}\label{1.6}
\textrm{d}\alpha^{+}=\frac{\alpha^{-}s_{+}^{2}}
{\alpha^{-}\rho^{+}s_{+}^{2}+\alpha^{+}\rho^{-}s_{-}^{2}}\textrm{d}R^{+}
-\frac{\alpha^{+}s_{-}^{2}}
{\alpha^{-}\rho^{+}s_{+}^{2}+\alpha^{+}\rho^{-}s_{-}^{2}}
\textrm{d}R^{-}.
\end{equation*}
Substituting  the above equality into \eqref{1.5}, we obtain
\begin{equation*}
\textrm{d}\rho^{+}=\frac{\rho^{+}\rho^{-}s_{-}^{2}}
{R^{-}(\rho^{+})^{2}s_{+}^{2}+R^{+}(\rho^{-})^{2}s_{-}^{2}}
\Big(\rho^{-}\textrm{d}R^{+}
+\rho^{+}\textrm{d}R^{-}\Big),
\end{equation*}
and
\begin{equation*}
\textrm{d}\rho^{-}=\frac{\rho^{+}\rho^{-}s_{+}^{2}}
{R^{-}(\rho^{+})^{2}s_{+}^{2}+R^{+}(\rho^{-})^{2}s_{-}^{2}}
\Big(\rho^{-}\textrm{d}R^{+}
+\rho^{+}\textrm{d}R^{-}\Big),
\end{equation*}
which give, for the pressure differential $\textrm{d}P^{\pm}$,
\begin{equation*}
\textrm{d}P^{+}=\mathcal{C}^{2}\big(\rho^{-}\textrm{d}R^{+}
+\rho^{+}\textrm{d}R^{-}\big),
\end{equation*}
and
\begin{equation*}
\textrm{d}P^{-}=\mathcal{C}^{2}\big(\rho^{-}\textrm{d}R^{+}
+\rho^{+}\textrm{d}R^{-}\big),
\end{equation*}
where
\begin{equation*}
\mathcal{C}^{2}\eqdefa\frac{s_{-}^{2}s_{+}^{2}}{\alpha^{-}\rho^{+}s_{+}^{2}
+\alpha^{+}\rho^{-}s_{-}^{2}}.
\end{equation*}
Recalling $\alpha^{+}+\alpha^{-}=1$, we get the following identity:
\begin{equation}\label{1.61}
\frac{R^{+}}{\rho^{+}}+\frac{R^{-}}{\rho^{-}}=1,\quad \hbox{and ~therefore~}
\rho^{-}=\frac{R^{-}\rho^{+}}{\rho^{+}-R^{+}}.
\end{equation}
Then it follows from the pressure relation $\eqref{equ:CTFS}_{4}$ that
\begin{equation}\label{1.7}
\varphi(\rho^{+}):=P^{+}(\rho^{+})-P^{-}(\frac{R^{-}\rho^{+}}
{\rho^{+}-R^{+}})=0.
\end{equation}
 Differentiating $\varphi$ with respect to $\rho^{+}$,   we have
$$\varphi^{'}(\rho^{+})=s^{2}_{+}+s^{2}_{-}\frac{R^{-}R^{+}}{(\rho^{+}-R^{+})^{2}}.$$
By the definition of $R^{+}$, it is natural to look for $\rho^{+}$ which belongs to $(R^{+},+\infty).$  Since $\varphi'>0$ in $(R^{+},+\infty)$  for any given $R^{\pm}>0,$ and $\varphi:(R^{+},+\infty)\longmapsto(-\infty,+\infty),$  this determines that $\rho^{+}=\rho^{+}(R^{+},R^{-})\in(R^{+},+\infty)$ is the unique solution of  equation \eqref{1.7}. Due to \eqref{1.5}, \eqref{1.61} and $\eqref{equ:CTFS}_{1}$, $\rho^{-}$ and $\alpha^{\pm}$  are defined as follows:
$$\rho^{-}(R^{+},R^{-})=\frac{R^{-}\rho^{+}(R^{+},R^{-})}{\rho^{+}(R^{+},R^{-})-R^{+}},$$
$$\alpha^{+}(R^{+},R^{-})=\frac{R^{+}}{\rho^{+}(R^{+},R^{-})},$$
$$\alpha^{-}(R^{+},R^{-})=1-\frac{R^{+}}{\rho^{+}(R^{+},R^{-})}=\frac{R^{-}}{\rho^{-}(R^{+},R^{-})}.$$
Based on  the above analysis,   system \eqref{equ:CTFS} is equivalent to the following form
\begin{align}\label{equ:CTFS1}
\left\{
\begin{aligned}
&\p_tR^{\pm}+\textrm{div}(R^{\pm}u^{\pm})=0, \\
&\p_t(R^{+}u^{+})+\textrm{div}(R^{+}u^{+}\otimes u^{+})+\alpha^{+}\mathcal{C}^{2}[\rho^{-}\nabla R^{+}+\rho^{+}\nabla R^{-}]-\sigma^{+}R^{+}\nabla \Delta R^{+}
\\&\qquad=\textrm{div}\big(\alpha^{+}[\mu^{+}(\nabla u^{+}+\nabla^{t}u^{+})+\lambda^{+}\textrm{div}u^{+}\textrm{Id}]\big), \\
&\p_t(R^{-}u^{-})+\textrm{div}(R^{-}u^{-}\otimes u^{-})+\alpha^{-}\mathcal{C}^{2}[\rho^{-}\nabla R^{+}+\rho^{+}\nabla R^{-}]-\sigma^{-}R^{-}\nabla \Delta R^{-}
\\&\qquad=\textrm{div}\big(\alpha^{-}[\mu^{-}(\nabla u^{-}+\nabla^{t}u^{-})+\lambda^{-}\textrm{div}u^{-}\textrm{Id}]\big).
\end{aligned}
\right.
\end{align}
Here, we are concerned with the Cauchy problem of  system
 \eqref{equ:CTFS1} in $\mathbb{R}_{+}\times \mathbb{R}^N$ subject to the
initial data
\begin{equation}\label{eq:initial data}
(R^{+},\,u^{+},\,R^{-},\,u^{-})(x,t)|_{t=0}=(R^{+}_{0},\,u^{+}_{0},\,R^{-}_{0},\,u^{-}_{0})(x), \quad x\in\mathbb{R}^{N},
\end{equation}
and
$$u^{+}(x,t)\rightarrow0,\quad u^{-}(x,t)\rightarrow0,\quad R^{+}\rightarrow R^{+}_{\infty}>0,\quad R^{-}\rightarrow R^{-}_{\infty}>0, \hbox{~as~} |x|\rightarrow\infty,$$
where $R^{\pm}_{\infty}$ denote the background doping profile, and in the present paper $R^{\pm}_{\infty}$ are taken as $1$ without losing generality.

 For simplicity, we take $\sigma^{+}=\sigma^{-}=1$. Set
$c^{\pm}=R^{\pm}-1$.
Then,  system \eqref{equ:CTFS1} can be rewritten  as
\begin{equation}\label{equ:CTFS2}
\left\{
\begin{aligned}{}
&\p_tc^{+}+\textrm{div}u^{+}=H_{1},\\
&\p_t{u}^{+}+\beta_{1}\nabla c^{+}
+\beta_{2}\nabla c^{-}-\nu_{1}^{+}\Delta u^{+}
-\nu_{2}^{+}\nabla\textrm{div}u^{+}-\nabla \Delta c^{+}=H_{2},
\\&\p_tc^{-}+\textrm{div}u^{-}=H_{3},
\\&\p_tu^{-}+\beta_{3}\nabla c^{+}
+\beta_{4}\nabla c^{-}-\nu_{1}^{-}\Delta u^{-}
-\nu_{2}^{-}\nabla\textrm{div}u^{-}-\nabla \Delta c^{-}=H_{4},
\end{aligned}
\right.
\end{equation}
with initial data
\begin{equation}\label{equ:CTFS3}(c^{+},\,u^{+},\,c^{-},\,u^{-})(x,t)|_{t=0}=(c^{+}_{0},\,u^{+}_{0},\,c^{-}_{0},\,u^{-}_{0})(x),\end{equation}
where $\beta_{1}=\frac{\mathcal{C}^{2}(1,1)\rho^{-}(1,1)}{\rho^{+}(1,1)},\quad
\beta_{2}=\beta_{3}=\mathcal{C}^{2}(1,1),\quad
\beta_{4}=\frac{\mathcal{C}^{2}(1,1)\rho^{+}(1,1)}{\rho^{-}(1,1)},\quad
\nu_{1}^{\pm}=\frac{\mu^{\pm}}{\rho^{\pm}(1,1)},\quad
\nu_{2}^{\pm}=\frac{\mu^{\pm}+\lambda^{\pm}}{\rho^{\pm}(1,1)}$
and the source terms are
\begin{align}\label{3.2}
H_{1}&=H_{1}(c^{+},u^{+})=-\textrm{div}(c^{+}u^{+}),\\
\label{3.3}
H_{2}^{i}&=H_{2}(c^{+},u^{+},c^{-})=-g_{+}(c^{+},c^{-})\partial_{i}c^{+}
-\tilde{g}_{+}(c^{+},c^{-})\partial_{i}c^{-}-(u^{+}\cdot\nabla)u_{i}^{+}\nonumber\\
&\quad+\mu^{+}h_{+}(c^{+},c^{-})\partial_{j}c^{+}\partial_{j}u^{+}_{i}
+\mu^{+}k_{+}(c^{+},c^{-})\partial_{j}c^{-}\partial_{j}u^{+}_{i}\nonumber\\
&\quad+\mu^{+}h_{+}(c^{+},c^{-})\partial_{j}c^{+}\partial_{i}u^{+}_{j}
+\mu^{+}k_{+}(c^{+},c^{-})\partial_{j}c^{-}\partial_{i}u^{+}_{j}\\
&\quad+\lambda^{+}h_{+}(c^{+},c^{-})\partial_{i}c^{+}\partial_{j}u^{+}_{j}
+\lambda^{+}k_{+}(c^{+},c^{-})\partial_{i}c^{-}\partial_{j}u^{+}_{j}\nonumber\\
&\quad+\mu^{+}l_{+}(c^{+},c^{-})\partial_{j}^{2}u_{i}^{+}+(\mu^{+}+\lambda^{+})l_{+}
(c^{+},c^{-})\partial_{i}\partial_{j}u^{+}_{j},\quad  i,j \in \{1,2,\cdots N\},\nonumber\\
\label{3.4}
H_{3}&=H_{3}(c^{-},u^{-})=-\textrm{div}(c^{-}u^{-}),\\
\label{3.5}
H_{4}^{i}&=H_{4}(c^{+},u^{-},c^{-})=-g_{-}(c^{+},c^{-})\partial_{i}c^{-}
-\tilde{g}_{-}(c^{+},c^{-})\partial_{i}c^{+}-(u^{-}\cdot\nabla)u_{i}^{-}\nonumber\\
&\quad+\mu^{-}h_{-}(c^{+},c^{-})\partial_{j}c^{+}\partial_{j}u^{-}_{i}
+\mu^{-}k_{-}(c^{+},c^{-})\partial_{j}c^{-}\partial_{j}u^{-}_{i}\nonumber\\
&\quad+\mu^{-}h_{-}(c^{+},c^{-})\partial_{j}c^{+}\partial_{i}u^{-}_{j}
+\mu^{-}k_{-}(c^{+},c^{-})\partial_{j}c^{-}\partial_{i}u^{-}_{j}\\
&\quad+\lambda^{-}h_{-}(c^{+},c^{-})\partial_{i}c^{+}\partial_{j}u^{-}_{j}
+\lambda^{-}k_{-}(c^{+},c^{-})\partial_{i}c^{-}\partial_{j}u^{-}_{j}\nonumber\\
&\quad+\mu^{-}l_{-}(c^{+},c^{-})\partial_{j}^{2}u_{i}^{-}+(\mu^{-}+\lambda^{-})l_{-}
(c^{+},c^{-})\partial_{i}\partial_{j}u^{-}_{j}, \quad  i,j \in \{1,2,\cdots N\},\nonumber
\end{align}
where we define the nonlinear functions of $(c^{+},c^{-})$ by
\begin{align}\label{3.6}
\left\{
\begin{aligned}
&g_{+}(c^{+},c^{-})=\frac{(\mathcal{C}^{2}\rho^{-})(c^{+}+1,c^{-}+1)}
{\rho^{+}(c^{+}+1,c^{-}+1)}
-\frac{(\mathcal{C}^{2}\rho^{-})(1,1)}
{\rho^{+}(1,1)},
\\&g_{-}(c^{+},c^{-})=\frac{(\mathcal{C}^{2}\rho^{+})(c^{+}+1,c^{-}+1)}
{\rho^{-}(c^{+}+1,c^{-}+1)}
-\frac{(\mathcal{C}^{2}\rho^{+})(1,1)}
{\rho^{-}(1,1)},
\end{aligned}
\right.
\end{align}
\begin{align}\label{3.7}
\left\{
\begin{aligned}
&h_{+}(c^{+},c^{-})=\frac{(\mathcal{C}^{2}\alpha^{-})(c^{+}+1,c^{-}+1)}
{[c^{+}+1]s_{-}^{2}(c^{+}+1,c^{-}+1)},
\\&h_{-}(c^{+},c^{-})=-\frac{\mathcal{C}^{2}(c^{+}+1,c^{-}+1)}
{(\rho^{-}s_{-}^{2})(c^{+}+1,c^{-}+1)},
\end{aligned}
\right.
\end{align}
\begin{align}\label{3.8}
\left\{
\begin{aligned}
&k_{+}(c^{+},c^{-})=-{\frac{\mathcal{C}^{2}(c^{+}+1,c^{-}+1)}
{[c^{-}+1](s_{+}^{2}\rho^{+})(c^{+}+1,c^{-}+1)}},
\\&k_{-}(c^{+},c^{-})={\frac{(\alpha^{+}\mathcal{C}^{2})(c^{+}+1,c^{-}+1)}
{[c^{-}+1]s_{+}^{2}(c^{+}+1,c^{-}+1)}},
\end{aligned}
\right.
\end{align}
\begin{align}\label{3.9}
\tilde{g}_{+}(c^{+},c^{-})=\tilde{g}_{-}(c^{+},c^{-})=\mathcal{C}^{2}(c^{+}+1,c^{-}+1)
-\mathcal{C}^{2}(1,1),
\end{align}
\begin{align}\label{3.10}
l_{\pm}(c^{+},c^{-})=\frac{1}{\rho^{\pm}(c^{+}+1,c^{-}+1)}-\frac{1}
{\rho^{\pm}(1,1)}.
\end{align}
Our first main result on the global well-posedness then reads  as follows.
\begin{Theorem}\label{th:main1}\,  Let $N\geq2$ and $p$ satisfy
$2\leq p\leq\min(4,\frac{2N}{N-2})$ and, additionally, $ p\neq 4$ if $ N=2$, and denote $c^{+}_0:=R^{+}_{0}-1$, $c^{-}_0:=R^{-}_{0}-1$.
There exists  a small enough  constant $\eta$ such that
if    $(c^{+h}_0,\,c^{-h}_0)\in \dot B^{\frac d{p}}_{p,1}\times\dot B^{\frac d{p}}_{p,1}$
and $(u^{+h}_0,\,u^{-h}_0)\in \dot B^{\frac d{p}-1}_{p,1}\times\dot B^{\frac d{p}-1}_{p,1}$ with
    besides  $(c_0^{+\ell},\,c_0^{-\ell})\in \dot B^{\frac d2-1}_{2,1}\times\dot B^{\frac d2-1}_{2,1}$ and $(u_0^{+\ell},\,u_0^{-\ell})\in \dot B^{\frac d2-1}_{2,1}\times\dot B^{\frac d2-1}_{2,1}$
   satisfy
\begin{equation}\label{eq:globalsmall1}
X_{0}\eqdefa \|(c_0^{+},\,u_0^{+},\,c_0^{-},\,u_0^{-})\|^\ell_{\dot B^{\frac N2-1}_{2,1}}+\|(c_0^{+},\,c_0^{-})\|^h_{\dot B^{\frac Np}_{p,1}}
+\|(u_0^{+},\,u_0^{-})\|^h_{\dot B^{\frac N{p}-1}_{p,1}}\leq \eta,\end{equation}
then the Cauchy  problem \eqref{equ:CTFS1}-\eqref{eq:initial data} admits a unique global-in-time solution $(c^{+},\,u^{+},\,c^{-},\,u^{-})$ in the space $\mathbb{X}^{\frac{N}{p}}$ defined by
\begin{equation*}\label{AA}
\begin{aligned}
&(c^{+},\,u^{+},\,c^{-},\,u^{-})^{\ell}\in \mathcal{C}(\mathbb{R}_+;\dot{B}_{2,1}^{\frac{N}{2}-1})\cap L^{1}(\mathbb{R}_+;\dot{B}_{2,1}^{\frac{N}{2}+1}),\\
&(c^{+},\,c^{-})^{h}\in \mathcal{C}(\mathbb{R}_+;\dot{B}_{p,1}^{\frac{N}{p}})\cap L^{1}(\mathbb{R}_+;\dot{B}_{p,1}^{\frac{N}{p}+2}),\\
&(u^{+},\,u^{-})^{h}\in \mathcal{C}(\mathbb{R}_+;\dot{B}_{p,1}^{\frac{N}{p}-1})\cap L^{1}(\mathbb{R}_+;\dot{B}_{p,1}^{\frac{N}{p}+1}).
\end{aligned}
\end{equation*}
Furthermore, we have for  $t\geq 0$,
\begin{equation}\label{1.6A}
\begin{split}
\,X(t)\lesssim X(0),
\end{split}
\end{equation}
where
\begin{equation}
\begin{split}
\label{1.6}
\,X(t)&\eqdefa\|(c^{+},\,u^{+},\,c^{-},\,u^{-})\|_{\wt L^\infty_t(\dot B^{\frac{N}{2}-1}_{2,1})}^\ell
+\|(c^{+},\,u^{+},\,c^{-},\,u^{-})\|_{L^1_t(\dot B^{\frac{N}{2}+1}_{2,1})}^\ell\\
&\quad+\|(u^{+},\,u^{-})\|_{\wt L^\infty_t(\dot B^{\frac{N}{p}-1}_{p,1})}^h
+\|(c^{+},\,c^{-})\|_{\wt L^\infty_t(\dot B^{\frac{N}{p}}_{p,1})}^h
\\
&\quad+\|(u^{+},\,u^{-})\|_{L^1_t(\dot B^{\frac{N}{p}+1}_{p,1})}^h
+\|(c^{+},\,c^{-})\|_{L^1_t(\dot B^{\frac{N}{p}+2}_{p,1})}^h.
\end{split}
\end{equation}
\end{Theorem}
Our second main result on the optimal time decay rates of strong solutions
states as follows.
\begin{Theorem}\label{th:decay}
Let the data $(c^{+}_{0},\,u^{+}_{0},\,c^{-}_{0},\,u^{-}_{0})$ satisfy the assumptions of Theorem \ref{th:main1}. Let $N\geq2$ and $p$ satisfy Theorem \ref{th:main1}. Denote $\langle t\rangle\eqdefa (1+t)$
and $\alpha\eqdefa\frac{N}{p}+\frac{1}{2}-\varepsilon$ with $\varepsilon>0$ arbitrarily small.
 There exists a positive constant $c$ such that if in addition
\begin{equation}\label{eq:D0}
D_0\eqdefa\|(c^{+}_{0},\,u^{+}_{0},\,c^{-}_{0},\,u^{-}_{0})\|^{\ell}_{\dot B^{-s_0}_{2,\infty}}\leq c \quad \text{with}\,\, s_0\eqdefa\frac{2N}{p}-\frac{N}{2}
\end{equation}
then  the global solution $(c^{+},\,u^{+},\,c^{-},\,u^{-})$
given by Theorem \ref{th:main1} satisfies for all $t\geq0,$
\begin{equation}
\label{1.8}
D(t)\lesssim\Big(D_0+\|(\nabla c^{+}_{0},\,u^{+}_{0},\,\nabla c^{-}_{0},\,u^{-}_{0})\|^h_{\dot B^{\frac Np-1}_{p,1}}\Big)
\end{equation}
with
\begin{equation}
\begin{split}
\label{1.9}
D(t)&\eqdefa \sup_{s\in[\varepsilon-s_0,\frac{N}{2}+1]}\|\langle\tau\rangle^{\frac{s+s_0}{2}}(c^{+},\,u^{+},\,c^{-},\,u^{-})\|_{L^\infty_t(\dot B^s_{2,1})}^\ell
\\ &\quad+\|\tau^{\alpha}(\nabla c^{+},\,u^{+},\,\nabla c^{-},\,u^{-})\|_{\wt L^\infty_t(\dot B^{\frac Np+1}_{p,1})}^h.
\end{split}
\end{equation}
\end{Theorem}
We would like to give some comments on our main results.
\begin{Remark}In Theorem \ref{th:main1},  the regularity indices for
the high frequency parts of $u^{\pm}_0$ may  be negative. Especially, this
allows us to obtain the global well-posedness of  system \eqref{equ:CTFS} for
the highly oscillating initial velocities $u^{\pm}_0$. For example, let
$$
u^{+}_0(x)=\sin\bigl(\f {x_1} {\varepsilon_1}\bigr)\phi(x),\quad u^{-}_0(x)=\sin\bigl(\f {x_1} {\varepsilon_2}\bigr)\phi(x),\quad \phi(x)\in \cS(\mathbb{R}^N).
$$
Thus for any $\varepsilon_i>0(i=1,2)$
$$
\|u^{+}_0\|_{\dot{ B}^{\f Np-1}_{p,1}}^{h}\leq C\varepsilon_{1}^{1-\frac{N}{p}}, \quad \|u^{-}_0\|_{\dot{ B}^{\f Np-1}_{p,1}}^{h}\leq C\varepsilon_{2}^{1-\frac{N}{p}} \quad \textrm{for}\quad p>N.
$$
Hence such data with small enough~$\varepsilon_i>0(i=1,2)$  generate global unique solutions
in dimensions $N=2,3$.
\end{Remark}
\begin{Remark}\label{1.2} Compared with \cite{CWYZ,LWWZ,Tan1,Tan2},   in Theorems \ref{th:main1} and \ref{th:decay}, we obtain the global well-posedness and     optimal time decay rates for multi-dimensional  non-conservative viscous
compressible two-fluid system
 \eqref{equ:CTFS} in critical $L^{p}$-framework respectively.   Additionally, in Theorem \ref{th:decay},
the regularity  index $s$ can take both negative and nonnegative values, rather than only nonnegative integers, which   improves the classical decay results  in high Sobolev regularity, such as \cite{CWYZ,LWWZ,Tan1,Tan2}. Moreover,  our results cover the case $N=2$.
 \end{Remark}
 \begin{Remark}\label{1.2-2} Compared with \cite{CWYZ}, in this paper,  we consider the case of general constant viscosities and   relax the special choice for viscosities in  \cite{CWYZ}.
 \end{Remark}

As a consequence of Theorem \ref{th:decay},  we can show the  following decay  rates of  $L^{p}$ norm of solution $(R^{+}-1,\,u^{+},\,R^{-}-1,\,u^{-})$.
\begin{Corollary}\label{cor1.1}  The solution $(R^{+}-1,\,u^{+},\,R^{-}-1,\,u^{-})$  constructed in  Theorem \ref{th:main1} satisfies
\begin{equation}\label{6.1}
\big\|\Lambda^{s}\big(R^{+}-1,R^{-}-1\big)\big\|_{L^p}
\lesssim \Big(D_0+\big\|\big(\nabla R^{+}_{0},\,u^{+}_{0},\,\nabla R^{-}_{0},\,u^{-}_{0}\big)\big\|^h_{\dot B^{\frac Np-1}_{p,1}}\Big)\langle t\rangle^{-\frac{s+s_0}{2}}
\quad  \hbox{ if }\, -s_0<s\leq N/p,
\end{equation}
\begin{equation}\label{6.2}
\big\|\Lambda^{s}\big(u^{+},u^{-}\big)\big\|_{L^p}
\lesssim \Big(D_0+\big\|\big(\nabla R^{+}_{0},\,u^{+}_{0},\,\nabla R^{-}_{0},\,u^{-}_{0}\big)\big\|^h_{\dot B^{\frac Np-1}_{p,1}}\Big)\langle t\rangle^{-\frac{s+s_0}{2}}
\quad \hbox{ if }\, -s_0<s\leq N/p-1,
\end{equation}
where the fractional derivative
 operator $\Lambda^{\ell}$ is defined by $\Lambda^{\ell}f\triangleq\mathcal{F}^{-1}(|\cdot|^{\ell}\mathcal{F}f)$.
\end{Corollary}
\begin{Remark}\label{1.2} In Corollary \ref{cor1.1},  taking $p=2$(hence $s_0=\frac{N}{2}$), $s=0$ leads back to  the standard optimal   $L^{1}$-$L^{2}$ time decay rates  which is a consistent with the optimal time decay rates  from  a single phase flow model in \cite{Tan1,Tan2}.
\end{Remark}

Before going into the heart of the proof of our main results, 
we make a brief interpretation of the main difficulties and techniques involved in the proof.  Due to the mixed
hyperbolic-parabolic property of system \eqref{equ:CTFS2},   
the system has to be handled differently in the low
and high frequencies respectively. Roughly speaking,  the first order terms
predominate in low frequencies, so that  system \eqref{equ:CTFS2} has to
be treated by means of hyperbolic energy methods, which implies that
we must treat the low frequencies regime only in spaces constructed
on $L^{2}$, as it is classical that hyperbolic systems are ill-posed
in general $L^{p}$ spaces.  In contrast, in the  high frequencies, a
$L^{ p}$ approach may be used. On the other hand, various important mathematical difficulties occur when we want to generalize well-known results of the compressible Navier-Stokes equations  to the two-phase system \eqref{equ:CTFS2} since the corresponding model is
non-conservative.

First, in the low frequencies,  we need deal with system \eqref{equ:CTFS2} in $L^{2}$-framework by means of hyperbolic energy methods. In general, the proof consists of spectral analysis of  Green's function for the corresponding linearized
system and energy estimates of the solutions to the nonlinear system, refer for
instance to \cite{CD,CMZ2,HS2}. However, we  encounter a fundamental obstacle that
 Green's function of the viscous compressible two-fluid
model \eqref{equ:CTFS2} is an 8-order matrix and is not self-adjoint so that we can not make some complicate analysis of Green's function. To get around this difficulty, we will follow
the main ideas of 
\cite{XC} (see  Lemma 3.1 for details) to exploit the maximal regularity estimates for  the corresponding linearized  system in   $L^{2}$-type critical Besov spaces by employing the energy argument of Godunov
\cite{GO} for partially dissipative first-order symmetric systems (further developed by
 \cite{KF}) and  Fourier-Plancherel
theorem.

Second,  we need to cover more general values of the integration parameter $p$ in the high frequencies.
Generally speaking, in order to solve this problem,  there are two fundamental effective methods concerning the standard barotropic compressible  Navier-Stokes equations. The first  is some suitable effective
velocity field (named viscous effective flux in Hoff's work \cite{DH}) in order to kill the relation of coupling between the velocity and the pressure from \cite{HS2}. The second is  the study of the paralinearized system combined with a Lagrangian change of coordinates (in the spirit of that introduced by Hmidi in
\cite{Hi} for the convection-diffusion equation) from \cite{CD,CMZ2} so as to counter the loss of regularity coming from the convection terms and
almost completely avoid the undesired coupling.  Unfortunately,   these two  methods fail to work for the current  system  \eqref{equ:CTFS2} due to the non-conservation. To explain  this problem,  one study the following non-conservative coupled hyperbolic-parabolic model
\begin{equation*}
\left\{
\begin{aligned}{}
&\p_tc^{+}+\textrm{div}\cQ u^{+}=0,\\
&\p_t\cQ u^{+}+\beta_{1}\nabla c^{+}
+\beta_{2}\nabla c^{-}-\nu^{+}\Delta\cQ u^{+}-\nabla \Delta c^{+}=0,
\\&\p_tc^{-}+\textrm{div}\cQ u^{-}=0,
\\&\p_t\cQ u^{-}+\beta_{3}\nabla c^{+}
+\beta_{4}\nabla c^{-}-\nu^{-}\Delta \cQ u^{-}-\nabla \Delta c^{-}=0,
\end{aligned}
\right.
\end{equation*}
where $\cQ$ is  potential vector fields. 
For simplicity, here we only consider the above non-conservative coupled hyperbolic-parabolic without capillarity effects. Following from Haspot's idea in \cite{HS2}, we introduce two new auxiliary functions
$W_{1}\eqdefa\mathcal{Q}{u}^{+}
+\frac{\gamma_{1}}{\nu^{+}}(-\Delta)^{-1}\nabla c^{+}
$ and $W_{2}\eqdefa\mathcal{Q}{u}^{-}+\frac{\gamma_{4}}{\nu^{-}}(-\Delta)^{-1}\nabla c^{-}$,
then we get
\begin{equation*}\label{w}
\left\{
\begin{aligned}{}
&\p_tW_{1}-\nu^{+}\Delta W_{1}
=\frac{\gamma_{1}^{2}}{\nu^{+}}W_{1}
-\frac{\gamma_{1}^{3}}{(\nu^{+})^{2}}
(-\Delta)^{-1}\nabla c^{+}-\gamma_{2}\nabla c^{-},
\\&\p_tW_{2}-\nu^{-}\Delta W_{2}
=\frac{\gamma_{4}^{2}}{\nu^{-}}W_{2}-\frac{\gamma_{4}^{3}}{(\nu^{-})^{2}}
(-\Delta)^{-1}\nabla c^{-}-\gamma_{3}\nabla c^{+},
\end{aligned}
\right.
\end{equation*}
and 
\begin{equation*}\label{trnsport}
\left\{
\begin{aligned}{}
&\partial_{t}c^{+}+\frac{\gamma_{1}^{2}}{\nu^{+}}c^{+}
+\frac{1}{\sqrt{\beta_1}}u^{+}\cdot\nabla c^{+}
=-\frac{1}{\sqrt{\beta_1}}c^{+}\textrm{div}u^{+}-\gamma_1\textrm{div}W_{1},
\\&\partial_{t}c^{-}+\frac{\gamma_{4}^{2}}{\nu^{-}}c^{-}
+\frac{1}{\sqrt{\beta_4}}u^{-}\cdot\nabla c^{-}
=-\frac{1}{\sqrt{\beta_4}}c^{-}\textrm{div}u^{-}-\gamma_4\textrm{div}W_{2}.
\end{aligned}
\right.
\end{equation*}
By simple calculation, we have  
 \begin{equation*}
\begin{split}\label{eq:W1}
&\|W_{1}\|_{\wt L^\infty_t(\dot B^{\frac Np-1}_{p,1})}^h+\nu^{+}\|W_{1}\|_{L^1_t(\dot B^{\frac Np+1}_{p,1})}^h
 \\& \leq C_1\Big(\|W_{10}\|_{\dot B^{\frac Np-1}_{p,1}}^h +
 \frac{\gamma_{1}^{2}}{\nu^{+}}\|W_1\|_{L^1_t(\dot B^{\frac Np-1}_{p,1})}^h+\frac{\gamma_{1}^{3}}{(\nu^{+})^{2}}\| c^{+}\|_{L^1_t(\dot B^{\frac Np-2}_{p,1})}^h+\gamma_{2}\| c^{-}\|_{L^1_t(\dot B^{\frac Np}_{p,1})}^h\Big),
\end{split}
\end{equation*}
\begin{equation*}
\begin{split}\label{eq:W2}
 &\|W_{2}\|_{\wt L^\infty_t(\dot B^{\frac Np-1}_{p,1})}^h+\nu^{+}\|W_{2}\|_{L^1_t(\dot B^{\frac Np+1}_{p,1})}^h
 \\&\leq C_2\Big(\|W_{20}\|_{\dot B^{\frac Np-1}_{p,1}}^h 
 +\frac{\gamma_{4}^{2}}{\nu^{-}}\|W_2\|_{L^1_t(\dot B^{\frac Np-1}_{p,1})}^h+\frac{\gamma_{4}^{3}}{(\nu^{-})^{2}}\| c^{-}\|_{L^1_t(\dot B^{\frac Np-2}_{p,1})}^h+\gamma_{3}\| c^{+}\|_{L^1_t(\dot B^{\frac Np}_{p,1})}^h\Big),
\end{split}
\end{equation*}
\begin{equation*}\label{transport3}
\begin{split}
\|c^{+}\|_{\wt L^\infty_t(\dot B^{\frac{N}{p}}_{p,1})}^{h}
&+\frac{\gamma_{1}^{2}}{\nu^{+}}\|c^{+}\|_{L^{1}_{t}(\dot B^{\frac{N}{p}}_{p,1})}^{h}
\\&\leq C_3\Big(\|c^{+}_{0}\|_{\dot B^{\frac{N}{p}}_{p,1}}^{h}
+\gamma_1\|\textrm{div}W_{1}\|_{L^{1}_{t}(\dot B^{\frac{N}{p}}_{p,1})}^{h}
+\frac{1}{\sqrt{\beta_1}}\int_{0}^{t}\|\nabla u^{+}\|_{\dot B^{\frac{N}{p}}_{p,1}}
\|c^{+}\|_{\dot B^{\frac{N}{p}}_{p,1}}d\tau\Big),
\end{split}
\end{equation*}
and
\begin{equation*}\label{transport4}
\begin{split}
\|c^{-}\|_{\wt L^\infty_t(\dot B^{\frac{N}{p}}_{p,1})}^{h}
&+\frac{\gamma_{4}^{2}}{\nu^{-}}\|c^{-}\|_{L^{1}_{t}(\dot B^{\frac{N}{p}}_{p,1})}^{h}
\\&\leq C_4\Big(\|c^{-}_{0}\|_{\dot B^{\frac{N}{p}}_{p,1}}^{h}
+\gamma_4\|\textrm{div}W_{2}\|_{L^{1}_{t}(\dot B^{\frac{N}{p}}_{p,1})}^{h}+\frac{1}{\sqrt{\beta_4}}\int_{0}^{t}\|\nabla u^{-}\|_{\dot B^{\frac{N}{p}}_{p,1}}
\|c^{-}\|_{\dot B^{\frac{N}{p}}_{p,1}}d\tau\Big).
\end{split}
\end{equation*}
Different from  the barotropic compressible Navier-Stokes equations \cite{HS2} and  the compressible Navier-Stokes
system with capillarity \cite{CDX,KSX}, there are
four  terms $\|\textrm{div}W_{1}\|_{L^{1}_{t}(\dot B^{\frac{N}{p}}_{p,1})}^{h}, \|\textrm{div}W_{2}\|_{L^{1}_{t}(\dot B^{\frac{N}{p}}_{p,1})}^{h},$ $\| c^{+}\|_{L^1_t(\dot B^{\frac Np}_{p,1})}^h$  and $\| c^{-}\|_{L^1_t(\dot B^{\frac Np}_{p,1})}^h$  appearing  in the righthand sides of the above estimates.
Obviously, it  seems impossible  to obtain the  desired  estimates of $(c_1,W_1,c_2,W_2)$ when  the four  terms are  simultaneously treated  as  source terms  in the high frequencies. Moreover, the method from \cite{CD,CMZ2} has also similar difficulty  for system  \eqref{equ:CTFS2}. To overcome this essential difficulty, 
the main idea here is that in view of a crucial observation according to the  nice mathematical structures of   system \eqref{equ:CTFS2}. 
More precisely,  applying operator $\Delta$  to the above equations of  $c^{\pm}$ and taking operator $\Dv$  to the  equations of  $\cQ u^{\pm}$, we immediately notice that $\dv \cQ u^{\pm}$  should have the same regularity as  $\Delta c^{\pm}$.
 Moreover, we also 
find that  $(\Delta c^{\pm},\dv \cQ u^{\pm})$ satisfy two similar linearized coupling  systems which likely  exhibit parabolic  properties in the high frequencies due to the presence  of the capillary terms by the spectral analysis (see \eqref{gama3-A-2-1} and Lemma \ref{lemma3.4} for details).  More importantly, based on the observation,   the perturbations  $\beta_{1}\Delta c^{+},\beta_{2}\Delta c^{-},\beta_{3}\Delta c^{+}$ and $\beta_{4}\Delta c^{-}$ can be treated
as  harmless source terms  in the high frequencies, which induces us to get  the desired optimal \emph{a priori} estimates of $(c^{+},\,u^{+},\,c^{-},\,u^{-})$ in critical $L^{p}$-framework.

Third, we also investigate  how global strong solutions constructed above look like for large time. In this part, our main ideas are based on  a refined time-weighted energy inequalities in the Fourier spaces and the benefit of low-frequency and high-frequency decomposition. In the low frequencies, making good use of
Fourier   localization  analysis to a linearized parabolic-hyperbolic system and Parseval's equality in order to obtain  smoothing effects of  Green's function  and  avoid some complicate   spectral analysis as in \cite{CWYZ}.  Consequently, it is possible to adapt the standard Duhamel's principle handling those nonlinear terms.  With the aid of the nonclassical product estimates in Besov spaces, one can obtain the desired estimates. In the high frequencies,  in order to close  the estimate from time-weighted energy
functional, using  the  similar  method of the proof for global existence in Section three together with elaborate nonlinear estimates, we further exploit some  decay estimates with gain of regularity   of $(\nabla c^{+},\,u^{+},\,\nabla c^{-},\,u^{-})$.

Finally,  we also have to deal  with some difficulties caused by much more complicate nonlinear terms by using the usual product estimates in  Besov spaces, Bony's decomposition and the low-high frequency decomposition, and  some new  composition of the  binary functions from harmonic analysis.


The rest of the paper unfolds as follows. In the next section,  we recall some basic facts about Littlewood-Paley
decomposition, Besov spaces  and some useful lemmas.
 Section 3 is devoted to the proof of the global well-posedness for initial data near equilibrium in critical Besov spaces.
In Section 4,  we present  the optimal time decay rates of the  global  strong solutions. 
Some material concerning
 paradifferential calculus and product  estimates in Besov spaces
is recalled in Appendix.

\noindent{\bf Notations.}  We assume $C$ be
a positive generic constant throughout this paper that may vary at
different places and
denote  $A\le CB$ by  $A\lesssim B$.
We shall also use the  following notations
 $$z^\ell\eqdefa\sum_{j\leq k_0}\dot{\Delta}_{j}z\quad\hbox{and}\quad z^h\eqdefa z-z^\ell, \quad\hbox{for some}   k_0.$$
$$\|z\|^\ell_{\dot B^s_{2,1}}\eqdefa \sum_{j\leq k_0}2^{js}\|\dot{\Delta}_{j}z\|_{L^{2}}\quad\hbox{and}\quad \|z\|^h_{\dot B^s_{2,1}}\eqdefa \sum_{j\geq k_0}2^{js}\|\dot{\Delta}_{j}z\|_{L^{2}}, \quad\hbox{for some }  k_0.$$
Noting the small overlap between low and high frequencies,  we have
$$\|z^\ell\|_{\dot B^s_{2,1}}\lesssim \|z\|^\ell_{\dot B^s_{2,1}}\quad\hbox{and}\quad \|z^h\|_{\dot B^s_{2,1}}\lesssim \|z\|^h_{\dot B^s_{2,1}}.$$
\par

\section{ Littlewood-Paley Theory and Some Useful Lemmas }
\ \ \ \ \ Let us introduce the Littlewood-Paley decomposition.
Choose a radial function  $\varphi \in {\cS}(\mathbb{R}^N)$
supported in ${\cC}=\{\xi\in\mathbb{R}^N,\,
\frac{3}{4}\le|\xi|\le\frac{8}{3}\}$ such that \beno \sum_{q\in
\mathbb{Z}}\varphi(2^{-q}\xi)=1 \quad \textrm{for all}\,\,\xi\neq 0.
\eeno The homogeneous frequency localization operators
$\dot{\Delta}_q$ and $\dot{S}_q$ are defined by
\begin{align}
\dot{\Delta}_qf=\varphi(2^{-q}D)f,\quad \dot{S}_qf=\sum_{k\le
q-1}\dot{\Delta}_kf\quad\mbox{for}\quad q\in \mathbb{Z}. \nonumber
\end{align}
With our choice of $\varphi$, one can easily verify that
\begin{equation*}\begin{split}
&\dot{\Delta}_q\dot{\Delta}_kf=0\quad \textrm{if}\quad|q-k|\ge
2\quad \textrm{and}
\quad \\
&\dot{\Delta}_q(\dot{S}_{k-1}f\dot{\Delta}_k f)=0\quad
\textrm{if}\quad|q-k|\ge 5.
\end{split}\end{equation*}
 We denote the space ${\cZ'}(\mathbb{R}^N)$ by the dual space of
${\cZ}(\mathbb{R}^N)=\{f\in {\cS}(\mathbb{R}^N);\,D^\alpha
\hat{f}(0)=0; \forall\alpha\in\mathbb{ N}^N \,\mbox{multi-index}\}.$
It also can be identified by the quotient space of
${\cS'}(\mathbb{R}^N)/{\cP}$ with the polynomials space ${\cP}$. The
formal equality \beno f=\sum_{q\in\mathbb{Z}}\dot{\Delta}_qf \eeno
holds true for $f\in {\cZ'}(\mathbb{R}^N)$ and is called the
homogeneous Littlewood-Paley decomposition.

Let us recall the definition of homogeneous
Besov spaces  (see \cite{BCD,Dan1}).
\begin{Definition}\label{def2.2} Let $s\in \mathbb{R}$, $1\le p,
r\le+\infty$. The homogeneous Besov space $\dot{B}^{s}_{p,r}$ is
defined by
$$\dot{B}^{s}_{p,r}=\Big\{f\in {\cZ'}(\mathbb{R}^N):\,\|f\|_{\dot{B}^{s}_{p,r}}<+\infty\Big\},$$
where \beno \|f\|_{\dot{B}^{s}_{p,r}}\eqdefa \Big\|2^{qs}
\|\dot{\Delta}_qf(t)\|_{p}\Big\|_{\ell^r}.\eeno
\end{Definition}
\begin{Remark}\label{2.3}\quad
Some properties about the  Besov spaces are as follows
\begin{itemize}
\item\,\, Derivation: $$\|f\|_{\dot{B}^{s}_{2,1}}\approx\|\nabla f\|_{\dot{B}^{s-1}_{2,1}};$$
\item\,\, Algebraic properties: for $s > 0$, $\dot{B}^{s}_{2,1} \cap L^{\infty}
$ is an algebra;
\item\,\,Interpolation: for
$s_1, s_2\in\mathbb{R}$ and $\theta\in[0,1]$,
we have $$\|f\|_{\dot{B}^{\theta s_1+(1-\theta)s_2}_{2,1}}\le
\|f\|^\theta_{\dot{B}^{s_1}_{2,1}}\|f\|^{1-\theta}_{\dot{B}^{s_2}_{2,1}}.$$
\end{itemize}
\end{Remark}
\begin{Definition} Let $s\in \mathbb{R}$,$(r,\rho,p)\in[1,\,+\infty]^3$ and $T\in(0,\,+\infty]$.
We say then that $f\in L^\rho_{T}(\dot{B}_{p,  r}^s),$ if
$$\|f\|_{L^\rho_{T}(\dot{B}_{p,  r}^s)}\eqdefa \Big\|\big\| 2^{qs}
\|\dot{\Delta}_q f\|_{L^p}\big\|_{\ell ^{r}}\Big\|_{L^\rho_{T}}<+\infty.$$
\end{Definition}

We next introduce the Besov-Chemin-Lerner space
$\widetilde{L}^q_T(\dot{B}^{s}_{p,r})$ which is initiated in
\cite{Che-Ler}.
\begin{Definition}Let $s\leq\frac{N}{p}$ (respectively $s\in\mathbb{R}$),
$(r,\rho,p)\in[1,\,+\infty]^3$ and $T\in(0,\,+\infty]$. We
define $\widetilde{L}^{\rho}_T(\dot B^s_{p\,r})$ as the completion of  $C([0, T]; \mathcal{S}_{h}')$  by the norm
$$
\|f\|_{\tilde{L}^\rho_{T}(\dot{B}_{p,  r}^s)}\eqdefa \Big\| 2^{js}
\|\dot{\Delta}_j f(t)\|_{L^\rho(0,T;L^{p})}\Big\|_{\ell^r}
<\infty,
$$
with the usual change if $r=\infty.$
\end{Definition}
Obviously, $
\widetilde{L}^1_T(\dot{B}^s_{p,1})=L^1_T(\dot{B}^s_{p,1}). $ By  a
direct application  of  Minkowski's inequality, we have the
following relations between these spaces
\begin{equation*}
L^\rho_{T}(\dot{B}_{p,r}^s)\hookrightarrow\widetilde
L^\rho_{T}(\dot{B}_{p,r}^s),\,\textnormal{if}\quad  r\geq
\rho,\end{equation*}
\begin{equation*}
\widetilde L^\rho_{T}(\dot{B}_{p,r}^s)\hookrightarrow
L^\rho_{T}(\dot{B}_{p,r}^s),\, \textnormal{if}\quad \rho\geq r.
\end{equation*}

The following Bernstein's inequalities will be frequently used.
\begin{Lemma}\cite{Che-book}\label{Lem:Bernstein}
Let $1\le p_{1}\le p_{2}\le+\infty$. Assume that $f\in L^{p_{1}}(\mathbb{R}^N)$,
then for any $\gamma\in(\mathbb{N}\cup\{0\})^N$, there exist
constants $C_1$, $C_2$ independent of $f$, $q$ such that \beno
&&{\rm supp}\hat f\subseteq \{|\xi|\le A_02^{q}\}\Rightarrow
\|\partial^\gamma f\|_{p_2}\le C_12^{q{|\gamma|}+q
N(\frac{1}{p_1}-\frac{1}{p_2})}\|f\|_{p_1},
\\
&&{\rm supp}\hat f\subseteq \{A_12^{q}\le|\xi|\le
A_22^{q}\}\Rightarrow \|f\|_{p_1}\le
C_22^{-q|\gamma|}\sup_{|\beta|=|\gamma|}\|\partial^\beta f\|_{p_1}.
\eeno
\end{Lemma}

We here recall basic nonlinear estimates in  Besov spaces  which will be used repeatedly in our proof.
\begin{Proposition}\label{p26}\cite{RS}
 For all $1\leq r,p, p_1, p_2\leq+\infty$,  there exists a positive universal
 constant such that
$$\|fg\|_{\dot{B}^{s}_{p,r}}\lesssim
\|f\|_{L^\infty}\|g\|_{\dot{B}^{s}_{p,r}}+\|g\|_{L^\infty}\|f\|_{\dot{B}^{s}_{p,r}},
\quad \text{if}\quad s>0;$$
$$\|fg\|_{\dot{B}^{s_1+s_2-\frac{N}{p}}_{p,r}}\lesssim
\|f\|_{\dot{B}^{s_1}_{p,r}}\|g\|_{\dot{B}^{s_2}_{p,\infty}}, \quad
\text{if}\quad s_1,s_2<\frac{N}{p},\quad \text{and}\quad
s_1+s_2>0;$$
$$\|fg\|_{\dot{B}^{s}_{p,r}}\lesssim
\|f\|_{\dot{B}^{s}_{p,r}}\|g\|_{\dot{B}^{\frac{N}{p}}_{p,\infty}\cap
L^{\infty}}, \quad \text{if}\quad |s|<\frac{N}{p};$$
$$\|fg\|_{\dot{B}^s_{2,1}}\lesssim \|f\|_{\dot{B}^{N/2}_{2,1}}\|g\|_{\dot{B}^s_{2,1}}, \quad \text{if}\quad
s\in (-N/2,N/2].$$
\end{Proposition}

The basic tool of the paradifferential calculus is Bony's decomposition \cite{Bony}.  Formally, the  product  of two  tempered distributions $u$ and $v$   may be decomposed
into
\begin{equation*}
uv=\dot{T}_{u}v+\dot{T}_{v}u+\dot{R}(u,v)
\end{equation*}
with
\begin{equation*}
\dot{T}_{u}v=\sum\limits_{j\in\mathbb{Z}}\dot{S}_{j-1}u\dot{\Delta}_{j}v,\quad
\dot{R}(u,v)=\sum\limits_{j\in\mathbb{Z}}\dot{\Delta}_{j}u\widetilde{\dot{\Delta_{j}}}v,\quad
\widetilde{\dot{\Delta_{j}}}v=\sum\limits_{|j-j'|\leq1}\dot{\Delta}_{j'}v.
\end{equation*}
As a consequence, the estimates for the paraproduct and remainder operators can be given by
\begin{Proposition}\cite{BH}\label{200}\ \  Let $N\geq 2$, $s\in \mathbb{R}$ and $2\leq p\leq \min(4,\frac{2N}{N-2})$, we have
\begin{align}\label{A7.1}\|T_{f}g\|_{L_{t}^{1}(\dot{B}_{2,1}^{s-1+\frac{N}{2}-\frac{N}{p}})}\leq C\|f\|_{L_{t}^{\infty}(\dot{B}_{p,1}^{\frac{N}{p}-1})}
\|g\|_{L_{t}^{1}(\dot{B}_{p,1}^{s})}.
\end{align} In particular,
for $s\in \mathbb{R}$, $m\geq0$,   we also have
\begin{align}\label{A7.2}\|(T_{f}g)^\ell\|_{L_{t}^{1}(\dot{B}_{2,1}^{s-1+\frac{N}{2}-\frac{d}{p}})}\leq C\|f\|_{\dot{B}_{p,1}^{\frac{N}{p}-1}}
\|g\|_{\dot{B}_{p,1}^{s-m}}.\end{align}
\end{Proposition}
\begin{Proposition}\cite{BH}\label{Pro:201}\ \  Let $N\geq2$,\quad $s>1-\min(\frac{N}{p},\frac{N}{p'})$ and $1\leq p\leq4$, we have
$$\|R(f,g)\|_{\dot{B}_{2,1}^{s-1+\frac{N}{2}-\frac{N}{p}}}\leq C\|f\|_{\dot{B}_{p,1}^{\frac{N}{p}-1}}
\|g\|_{\dot{B}_{p,1}^{s}}.$$
\end{Proposition}
\begin{Proposition}\cite{DX}\label{5.5}\quad Let the real numbers $\sigma_{1}$, $\sigma_{2}$, $p_{1}$ and $p_{2}$  satisfy
$$\sigma_{1}+\sigma_{2}>0,\quad \sigma_{1}\leq\frac{N}{p_{1}},\quad \sigma_{2}\leq\frac{N}{p_{2}},\quad \sigma_{1}\geq\sigma_{2},\quad \frac{1}{p_{1}}+\frac{1}{p_{2}}\leq 1,$$
then
$$\|fg\|_{\dot{B}_{q,1}^{\sigma_{2}}}\lesssim\|f\|_{\dot{B}_{p_{1},1}^{\sigma_{1}}}\|g\|_{\dot{B}_{p_{2},1}^{\sigma_{2}}}\quad with\quad\frac{1}{q}=\frac{1}{p_{1}}+\frac{1}{p_{2}}-\frac{\sigma_{1}}{N}.$$
Let the exponents $\sigma>0$ and $1\leq p_{1},p_{2},q\leq\infty$  satisfy
$$\frac{N}{p_{1}}+\frac{N}{p_{2}}-d\leq\sigma\leq \min(\frac{N}{p_{1}},\frac{N}{p_{2}})\quad with\quad\frac{1}{q}=\frac{1}{p_{1}}+\frac{1}{p_{2}}-\frac{\sigma}{N},$$
then
$$\|fg\|_{\dot{B}_{q,\infty}^{-\sigma}}\lesssim\|f\|_{\dot{B}_{p_{1},1}^{\sigma}}
\|g\|_{\dot{B}_{p_{2},\infty}^{-\sigma}}.$$
\end{Proposition}
\begin{Corollary}\label{cor7.4} Let $N\geq2$ and $p$ satisfy
$2\leq p\leq\min(4,\frac{2N}{N-2})$ and, additionally, $ p\neq 4$ if $ N=2$, then
\begin{equation}\label{low.7}
\|fg\|_{\dot B^{-s_0}_{2,\infty}}
\lesssim\|f\|_{\dot B^{1-\frac{N}{p}}_{p,1}}
\|g\|_{\dot B^{\frac{N}{2}-1}_{2,1}}.
\end{equation}
\end{Corollary}
\begin{Corollary}\label{cor7.5}Let $N\geq2$ and $p$ satisfy
$2\leq p\leq N$, then
\begin{equation}\label{low.6.23}
\|fg^{h}\|_{\dot B^{-s_0}_{2,\infty}}^{\ell}
\lesssim\|f\|_{\dot B^{\frac{N}{p}-1}_{p,1}}
\|g^{h}\|_{\dot B^{\frac{N}{p}-1}_{p,1}}.
\end{equation}
\end{Corollary}
\begin{Proposition}\cite{DX}
Let $q_0\in\mathbb{Z},$ and denote $\dot{S}_{q_0}u\triangleq u^{\ell}$ and for any $s\in\mathbb{R},$
there exists a universal integer $N_0$ such that for any $2\leq p\leq 4$ and $\sigma>0,$  then
\begin{equation}\label{prolow1}
\|uv^{h}\|_{\dot B^{-s_0}_{2,\infty}}^{\ell}
\leq C(\|u\|_{\dot B^{\sigma}_{p,1}}+\|\dot{S}_{q_0+N_0}u\|_{L^{p^{*}}})
\|v^{h}\|_{\dot B^{-\sigma}_{p,\infty}}
\end{equation}
and
\begin{equation}\label{prolow2}
\|u^{h}v\|_{\dot B^{-s_0}_{2,\infty}}^{\ell}
\leq C(\|u^{h}\|_{\dot B^{\sigma}_{p,1}}+\|\dot{S}_{q_0+N_0}u^{h}\|_{L^{p^{*}}})
\|v\|_{\dot B^{-\sigma}_{p,\infty}}
\end{equation}
with $s_0:=\frac{2N}{p}-\frac{N}{2}$ and $\frac{1}{p^{*}}:=\frac{1}{2}-\frac{1}{p}$, and $C$ depending only on $k_0, N$ and $\sigma.$
\end{Proposition}

For the composition of the  binary functions, we have the following estimates.
\begin{Proposition}\label{p27}\cite{XC}
\mbox{(i) } Let $I$ be an open interval of $\mathbb{R}$ containing 0, $s>0$, $t\geq0$, $1\le p, q,r\le \infty$ and $(f,g)\in \big(\tilde{L}_{t}^{q}(\dot{B}^{s}_{p,r})\cap L_{t}^{\infty}(L^{\infty})\big)\times\big(\tilde{L}_{t}^{q}(\dot{B}^{s}_{p,r})\cap L_{t}^{\infty}(L^{\infty})\big)$. If $F\in W_{loc}^{[s]+2,\infty}(I)\times W_{loc}^{[s]+2,\infty}(I)$ with
$F(0,0)=0$, then $F(f,g)\in \tilde{L}_{t}^{q}(\dot{B}^{s}_{p,r})$. Moreover,  there exists a positive constant
 $C$ depending only on $s,p,N$ and $F$
such that
\begin{equation}\label{eq:2.1}
\|F(f,g)\|_{\tilde{L}_{t}^{q}(\dot{B}^{s}_{p,r})}
\leq C\Big(1+\|f\|_{L_{t}^{\infty}(L^{\infty})}\|g\|_{L_{t}^{\infty}(L^{\infty})}\Big)^{[s]+1}
\|(f,g)\|_{\tilde{L}_{t}^{q}(\dot{B}^{s}_{p,r})}.
\end{equation}
\mbox{ (ii)} If $(f_1,g_1)\in \tilde{L}_{t}^{\infty}(\dot{B}^{\frac{N}{p}}_{p,1})\times\tilde{L}_{t}^{\infty}(\dot{B}^{\frac{N}{p}}_{p,1})$ and $(f_2,g_2)\in \tilde{L}_{t}^{\infty}(\dot{B}^{\frac{N}{p}}_{p,1})\times\tilde{L}_{t}^{\infty}(\dot{B}^{\frac{N}{p}}_{p,1})$, and  $(f_2-f_1,g_2-g_1)$ belongs to $\tilde{L}_{t}^{q}(\dot{B}^{s}_{p,r})\times \tilde{L}_{t}^{q}(\dot{B}^{s}_{p,r})$ with $s\in (-\frac{N}{p},\frac{N}{p}]$ . If $F\in W_{loc}^{[\frac{N}{p}]+2,\infty}(\mathbb{R})\times W_{loc}^{[\frac{N}{p}]+2,\infty}(\mathbb{R})$ with
$\partial_{1}F(0,0)=0$ and $\partial_{2}F(0,0)=0$. Then, there exists a positive constant
 $C$ depending only on $s, p, N$ and $F$
such that
\begin{equation}\label{eq:2.3}\begin{split}
&\|F(f_2,g_2)-F(f_1,g_1)\|_{\tilde{L}_{t}^{q}(\dot{B}^{s}_{p,r})}
\\&\leq C\Big(1+\|(f_1,f_2)\|_{L_{t}^{\infty}(L^{\infty})}\|(g_1,g_2)\|_{L_{t}^{\infty}(L^{\infty})}\Big)^{[\frac{N}{p}]+1} \Big(\|(f_1,g_1)\|_{\tilde{L}_{t}^{\infty}(\dot{B}^{\frac{N}{p}}_{p,1})} +\|(f_2,g_2)\|_{\tilde{L}_{t}^{\infty}(\dot{B}^{\frac{N}{p}}_{p,1})}\Big)
\\&\quad\times\Big(\|f_2-f_1\|_{\tilde{L}_{t}^{q}(\dot{B}^{s}_{p,r})}+\|g_2-g_1\|_{\tilde{L}_{t}^{q}(\dot{B}^{s}_{p,r})}\Big).
\end{split}
\end{equation}
\end{Proposition}
In  $L^{p}$-framework, we need deal with  the composition of the  binary functions the case $s<0$.
\begin{Proposition}\label{p28}
Let $I$ be an open interval of $\mathbb{R}$ containing $0$, $s>-\min\{\frac{N}{p},\frac{N}{p'}\}$, $t\geq0$, $1\le p, q,r\le \infty$ and $(f,g)\in \big(\tilde{L}_{t}^{q}(\dot{B}^{s}_{p,r})\cap L_{t}^{\infty}(\dot{B}^{\frac{N}{p}}_{p,1})\big)\times\big(\tilde{L}_{t}^{q}(\dot{B}^{s}_{p,r})\cap L_{t}^{\infty}(\dot{B}^{\frac{N}{p}}_{p,1})\big)$. If $F\in W_{loc}^{[s]+2,\infty}(I)\times W_{loc}^{[s]+2,\infty}(I)$ with
$F(0,0)=0$, then $F(f,g)\in \tilde{L}_{t}^{q}(\dot{B}^{s}_{p,r})$. Moreover,  there exists a positive constant
 $C$ depending only on $s,p,N$ and $F$
such that
\begin{equation}\label{eq:2.1-A}\begin{split}
\|F(f,g)\|_{\tilde{L}_{t}^{q}(\dot{B}^{s}_{p,r})}
&\leq C\Big(|F'_{f}(0,0)+|F'_{g}(0,0)|\Big)\|(f,g)\|_{\tilde{L}_{t}^{q}(\dot{B}^{s}_{p,r})}
\\&\quad+\Big(1+\|f\|_{L_{t}^{\infty}(\dot{B}^{\frac{N}{p}}_{p,1})}\|g\|_{L_{t}^{\infty}(\dot{B}^{\frac{N}{p}}_{p,1})}\Big)^{[s]+1}\|(f,g)\|_{\tilde{L}_{t}^{\infty}(\dot{B}^{\frac{N}{p}}_{p,1})}
\|(f,g)\|_{\tilde{L}_{t}^{q}(\dot{B}^{s}_{p,r})}.
\end{split}\end{equation}
\end{Proposition}

\noindent{\bf Proof.}\, Employing  Taylor's expansion, we have
\begin{equation}
F(f,g)=fF'_{f}(0,0)+gF'_{g}(0,0)+f \tilde{F}_1(f,g)+g \tilde{F}_2(f,g),
\end{equation}
where $\tilde{F}_i: I\times I\rightarrow \mathbb{R}(i=1,2)$ are smooth and vanishes at $(0,0)$. Furthermore, using Propositions \ref{p26} and \ref{p27}, we can conclude \eqref{eq:2.1-A}.\endproof
We  also recall some maximal regularity properties for the heat equation.
\begin{Proposition}\cite{BCD}
\label{Pro:3}
Assume $\mu>0$, $\sigma\in\mathbb{R}, (p,r)\in[1,\infty]^{2}$ and $1\leq \rho_{2}\leq \rho_{1}\leq\infty$. Let $u$ satisfy
\begin{align}\label{eq:heat}
\left\{
\begin{aligned}
&\partial_{t}u-\mu\Delta u=f,\\
&u\mid_{t=0}=u_{0}. \end{aligned} \right.
\end{align}
Then for all $T>0$ the following a priori estimate is fulfilled
\begin{equation}\label{eq:heat1}\mu^{\frac{1}{\rho_{1}}}\|u\|_{\wt L^{\rho_{1}}_T(\dot B^{\sigma+\frac 2\rho_{1}}_{p,r})}\leq C\Big(\|u_{0}\|_{\dot B^{\sigma}_{p,r}}
+\mu^{\frac{1}{\rho_{2}}-1}\|f\|_{\wt L^{\rho_{2}}_T(\dot B^{\sigma-2+\frac 2\rho_{2}}_{p,r})}\Big).
\end{equation}
\end{Proposition}

We finish this subsection by listing an elementary but useful
inequality.
\begin{Lemma}\cite{MN2}\label{lemma2.13}\quad  Let $r_1,r_2>0$ satisfy $\max\{r_1,r_2\}>1$. Then
\begin{equation}\label{99111AA}\int_0^t(1+t-\tau)^{-r_1}(1+\tau)^{-r_2}d\tau\leq C(r_1,r_2)(1+t)^{-\min\{r_1,r_2\}}.\end{equation}
\end{Lemma}
\begin{Lemma}\cite{DX}\label{lemma2.14}\quad  Let $0\leq r_1\leq r_2$ with $ r_2>1$. Then
\begin{equation}\label{99111A}\int_0^t(1+t-\tau)^{-r_1}\tau^{-\theta}(1+\tau)^{\theta-r_2}d\tau\leq C(r_1,r_2)(1+t)^{-r_1}\quad \text{for}\quad 0\leq\theta<1.\end{equation}
\end{Lemma}
\section{The proof of Theorem \ref{th:main1}}
\ \ \ \ \  In this section,  we shall exhibit the  proof of Theorem \ref{th:main1}. We divide it into the following three parts.
\subsection{Maximal regularity estimates in the low frequencies}\ \ \ \ \
Here, we shall establish  the following \emph{a priori}  estimates based on the $L^{2}$-framework of  the  Cauchy problem \eqref{equ:CTFS2}-\eqref{equ:CTFS3}.
\begin{Proposition}\label{low part} Let $T\geq0$,
$N\geq2$, $p$ satisfy $2\leq p\leq\min(4,\frac{2N}{N-2})$ and, additionally, $ p\neq 4$ if $ N=2$.
Assume that $(c^{+},\, u^{+},\, c^{-},\, u^{-})$ is a solution to the  Cauchy problem \eqref{equ:CTFS2}-\eqref{equ:CTFS3}  on $[0,T]\times \mathbb{R}^{N}$,  then
\begin{equation}\label{low part1}\begin{split}\|(c^{+},\, u^{+},\, c^{-},\, &u^{-})\|_{\wt L^\infty_t(\dot B^{\frac{N}{2}-1}_{2,1})}^\ell+\|(c^{+},\, u^{+},\, c^{-},\, u^{-})\|_{\wt L^{1}_t(\dot B^{\frac{N}{2}+1}_{2,1})}^\ell
\\&\lesssim X(0)+ X^{2}(t)+\Big(1+X^{2}(t)\Big)^{[\frac{N}{p}]+1}\Big(X^{2}(t)+X^{3}(t)\Big) \quad \hbox{for all} \quad  t\in[0,T].\end{split}\end{equation}
\end{Proposition}

\noindent{\bf Proof.}\, Employing the energy argument of Godunov
\cite{GO} for partially dissipative first-order symmetric systems (further developed by
 \cite{KF}),  by a similar derivation in Lemma 3.1 in  \cite{XC}, we conclude that
 \begin{equation}\label{gama8}
\big|(\widehat{c^{+}},\widehat{u^{+}}, \widehat{c^{-}}, \widehat{u^{-}})\big|
\leq Ce^{-c_{0}|\xi|^{2}t}\big|(\widehat{c^{+}}, \widehat{u^{+}}, \widehat{c^{-}}, \widehat{u^{-}})\big|(0)\quad \text{for}\quad \xi\leq \xi_0,
\end{equation}
 which together with  Fourier-Plancherel
theorem and Duhamel's formula, implies that
\begin{equation}\label{gama10}\begin{split}
&\|(c^{+},u^{+},c^{-},u^{-})\|
_{\wt L^\infty_t(\dot B^{\frac{N}{2}-1}_{2,1})}^\ell
+\|(c^{+},u^{+},c^{-},u^{-})
\|_{L^{1}_{t}(\dot{B}^{\f{N}{2}+1}_{2,1})}^{\ell}
\\&\lesssim \|(c^{+}_0,u^{+}_0,c^{-}_0,u^{-}_0)\|
_{\dot B^{\frac{N}{2}-1}_{2,1}}^\ell+\|(H_{1},H_{2},H_{3},H_{4})
\|_{L^{1}_{t}(\dot{B}^{\f{N}{2}-1}_{2,1})}^{\ell}.
\end{split}\end{equation}
In what follows, we derive some estimates for the nonlinear terms
$\|(H_{1},H_{2},H_{3},H_{4})
\|_{L^{1}_{t}(\dot{B}^{\f{N}{2}-1}_{2,1})}^{\ell}$.

For the sake of simplicity, we first show the following five  important estimates of the paradifferential calculus  from Propositions \ref{200} and \ref{Pro:201} respectively  which will be frequently used in our process later.
\begin{equation}\label{estmates1}T:\dot B^{\frac Np-1}_{p,1}\times \dot B^{\frac Np+1}_{p,1}\rightarrow \dot B^{\frac N2}_{2,1}\quad  \text{for}\quad  2\leq p\leq\min\bigl(4,\frac{2N}{N-2}\bigr),
\end{equation}
\begin{equation}\label{estmates2}R: \dot B^{\frac Np}_{p,1}\times \dot B^{\frac Np}_{p,1}\rightarrow \dot B^{\frac N2}_{2,1} \quad  \text{for}\quad  2\leq p\leq4,
\end{equation}
\begin{equation}\label{estmates3}T: \dot B^{\frac Np-1}_{p,1}\times \dot B^{\frac Np}_{p,1}\rightarrow\dot B^{\frac N2-1}_{2,1}\quad  \text{for}\quad  2\leq p\leq\min\bigl(4,\frac{2N}{N-2}\bigr),
\end{equation}
\begin{equation}\label{estmates4}R:\dot B^{\frac Np-1}_{p,1}\times \dot B^{\frac Np}_{p,1}\rightarrow\dot B^{\frac N2-1}_{2,1}\quad  \text{for}\quad  2\leq p\leq4,
\end{equation} and
\begin{equation}\label{estmates5}T:\dot B^{\frac Np-1}_{p,1}\times \dot B^{\frac Np-1}_{p,1}\rightarrow\dot B^{\frac N2-2}_{2,1}\quad  \text{for}\quad  2\leq p\leq\min\bigl(4,\frac{2N}{N-2}\bigr).
\end{equation}
For the term $H_{1}$,  using Bony's decomposition  we see that
$$(c^{+}u^{+})^{\ell}=(T_{c^{+}}u^{+})^{\ell}+(R(c^{+}u^{+}))^{\ell}
+(T_{u^{+}}c^{+})^{\ell}.$$
Employing \eqref{estmates1}
and  the embedding relation  $\dot B^{\frac N2+s}_{2,1}\times \dot B^{\frac Np+s}_{p,1}$ for $s\in \mathbb{R}$, $p\geq2$ yields that
\begin{equation*}\label{eq:c-u-1}
\begin{split}
&\|(T_{c^{+}}u^{+})^{\ell}\|_{L^1_t(\dot B^{\frac N2}_{2,1})}\lesssim\|c^{+}\|_{L^\infty_t(\dot B^{\frac Np-1}_{p,1})}\|u^{+}\|_{L^1_t(\dot B^{\frac Np+1}_{p,1})}\\
&\qquad\lesssim\Big(\|c^{+}\|^\ell_{L^\infty_t(\dot B^{\frac Np-1}_{p,1})}+\|c^{+}\|^h_{L^\infty_t(\dot B^{\frac Np-1}_{p,1})}\Big)\Big(\|u^{+}\|^\ell_{L^1_t(\dot B^{\frac Np+1}_{p,1})}+\|u^{+}\|^h_{L^1_t(\dot B^{\frac Np+1}_{p,1})}\big)\\
&\qquad\lesssim\Big(\|c^{+}\|^\ell_{L^\infty_t(\dot B^{\frac N2-1}_{2,1})}+\|c^{+}\|^h_{L^\infty_t(\dot B^{\frac Np}_{p,1})}\Big)\Big(\|u^{+}\|^\ell_{L^1_t(\dot B^{\frac N2+1}_{2,1})}+\|u^{+}\|^h_{L^1_t(\dot B^{\frac Np+1}_{p,1})}\Big)\\
&\qquad\lesssim X^{2}(t),
\end{split}
\end{equation*}
and
\begin{equation*}\label{eq:c-u-1}
\begin{split}
&\|(T_{u^{+}}c^{+})^{\ell}\|_{L^1_t(\dot B^{\frac N2}_{2,1})}
\lesssim\|u^{+}\|_{L^\infty_t(\dot B^{\frac Np-1}_{p,1})}\|c^{+}\|_{L^1_t(\dot B^{\frac Np+1}_{p,1})}\\
&\qquad\lesssim\Big(\|u^{+}\|^\ell_{L^\infty_t(\dot B^{\frac Np-1}_{p,1})}+\|u^{+}\|^h_{L^\infty_t(\dot B^{\frac Np-1}_{p,1})}\Big)\Big(\|c^{+}\|^\ell_{L^1_t(\dot B^{\frac Np+1}_{p,1})}+\|c^{+}\|^h_{L^1_t(\dot B^{\frac Np+1}_{p,1})}\big)\\
&\qquad\lesssim\Big(\|u^{+}\|^\ell_{L^\infty_t(\dot B^{\frac N2-1}_{2,1})}+\|u^{+}\|^h_{L^\infty_t(\dot B^{\frac Np-1}_{p,1})}\Big)\Big(\|c^{+}\|^\ell_{L^1_t(\dot B^{\frac N2+1}_{2,1})}+\|c^{+}\|^h_{L^1_t(\dot B^{\frac Np+2}_{p,1})}\Big)\\
&\qquad\lesssim X^{2}(t).
\end{split}
\end{equation*}
For the remainder term, using \eqref{estmates2}  and the embedding $\dot B^{\frac N2+s}_{2,1}\times \dot B^{\frac Np+s}_{p,1}$ for $s\in \mathbb{R}$, $p\geq2$ gives rise to
\begin{equation*}
\begin{split}
&\|(R(c^{+}u^{+}))^{\ell}\|_{{L^1_t}
(\dot{B}^{\f{N}{2}}_{2,1})}
\\&\lesssim\|c^{+}\|_{L^2_{t}
(\dot{B}^{\f{N}{p}}_{p,1})}
\|u^{+}\|_{L^2_t(\dot{B}^{\f{N}{p}}_{p,1})}
\\&\lesssim\Big(\|c^{+}\|_{L^2_{t}
(\dot{B}^{\f{N}{p}}_{p,1})}^{\ell}+\|c^{+}\|_{L^2_{t}
(\dot{B}^{\f{N}{p}}_{p,1})}^{h}\Big)
\Big(\|u^{+}\|_{L^2_t(\dot{B}^{\f{N}{p}}_{p,1})}^{\ell}+\|u^{+}\|_{L^2_t(\dot{B}^{\f{N}{p}}_{p,1})}^{h}\Big)
\\&\lesssim\Big(\|c^{+}\|_{L^2_{t}
(\dot{B}^{\f{N}{2}}_{2,1})}^{\ell}+\|c^{+}\|_{L^2_{t}
(\dot{B}^{\f{N}{p}+1}_{p,1})}^{h}\Big)
\Big(\|u^{+}\|_{L^2_t(\dot{B}^{\f{N}{2}}_{2,1})}^{\ell}+\|u^{+}\|_{L^2_t(\dot{B}^{\f{N}{p}}_{p,1})}^{h}\Big)
\\&\lesssim X^{2}(t),
\end{split}
\end{equation*}
where we have used the following interpolation inequalities (which will be frequently used later),
$$\|c^{+}\|^\ell_{\tilde{L}^2_t(\dot B^{\frac N2}_{2,1})}\leq \Big(\|c^{+}\|^\ell_{\tilde{L}^\infty_t(\dot B^{\frac N2-1}_{2,1})}\Big)^{\frac{1}{2} }\Big(\|c^{+}\|^\ell_{L^1_t(\dot B^{\frac N2+1}_{2,1})}\Big)^{\frac{1}{2}}, $$
$$\|c^{+}\|^h_{\tilde{L}^2_t(\dot B^{\frac Np+1}_{p,1})}\leq \Big(\|c^{+}\|^h_{\tilde{L}^\infty_t(\dot B^{\frac Np}_{p,1})}\Big)^{\frac{1}{2} }\Big(\|c^{+}\|^h_{L^1_t(\dot B^{\frac Np+2}_{p,1})}\Big)^{\frac{1}{2}}, $$
$$\|u^{+}\|^\ell_{\tilde{L}^2_t(\dot B^{\frac N2}_{2,1})}\leq \Big(\|u^{+}\|^\ell_{\tilde{L}^\infty_t(\dot B^{\frac N2-1}_{2,1})}\Big)^{\frac{1}{2} }\Big(\|u^{+}\|^\ell_{L^1_t(\dot B^{\frac N2+1}_{2,1})}\Big)^{\frac{1}{2}},$$
and $$\|u^{+}\|^h_{\tilde{L}^2_t(\dot B^{\frac Np}_{p,1})}\leq \Big(\|u^{+}\|^h_{\tilde{L}^\infty_t(\dot B^{\frac Np-1}_{p,1})}\Big)^{\frac{1}{2} }\Big(\|u^{+}\|^h_{L^1_t(\dot B^{\frac Np+1}_{p,1})}\Big)^{\frac{1}{2}}.$$
Hence
\begin{equation*}
\begin{split}\|H_{1}\|_{L^{1}([0,t];\dot{B}^{\f{N}{2}-1}_{2,1})}^{\ell}
\lesssim\|(c^{+}u^{+})^{\ell}\|_{L^{1}([0,t];\dot{B}^{\f{N}{2}}_{2,1})}
\lesssim X^{2}(t).
\end{split}
\end{equation*}
The term  $H_{3}$  may be treated along the same lines,  and we omit it.

For the  term  $(g_{+}(c^{+},c^{-})\partial_{i} c^{+})^{\ell}$ in  $H_{2}^{i}$,  employing Bony's decomposition and splitting $c^{+}$ into ${c^{+}}^{\ell}+{c^{+}}^{h}$
 we have
\begin{equation*}
\begin{split}(g_{+}(c^{+},c^{-})\partial_{i} c^{+})^{\ell}
&=(T_{\partial_{i} c^{+}}g_{+}(c^{+},c^{-}))^{\ell}
+(R(g_{+}(c^{+},c^{-}),\partial_{i} c^{+}))^{\ell}
\\&\quad+(T_{g_{+}(c^{+},c^{-})}\partial_{i} {c^{+}}^{\ell})^{\ell}
+(T_{g_{+}(c^{+},c^{-})}\partial_{i} {c^{+}}^{h})^{\ell}.
\end{split}
\end{equation*}
We bound  the above  items one by one.  To bound the first two terms,  it  suffices to notice that  \eqref{estmates3}  and \eqref{estmates4} and use  Proposition  \ref{p27}\mbox{(i)}. Hence
\begin{equation*}
\begin{split}
&\|(T_{\partial_{i} c^{+}}g_{+}(c^{+},c^{-}))^{\ell}
+(R(g_{+}(c^{+},c^{-}),\partial_{i} c^{+}))^{\ell}\|_{{L^1_t}
(\dot{B}^{\f{N}{2}-1}_{2,1})}
\\&\lesssim\|\nabla c^{+}\|_{L^2_{t}
(\dot{B}^{\f{N}{p}-1}_{p,1})}
\|g_{+}(c^{+},c^{-})\|_{L^2_t(\dot{B}^{\f{N}{p}}_{p,1})}
\\&\lesssim \Big(1+X^{2}(t)\Big)^{[\frac{N}{p}]+1}\| (c^{+},c^{-})\|_{\tilde{L}^2_{t}
(\dot{B}^{\f{N}{p}}_{p,1})}^{2}
\\&\lesssim \Big(1+X^{2}(t)\Big)^{[\frac{N}{p}]+1}X^{2}(t).
\end{split}
\end{equation*}
 Employing $\|T_{f}g\|_{
\dot{B}^{s}_{p,r}}
\lesssim\|f\|_{L^{\infty}}\|g\|_{\dot{B}^{s}_{p,r}}$, imbedding relation $\dot{B}^{\f{N}{p}}_{p,1}\hookrightarrow L^{\infty}$  and Proposition \ref{p27} \mbox{(i)}, we handle  the third term  as follows
\begin{equation*}
\begin{split}
&\|(T_{g_{+}(c^{+},c^{-})}\partial_{i} {c^{+}}^{\ell})^{\ell}\|_{{L^1_t}
(\dot{B}^{\f{N}{2}-1}_{2,1})}
\\&\lesssim\|\nabla c^{+\ell}\|_{L^2_{t}
(\dot{B}^{\f{N}{2}-1}_{2,1})}
\|g_{+}(c^{+},c^{-})\|_{L^2_t(L^{\infty})}
\\&\lesssim\|c^{+\ell}\|_{L^2_{t}
(\dot{B}^{\f{N}{2}}_{2,1})}
\|g_{+}(c^{+},c^{-})\|_{L^2_t(\dot{B}^{\f{N}{p}}_{p,1})}
\\&\lesssim \Big(1+X^{2}(t)\Big)^{[\frac{N}{p}]+1}X^{2}(t).
\end{split}
\end{equation*}
The last term may be bounded thanks to \eqref{estmates3}. Moreover,  we just have to use Proposition \ref{p28} as it may happen that $\frac Np-1<0$.  Thus
\begin{equation*}
\begin{split}
&\|(T_{g_{+}(c^{+},c^{-})}\partial_{i} {c^{+}}^{h})^{\ell}\|_{{L^1_t}
(\dot{B}^{\f{N}{2}-1}_{2,1})}
\\&\lesssim\|g_{+}(c^{+},c^{-})\|_{L^\infty_t(\dot{B}^{\f{N}{p}-1}_{p,1})}
\|\nabla c^{+h}\|_{L^1_{t}
(\dot{B}^{\f{N}{p}-1}_{p,1})}
\\&\lesssim\Big(\|(c^{+},c^{-})\|_{\tilde{L}^\infty_t(\dot{B}^{\f{N}{p}-1}_{p,1})}
+\big(1+\|(c^{+},c^{-})\|^{2}_{\tilde{L}^\infty_t(\dot{B}^{\f{N}{p}}_{p,1})}\big)^{[\frac{N}{p}+1]}\|(c^{+},c^{-})\|_{\tilde{L}^\infty_t(\dot{B}^{\f{N}{p}}_{p,1})}\|(c^{+},c^{-})\|_{\tilde{L}^\infty_t(\dot{B}^{\f{N}{p}-1}_{p,1})}\Big)
\\&\quad\times\| c^{+}\|^{h}_{L^1_{t}
(\dot{B}^{\f{N}{p}+2}_{p,1})}
\\&\lesssim X^{2}(t)+\Big(1+X^{2}(t)\Big)^{[\frac{N}{p}]+1}X^{3}(t).
\end{split}
\end{equation*}
Hence
\begin{equation*}
\|(g_{+}(c^{+},c^{-})\partial_{i} c^{+})^{\ell}\|_{L^{1}([0,t];\dot{B}^{\f{N}{2}-1}_{2,1})}
\lesssim X^{2}(t)+\Big(1+X^{2}(t)\Big)^{[\frac{N}{p}]+1}\Big(X^{2}(t)+X^{3}(t)\Big).
\end{equation*}
Bounding the  term   $(\tilde{g}_{+}(c^{+},c^{-})\partial_{i} c^{-})^{\ell}$ is totally similar, we omit it.

For the term with $((u^{+}\cdot\nabla)u_{i}^{+})^{\ell}$ in  $H_{2}^{i}$, employing Bony's decomposition implies that
$$((u^{+}\cdot\nabla)u_{i}^{+})^{\ell}=(T_{\nabla u_{i}^{+}}u^{+})^{\ell}
+\big(\sum_{j=1}^{N}(R(u_{i}^{+},\partial_{j}u_{i}^{+}))\big)^{\ell}
+(T_{u^{+}}\nabla u_{i}^{+})^{\ell}.
$$
Then it follows from \eqref{estmates3}  and \eqref{estmates4}  that
\begin{equation*}
\begin{split}
&\|(T_{\nabla u_{i}^{+}}u^{+})^{\ell}
+(R(u_{i}^{+},\partial_{j}u_{i}^{+}))^{\ell}\|_{{L^1_t}
(\dot{B}^{\f{N}{2}-1}_{2,1})}+\|(T_{u^{+}}\nabla u_{i}^{+})^{\ell}\|_{{L^1_t}
(\dot{B}^{\f{N}{2}-1}_{2,1})}
\\&\lesssim\|\nabla u^{+}\|_{L^2_{t}
(\dot{B}^{\f{N}{p}-1}_{p,1})}
\|u^{+}\|_{L^2_t(\dot{B}^{\f{N}{p}}_{p,1})}
\\&\lesssim\|u^{+}\|_{L^2_{t}
(\dot{B}^{\f{N}{p}}_{p,1})}^{2}+\|u^{+}\|_{L^\infty_t(\dot{B}^{\f{N}{p}-1}_{p,1})}
\|\nabla u^{+}\|_{L^1_{t}
(\dot{B}^{\f{N}{p}}_{p,1})}
\\&\lesssim\|u^{+}\|_{L^\infty_{t}
(\dot{B}^{\f{N}{p}-1}_{p,1})}
\|u^{+}\|_{L^1_{t}
(\dot{B}^{\f{N}{p}+1}_{p,1})}
\\&\lesssim X^{2}(t).
\end{split}
\end{equation*}
Hence
\begin{equation*}
\|((u^{+}\cdot\nabla)u_{i}^{+})^{\ell}\|_{L^{1}([0,t];\dot{B}^{\f{N}{2}-1}_{2,1})}
\lesssim X^{2}(t).
\end{equation*}

For the term  $(h_{+}(c^{+},c^{-})\partial_{j}c^{+}\partial_{j}u^{+}_{i})^{\ell}$ in  $H_{2}^{i}$, we decompose it into
\begin{equation*}
\begin{split}
(h_{+}(c^{+},c^{-})\partial_{j}c^{+}\partial_{j}u^{+}_{i})^{\ell}
&=(T_{h_{+}(c^{+},c^{-})\partial_{j}c^{+}}\partial_{j}u^{+}_{i})^{\ell}
+(R(h_{+}(c^{+},c^{-})\partial_{j}c^{+},\partial_{j}u^{+}_{i}))^{\ell}
\\&+(T_{\partial_{j}u^{+}_{i}}h_{+}(c^{+},c^{-})\partial_{j}c^{+\ell})^{\ell}
+(T_{\partial_{j}u^{+}_{i}}h_{+}(c^{+},c^{-})\partial_{j}c^{+h})^{\ell}.
\end{split}
\end{equation*}
Thanks to \eqref{estmates3}  and \eqref{estmates4}, Proposition \ref{p27} \mbox{(i)}, we deduce that
\begin{equation*}
\begin{split}
&\|(T_{h_{+}(c^{+},c^{-})\partial_{j}c^{+}}\partial_{j}u^{+}_{i})^{\ell}
+(R(h_{+}(c^{+},c^{-})\partial_{j}c^{+},\partial_{j}u^{+}_{i}))^{\ell}\|_{{L^1_t}
(\dot{B}^{\f{N}{2}-1}_{2,1})}
\\&\lesssim\|h_{+}(c^{+},c^{-})\partial_{j}c^{+}\|_{L^\infty_{t}
(\dot{B}^{\f{N}{p}-1}_{p,1})}
\|\nabla u^{+}\|_{L^1_t(\dot{B}^{\f{N}{p}}_{p,1})}
\\&\lesssim\|h_{+}(c^{+},c^{-})\|_{L^\infty_{t}
(\dot{B}^{\f{N}{p}}_{p,1})}
\|\nabla c^{+}\|_{L^\infty_{t}
(\dot{B}^{\f{N}{p}-1}_{p,1})}
\|\nabla u^{+}\|_{L^1_t(\dot{B}^{\f{N}{p}}_{p,1})}
\\&\lesssim\Big(1+\big(1+\|(c^{+},c^{-})\|^{2}_{\tilde{L}^\infty_t(\dot{B}^{\f{N}{p}}_{p,1})}\big)^{[\frac{N}{p}+1]}\|(c^{+},c^{-})\|_{\tilde{L}^\infty_t(\dot{B}^{\f{N}{p}}_{p,1})}\Big)
\| c^{+}\|_{L^\infty_{t}
(\dot{B}^{\f{N}{p}}_{p,1})}
\|u^{+}\|_{L^1_{t}
(\dot{B}^{\f{N}{p}+1}_{p,1})}
\\&\lesssim X^{2}(t)+\Big(1+X^{2}(t)\Big)^{[\frac{N}{p}]+1}X^{3}(t),
\end{split}
\end{equation*}
\begin{equation*}
\begin{split}
&\|(T_{\partial_{j}u^{+}_{i}}h_{+}(c^{+},c^{-})
\partial_{j}c^{+\ell})^{\ell}\|_{{L^1_t}
(\dot{B}^{\f{N}{2}-1}_{2,1})}
\\&\lesssim\|\nabla u^{+}\|_{L^\infty_{t}
(\dot{B}^{\f{N}{p}-1}_{p,1})}
\|h_{+}(c^{+},c^{-})
\nabla c^{+\ell}\|_{L^1_t(\dot{B}^{\f{N}{p}}_{p,1})}
\\&\lesssim\|\nabla u^{+}\|_{L^\infty_{t}
(\dot{B}^{\f{N}{p}-1}_{p,1})}
\Big(1+\big(1+\|(c^{+},c^{-})\|^{2}_{\tilde{L}^\infty_t(\dot{B}^{\f{N}{p}}_{p,1})}\big)^{[\frac{N}{p}+1]}\|(c^{+},c^{-})\|_{\tilde{L}^\infty_t(\dot{B}^{\f{N}{p}}_{p,1})}\Big)
\|\nabla c^{+\ell}\|_{L^1_t(\dot{B}^{\f{N}{p}}_{p,1})}
\\&\lesssim X^{2}(t)+\Big(1+X^{2}(t)\Big)^{[\frac{N}{p}]+1}X^{3}(t),
\end{split}
\end{equation*}and
\begin{equation*}
\begin{split}
&\|(T_{h_{+}(c^{+},c^{-})\partial_{j}u^{+}_{i}}\partial_{j}c^{h})^{\ell}
\|_{{L^1_t}(\dot{B}^{\f{N}{2}-1}_{2,1})}
\\&\lesssim\|h_{+}(c^{+},c^{-})\partial_{j}u^{+}_{i}\|_{L^1_{t}
(\dot{B}^{\f{N}{p}}_{p,1})}
\|\partial_{j}c^{+h}\|_{L^\infty_t(\dot{B}^{\f{N}{p}-1}_{p,1})}
\\&\lesssim
\Big(1+\big(1+\|(c^{+},c^{-})\|^{2}_{\tilde{L}^\infty_t(\dot{B}^{\f{N}{p}}_{p,1})}\big)^{[\frac{N}{p}+1]}\|(c^{+},c^{-})\|_{\tilde{L}^\infty_t(\dot{B}^{\f{N}{p}}_{p,1})}\Big)
\|\nabla u^{+}\|_{L^1_{t}
(\dot{B}^{\f{N}{p}}_{p,1})}
\|c^{+h}\|_{L^\infty_t(\dot{B}^{\f{N}{p}}_{p,1})}
\\&\lesssim X^{2}(t)+\Big(1+X^{2}(t)\Big)^{[\frac{N}{p}]+1}X^{3}(t).
\end{split}
\end{equation*}
Hence
\begin{equation*}
\|(h_{+}(c^{+},c^{-})\partial_{j}c^{+}\partial_{j}u^{+}_{i})^{\ell}\|_{L^{1}([0,t];\dot{B}^{\f{N}{2}-1}_{2,1})}
\lesssim X^{2}(t)+\Big(1+X^{2}(t)\Big)^{[\frac{N}{p}]+1}X^{3}(t).
\end{equation*}
Similarly,  we also obtain  the corresponding estimates of other terms
$(\mu^{+}k_{+}(c^{+},c^{-})\partial_{j}c^{-}\partial_{j}u^{+}_{i})^{\ell},$\\
$(\mu^{+}h_{+}(c^{+},c^{-})\partial_{j}c^{+}\partial_{i}u^{+}_{j})^{\ell},$
$(\mu^{+}k_{+}(c^{+},c^{-})\partial_{j}c^{-}\partial_{i}u^{+}_{j})^{\ell},
(\lambda^{+}h_{+}(c^{+},c^{-})\partial_{i}c^{+}\partial_{j}u^{+}_{j})^{\ell}$ and\\
$(\lambda^{+}k_{+}(c^{+},c^{-})\partial_{i}c^{-}\partial_{j}u^{+}_{j})^{\ell}$. Here, we omit the details.

Finally,  for the term $(l_{+}(c^{+},c^{-})\partial_{j}^{2}u_{i}^{+})^{\ell}$ in  $H_{2}^{i}$,
we  may decompose  it into
\begin{equation*}
\begin{split}
(l_{+}(c^{+},c^{-})\partial_{j}^{2}u_{i}^{+})^{\ell}
&=(T_{\partial_{j}^{2}u_{i}^{+}}l_{+}(c^{+},c^{-}))^{\ell}
+(R(l_{+}(c^{+},c^{-}),\partial_{j}^{2}u_{i}^{+}))^{\ell}
\\&\quad+(T_{l_{+}(c^{+},c^{-})}\partial_{j}^{2}u_{i}^{+\ell})^{\ell}
+(T_{l_{+}(c^{+},c^{-})}\partial_{j}^{2}u_{i}^{+h})^{\ell}.
\end{split}
\end{equation*}
To handle the first two terms, by virtue of  \eqref{estmates3}  and \eqref{estmates4},  Proposition \ref{p27}\mbox{(i)},
we get
\begin{equation*}
\begin{split}
&\|(T_{\partial_{j}^{2}u_{i}^{+}}l_{+}(c^{+},c^{-}))^{\ell}
+(R(l_{+}(c^{+},c^{-}),\partial_{j}^{2}u_{i}^{+}))^{\ell}\|_{{L^1_t}
(\dot{B}^{\f{N}{2}-1}_{2,1})}
\\&\lesssim\|\nabla^{2} u^{+}\|_{L^1_{t}
(\dot{B}^{\f{N}{p}-1}_{p,1})}
\|l_{+}(c^{+},c^{-})\|_{L^\infty_t(\dot{B}^{\f{N}{p}}_{p,1})}
\\&\lesssim\| u^{+}\|_{L^1_{t}
(\dot{B}^{\f{N}{p}+1}_{p,1})}
\big(1+\|(c^{+},c^{-})\|^{2}_{\tilde{L}^\infty_t(\dot{B}^{\f{N}{p}}_{p,1})}\big)^{[\frac{N}{p}+1]}\|(c^{+},c^{-})\|_{\tilde{L}^\infty_t(\dot{B}^{\f{N}{p}}_{p,1})}
\\&\lesssim X^{2}(t)+\Big(1+X^{2}(t)\Big)^{[\frac{N}{p}]+1}X^{2}(t).
\end{split}
\end{equation*}
To bound the third therm,  according to $\|T_{f}g\|_{
\dot{B}^{s}_{p,r}}
\lesssim\|f\|_{L^{\infty}}\|g\|_{\dot{B}^{s}_{p,r}}$, imbedding relation $\dot{B}^{\f{N}{p}}_{p,1}\hookrightarrow L^{\infty}$  and Proposition \ref{p27}\mbox{(i)}, we infer that
\begin{equation*}
\begin{split}
&\|(T_{l_{+}(c^{+},c^{-})}\partial_{j}^{2}u_{i}^{+\ell})^{\ell}
\|_{{L^1_t}(\dot{B}^{\f{N}{2}-1}_{2,1})}
\\&\lesssim\|l_{+}(c^{+},c^{-})\|_{L^\infty_{t}
(L^\infty)}
\|\nabla^{2}u^{+\ell}\|_{L^1_t(\dot{B}^{\f{N}{2}-1}_{2,1})}
\\&\lesssim
\|(c^{+},c^{-})\|_{L^\infty_t(\dot{B}^{\f{N}{p}}_{p,1})}
\|u^{+\ell}\|_{L^1_{t}(\dot{B}^{\f{N}{2}+1}_{2,1})}
\\&\lesssim \big(1+\|(c^{+},c^{-})\|^{2}_{\tilde{L}^\infty_t(\dot{B}^{\f{N}{p}}_{p,1})}\big)^{[\frac{N}{p}+1]}\|(c^{+},c^{-})\|_{\tilde{L}^\infty_t(\dot{B}^{\f{N}{p}}_{p,1})}\|u^{+}\|^{\ell}_{L^1_{t}(\dot{B}^{\f{N}{2}+1}_{2,1})},
\\&\lesssim X^{2}(t)+\Big(1+X^{2}(t)\Big)^{[\frac{N}{p}]+1}X^{2}(t).
\end{split}
\end{equation*}
Using  \eqref{estmates5} and  Proposition \ref{p28}, we bound  the last term as follows
\begin{equation*}
\begin{split}
&\|(T_{l_{+}(c^{+},c^{-})}\partial_{j}^{2}u_{i}^{+h})^{\ell}
\|_{{L^1_t}(\dot{B}^{\f{N}{2}-1}_{2,1})}
\\&\lesssim\|(T_{l_{+}(c^{+},c^{-})}\partial_{j}^{2}u_{i}^{+h})^{\ell}
\|_{{L^1_t}(\dot{B}^{\f{N}{2}-2}_{2,1})}
\\&\lesssim\|l_{+}(c^{+},c^{-})\|_{L^\infty_{t}
(\dot{B}^{\f{N}{p}-1}_{p,1})}
\|\nabla^{2}u^{+h}\|_{L^1_t(\dot{B}^{\f{N}{p}-1}_{p,1})}
\\&\lesssim
\Big(\|(c^{+},c^{-})\|_{\tilde{L}^\infty_t(\dot{B}^{\f{N}{p}-1}_{p,1})}
+\big(1+\|(c^{+},c^{-})\|^{2}_{\tilde{L}^\infty_t(\dot{B}^{\f{N}{p}}_{p,1})}\big)^{[\frac{N}{p}+1]}
\|(c^{+},c^{-})\|_{\tilde{L}^\infty_t(\dot{B}^{\f{N}{p}}_{p,1})}\|(c^{+},c^{-})\|_{\tilde{L}^\infty_t(\dot{B}^{\f{N}{p}-1}_{p,1})}\Big)
\\&\qquad\times\| u^{+}\|^{h}_{L^1_{t}
(\dot{B}^{\f{N}{p}+1}_{p,1})}
\\&\lesssim X^{2}(t)+\Big(1+X^{2}(t)\Big)^{[\frac{N}{p}]+1}X^{3}(t).
\end{split}
\end{equation*}
Hence
\begin{equation*}
\|H_{2}\|_{L^{1}([0,t];\dot{B}^{\f{N}{2}-1}_{2,1})}^{\ell}
\lesssim  X^{2}(t)+\Big(1+X^{2}(t)\Big)^{[\frac{N}{p}]+1}\Big(X^{2}(t)+X^{3}(t)\Big).
\end{equation*}
The term  $H_{4}$ may be treated along the same lines, we omit it.

Putting together all the above estimates for the terms of $H_{1}$-$H_{4}$, we  deduce that \eqref{low part1}. This completes the proof of Proposition \ref{low part}. \endproof
\subsection{Maximal regularity estimates in the high frequencies}\ \ \ \ \
In the following proposition, we shall  exploit the parabolic properties of  the  Cauchy problem \eqref{equ:CTFS2}-\eqref{equ:CTFS3} in the high frequencies and construct  \emph{a priori}  estimates based on the general $L^{p}$-framework.
\begin{Proposition}\label{high part} Let $T\geq0$,
$N\geq2$, $p$ satisfy $2\leq p\leq\min(4,\frac{2N}{N-2})$ and, additionally, $ p\neq 4$ if $ N=2$.
Assume that $(c^{+},\, u^{+},\, c^{-},\, u^{-})$ is a solution to the  Cauchy problem \eqref{equ:CTFS2}-\eqref{equ:CTFS3} on $[0,T]\times \mathbb{R}^{N}$,  then
\begin{equation}
\begin{split}
\label{high part1}
\|(u^{+},\,u^{-})\|_{\wt L^\infty_t(\dot B^{\frac{N}{p}-1}_{p,1})}^h
&+\|(u^{+},\,u^{-})\|_{L^1_t(\dot B^{\frac{N}{p}+1}_{p,1})}^h
+\|(c^{+},\,c^{-})\|_{\wt L^\infty_t(\dot B^{\frac{N}{p}}_{p,1})}^h
+\|(c^{+},\,c^{-})\|_{L^1_t(\dot B^{\frac{N}{p}+2}_{p,1})}^h
\\&\lesssim X(0)+X^{2}(t)+\Big(1+X^{2}(t)\Big)^{[\frac{N}{p}]+1}\Big(X^{2}(t)+X^{3}(t)\Big) \quad \hbox{for all} \quad  t\in[0,T].
\end{split}
\end{equation}
\end{Proposition}

 In order to prove Proposition \ref{high part}, we  introduce the orthogonal projectors over divergence-free $\cP$ and potential vector fields $\cQ$ satisfying the identity $I =$\cP$ +\cQ$. Then, applying the orthogonal projectors    $\cP$ and $\cQ$  over divergence-free and potential vector-fields, respectively,
 to \eqref{equ:CTFS2}$_{2}$  and \eqref{equ:CTFS2}$_{4}$,
and setting  $\nu^{\pm}\eqdefa\nu_{1}^{\pm}+\nu_{2}^{\pm}$ yield that
\begin{equation}\label{gama2-A}
\left\{
\begin{aligned}{}
&\p_t\cP u^{+}-\nu_{1}^{+}\Delta \cP u^{+}=\cP H_{2},
\\&\p_t\cP u^{-}-\nu_{1}^{-}\Delta \cP u^{-}=\cP H_{4},
\end{aligned}
\right.
\end{equation}
and
\begin{equation}\label{gama3-A}
\left\{
\begin{aligned}{}
&\p_tc^{+}+\textrm{div}\cQ u^{+}=H_{1},\\
&\p_t\cQ u^{+}+\beta_{1}\nabla c^{+}
+\beta_{2}\nabla c^{-}-\nu^{+}\Delta\cQ u^{+}-\nabla \Delta c^{+}=\cQ H_{2},
\\&\p_tc^{-}+\textrm{div}\cQ u^{-}=H_{3},
\\&\p_t\cQ u^{-}+\beta_{3}\nabla c^{+}
+\beta_{4}\nabla c^{-}-\nu^{-}\Delta \cQ u^{-}-\nabla \Delta c^{-}=\cQ H_{4}.
\end{aligned}
\right.
\end{equation}
 For system \eqref{gama2-A}, according to Proposition \ref{Pro:3} (restricted to the high frequencies), we have
\begin{equation}\label{eq:HFo}
\begin{split}
\|&\cP  u^{\pm}\|^h_{\wt L^\infty_t(\dot B^{\frac Np-1}_{p,1})}+\nu_{1}^{\pm} \|\cP  u^{\pm}\|^h_{L^1_t(\dot B^{\frac Np+1}_{p,1})}\\
&\lesssim\Big( \| \cP  u^{\pm}_0\|_{\dot B^{\frac Np-1}_{p,1}}^h+\|(\cP  H_{2},\cP  H_{4})\|^h_{L^1_t(\dot B^{\frac Np-1}_{p,1})}\Big).
\end{split}
\end{equation}
Next, to handle the coupling system \eqref{gama3-A} in the high frequencies,   applying operator $\Delta$  to  \eqref{gama3-A}$_{1}$   and  \eqref{gama3-A}$_{3}$,     and taking   operator  $\dv$   to  \eqref{gama3-A}$_{2}$ and \eqref{gama3-A}$_{4}$, respectively,   we get
\begin{equation}\label{gama3-A-1}
\left\{
\begin{aligned}{}
&\p_t\Delta c^{+}+\Delta\textrm{div}\cQ u^{+}=\Delta H_{1},\\
&\p_t\textrm{div}\cQ u^{+}-\nu^{+}\Delta\textrm{div}\cQ u^{+}-\Delta \Delta c^{+}=\textrm{div}\cQ H_{2}-\beta_{1}\Delta c^{+}
-\beta_{2}\Delta c^{-},
\\&\p_t\Delta c^{-}+\Delta\textrm{div}\cQ u^{-}=\Delta H_{3},
\\&\p_t\textrm{div}\cQ u^{-}-\nu^{-}\Delta \textrm{div}\cQ u^{-}-\Delta \Delta c^{-}=\textrm{div}\cQ H_{4}-\beta_{3}\Delta c^{+}
-\beta_{4}\Delta c^{-}.
\end{aligned}
\right.
\end{equation}
Obviously, $(\Delta c^{\pm}, \textrm{div}\cQ u^{\pm})$  in system \eqref{gama3-A-1} satisfies the linearized coupling  system
\begin{equation}\label{gama3-A-2-1}
\left\{
\begin{aligned}{}
&\p_ta+\Delta v=F,\\
&\p_tv-\nu\Delta v-\Delta  a=G.
\end{aligned}
\right.
\end{equation}
 In what follows, we will prove   system \eqref{gama3-A-2-1} has   the parabolic properties  in  $L^{p}$-framework in the high frequencies,  which implies that the terms  $\beta_{1}\Delta c^{+}+\beta_{2}\Delta c^{-}$ and $\beta_{3}\Delta c^{+}+\beta_{4}\Delta c^{-}$ on the  right hands of system \eqref{gama3-A-1} can be treated
as  harmless perturbations  in the high frequencies. We first make some  analysis for  Green's matrix
$\mathcal{G}(x,t)$ of the following linearized system without outer forces, namely
\begin{equation}\label{gama3-A-2-11}
\left\{
\begin{aligned}{}
&\p_ta+\Delta v=0,\\
&\p_tv-\nu\Delta v-\Delta  a=0.
\end{aligned}
\right.
\end{equation}
\begin{Lemma}\label{Lem:Greenmatrix}
Let $\mathcal {G}$ be  Green's matrix of  system \eqref{gama3-A-2-11}.
Then we have the following explicit expression for $\widehat{\mathcal {G}}$:

\noindent{(i)}\, when $\nu^{2}\neq4$,  \begin{equation}\label{equ:Green matrix1}
\widehat{\mathcal{G}}(\xi,t)=
\begin{pmatrix}
\frac{\lambda_{+}e^{\lambda_{-}t}-\lambda_{-}e^{\lambda_{+}t}}{\lambda_{+}-\lambda_{-}}&
\frac{e^{\lambda_{+}t}-e^{\lambda_{-}t}}{\lambda_{+}-\lambda_{-}}|\xi|^{2}\\
-\frac{e^{\lambda_{+}t}-e^{\lambda_{-}t}}{\lambda_{+}-\lambda_{-}}|\xi|^{2}&
\frac{\lambda_{+}e^{\lambda_{+}t}-\lambda_{-}e^{\lambda_{-}t}}{\lambda_{+}-\lambda_{-}}
\end{pmatrix},
\end{equation}
where $$\lambda_{\pm}=-\frac{\upsilon}{2}|\xi|^{2}\pm|\xi|^{2}\sqrt{\frac{\upsilon^{2}}{4}-1}.$$

\noindent{(ii)}\, when $\nu^{2}=4$,
\begin{equation}\label{equ:Green matrix2}
\widehat{\mathcal{G}}(\xi,t)=e^{-\frac{\upsilon}{2}|\xi|^{2}t}
\begin{pmatrix}
1+\frac{\upsilon}{2}t|\xi|^{2}&|\xi|^{2}t-\frac{2}{\upsilon}\\
-|\xi|^{2}t&1-\frac{\upsilon}{2}|\xi|^{2}t
\end{pmatrix},
\end{equation}where $$\lambda_{\pm}=-\frac{\upsilon}{2}|\xi|^{2}.$$
Moreover, there is  a positive
constant $\theta$  such that for $\xi\in \mathbb{R}^{N}$,
\begin{equation}\label{3.41}
|\widehat{\mathcal{G}}(\xi,t)|\leq Ce^{-\theta|\xi|^{2}t}.
\end{equation}
\end{Lemma}

\noindent{\bf Proof.}\  Taking Fourier transforms to  system \eqref{gama3-A-2-11} yields that
\begin{align} \label{1.2-A}
\left\{
\begin{aligned}
&\partial_{t}\hat{a}-|\xi|^{2}\hat{v}=0,\\
&\partial_{t}\hat{v}+\nu|\xi|^{2}\hat{v}+|\xi|^{2}\hat{a}=0.
\end{aligned}
\right.
\end{align}
Differentiating with respect to the time variable $t$  of \eqref{1.2-A}$_{2}$ yields that
\begin{equation}\label{1.3-A}
\hat{v}_{tt}+\nu|\xi|^{2}\hat{v}_{t}+|\xi|^{2}\hat{a}_{t}=0.
\end{equation}
Plugging \eqref{1.2-A}$_{1}$ into \eqref{1.3-A} gives rise to
\begin{align} \label{1.4-A}
\left\{
\begin{aligned}
&\hat{v}_{tt}+\nu|\xi|^{2}\hat{v}_{t}+|\xi|^{2}\hat{v}=0,\\
&\hat{v}(\xi,0)=\hat{v}_{0}(\xi),\hat{v}_{t}(\xi,0)=-\nu|\xi|^{2}\hat{v}_{0}-|\xi|^{2}\hat{a}_{0}.
\end{aligned}
\right.
\end{align}
It is easy to check that \begin{equation}\label{1.11-A}
\lambda_{\pm}=-\frac{\nu}{2}|\xi|^{2}\pm|\xi|^{2}\sqrt{\frac{\nu^{2}}{4}-1}
\end{equation}
are two roots of the corresponding characteristic equation of  \eqref{1.4-A}.

\noindent{ Case 1.}\,  For $\nu^{2}\neq4$. We assume that the solution of \eqref{1.4-A} has the following  form
\begin{equation}\label{1.5-A}
\hat{v}(\xi,t)=A(\xi)e^{\lambda_{-}(\xi)t}+B(\xi)e^{\lambda_{+}(\xi)t}.
\end{equation}
Using the initial conditions, we get
\begin{equation}\label{1.6-A}
\begin{split}
&A(\xi)=\frac{-(\lambda_{+}+\nu|\xi|^{2})\hat{v}_{0}-|\xi|^{2}\hat{a}_{0}}
{\lambda_{-}-\lambda_{+}},\\
&B(\xi)=\frac{(\lambda_{-}+\nu|\xi|^{2})\hat{v}_{0}+|\xi|^{2}\hat{a}_{0}}
{\lambda_{-}-\lambda_{+}},
\end{split}
\end{equation}
which implies that
\begin{equation}\label{1.7-A}
\hat{v}(\xi,t)=\frac{e^{\lambda_{+}t}-e^{\lambda_{-}t}}{\lambda_{-}-\lambda_{+}}|\xi|^{2}\hat{a}_{0}(\xi)+
\frac{\lambda_{-}e^{\lambda_{-}t}-\lambda_{+}e^{\lambda_{+}t}}{\lambda_{-}-\lambda_{+}}\hat{v}_{0}(\xi).
\end{equation}
On the other hand, from \eqref{1.2-A}$_{1}$, we obtain
\begin{equation}\label{1.8-A}
\hat{a}(\xi,t)=\hat{a}(\xi,0)+|\xi|^{2}\int_{0}^{t}\hat{v}(\xi,\tau)d\tau.
\end{equation}
Putting \eqref{1.7-A} into the above equality and using the following relations
\begin{equation}\label{1.9-A}
\lambda_{\pm}+\lambda_{\mp}=-\nu|\xi|^{2},\lambda_{+}\lambda_{-}=|\xi|^{4},
\end{equation}
we finally get
\begin{equation}\label{1.10-A}
\hat{a}(\xi,t)=\frac{\lambda_{-}e^{\lambda_{+}t}-\lambda_{+}e^{\lambda_{-}t}}{\lambda_{-}-\lambda_{+}}\hat{a}_{0}(\xi)+
\frac{e^{\lambda_{-}t}-e^{\lambda_{+}t}}{\lambda_{-}-\lambda_{+}}|\xi|^{2}\hat{v}_{0}(\xi).
\end{equation}
Thus, we get  an explicit derivation of the Fourier transform $\widehat{\mathcal{G}}(\xi,t)$ of  Green's
matrix corresponding  linearized system  \eqref{gama3-A-2-11} as follows
\begin{equation}\label{equ:Green matrix11}
\widehat{\mathcal{G}}(\xi,t)=
\begin{pmatrix}
\frac{\lambda_{+}e^{\lambda_{-}t}-\lambda_{-}e^{\lambda_{+}t}}{\lambda_{+}-\lambda_{-}}&
\frac{e^{\lambda_{+}t}-e^{\lambda_{-}t}}{\lambda_{+}-\lambda_{-}}|\xi|^{2}\\
-\frac{e^{\lambda_{+}t}-e^{\lambda_{-}t}}{\lambda_{+}-\lambda_{-}}|\xi|^{2}&
\frac{\lambda_{+}e^{\lambda_{+}t}-\lambda_{-}e^{\lambda_{-}t}}{\lambda_{+}-\lambda_{-}}
\end{pmatrix}.
\end{equation}
For $\nu^{2}<4$, we denote $h=\sqrt{1-\frac{\nu^{2}}{4}}$, thus $h>0$ and $\lambda_{\pm}=-\frac{\nu}{2}|\xi|^{2}\pm ih|\xi|^{2}$. Employing Euler's formula, we have
\begin{equation*}
\frac{e^{\lambda_{+}t}-e^{\lambda_{-}t}}{\lambda_{+}-\lambda_{-}}=\frac{\sin(h|\xi|^{2}t)}{h|\xi|^{2}}
e^{-\frac{\nu}{2}|\xi|^{2}t},
\end{equation*}
\begin{equation*}
\frac{\lambda_{+}e^{\lambda_{-}t}-\lambda_{-}e^{\lambda_{+}t}}{\lambda_{+}-\lambda_{-}}=
\big[\cos(h|\xi|^{2}t)+\frac{\nu}{2h}\sin(h|\xi|^{2}t)
\big]e^{-\frac{\nu}{2}|\xi|^{2}t},
\end{equation*}
\begin{equation*}
\frac{\lambda_{+}e^{\lambda_{+}t}-\lambda_{-}e^{\lambda_{-}t}}{\lambda_{+}-\lambda_{-}}=
\big[\cos(h|\xi|^{2}t)-\frac{\nu}{2h}\sin(h|\xi|^{2}t)
\big]e^{-\frac{\nu}{2}|\xi|^{2}t}.
\end{equation*}
Thus, we can easily  verify that there is  a positive
constant $\theta$  such that for $\xi\in \mathbb{R}^{N}$,
\begin{equation}\label{3.41-11}
|\widehat{\mathcal{G}}(\xi,t)|\leq Ce^{-\theta|\xi|^{2}t}.
\end{equation}
For $\nu^{2}>4$, we denote $\upsilon'=\sqrt{\frac{\nu^{2}}{4}-1}$, thus $\upsilon'>0$ and $\lambda_{\pm}=-\frac{\nu}{2}|\xi|^{2}\pm \upsilon'|\xi|^{2}$. Then,
\begin{equation*}
\frac{e^{\lambda_{+}t}-e^{\lambda_{-}t}}{\lambda_{+}-\lambda_{-}}=\frac{\sinh(\upsilon'|\xi|^{2}t)}{\upsilon'|\xi|^{2}}
e^{-\frac{\nu}{2}|\xi|^{2}t},
\end{equation*}
\begin{equation*}
\frac{\lambda_{+}e^{\lambda_{-}t}-\lambda_{-}e^{\lambda_{+}t}}{\lambda_{+}-\lambda_{-}}=
[\frac{\nu}{2\upsilon'}\sinh(|\xi|^{2}\upsilon't)+\cosh(|\xi|^{2}\upsilon't)]
e^{-\frac{\nu}{2}|\xi|^{2}t},
\end{equation*}
\begin{equation*}
\frac{\lambda_{+}e^{\lambda_{+}t}-\lambda_{-}e^{\lambda_{-}t}}{\lambda_{+}-\lambda_{-}}=
[-\frac{\nu}{2\upsilon'}\sinh(|\xi|^{2}\upsilon't)+\cosh(|\xi|^{2}\upsilon't)
]e^{-\frac{\nu}{2}|\xi|^{2}t},
\end{equation*}
which  implies that \eqref{3.41} holds.

\noindent{ Case 2.}\,  For $\nu^{2}=4$.
We assume that the solution of \eqref{1.4-A} has the  following form
\begin{equation}\label{1.12-A}
\hat{v}(\xi,t)=(A(\xi)+B(\xi)t)e^{-\frac{\nu}{2}|\xi|^{2}t}.
\end{equation}
It follows from the initial conditions that
\begin{equation}\label{1.13}
\begin{split}
&A(\xi)=\hat{v}_{0}(\xi),\\
&B(\xi)=-\nu|\xi|^{2}\hat{v}_{0}-|\xi|^{2}\hat{a}_{0}-\frac{\nu}{2}|\xi|^{2}\hat{v}_{0},
\end{split}
\end{equation}
which give rise to
\begin{equation}\label{1.14}
\hat{v}(\xi,t)=-|\xi|^{2}te^{-\frac{\nu}{2}|\xi|^{2}t}\hat{a}_{0}+(1-\frac{\nu}{2}|\xi|^{2}t)
e^{-\frac{\nu}{2}|\xi|^{2}t}\hat{v}_{0}.
\end{equation}
Furthermore, from  \eqref{1.2-A}$_{1}$, we obtain
\begin{equation}\label{1.15-A}
\hat{a}(\xi,t)=\hat{a}(\xi,0)+|\xi|^{2}\int_{0}^{t}\hat{v}(\xi,\tau)d\tau.
\end{equation}
Plugging \eqref{1.7-A} into the above equality, we finally  conclude that
\begin{equation}\label{1.16-A}
\hat{a}(\xi,t)=
(1+\frac{\nu}{2}t|\xi|^{2})e^{-\frac{\nu}{2}|\xi|^{2}t}
\hat{a}_{0}(\xi)+(|\xi|^{2}t-\frac{2}{\nu})
e^{-\frac{\nu}{2}|\xi|^{2}t}\hat{v}_{0}(\xi).
\end{equation}
Then
\begin{equation}\label{equ:Green matrix22}
\widehat{\mathcal{G}}(\xi,t)=e^{-\frac{\nu}{2}|\xi|^{2}t}
\begin{pmatrix}
1+\frac{\nu}{2}t|\xi|^{2}&|\xi|^{2}t-\frac{2}{\nu}\\
-|\xi|^{2}t&1-\frac{\nu}{2}|\xi|^{2}t
\end{pmatrix}.
\end{equation}
By a simple  computation,  we also conclude that \eqref{3.41} holds.The proof of Lemma \ref{Lem:Greenmatrix} is complete.\endproof

With Lemma \ref{Lem:Greenmatrix} at hand,  we can  exploit the following parabolic properties of system \eqref{gama3-A-2-1} in  $L^{p}$-framework in the high frequencies.
\begin{Lemma}\label{lemma3.4}
Let $s\in \mathbb{R}, p\in[1,+\infty]$, $1\leq \rho_{2}\leq \rho_{1}\leq\infty$, $\nu>0$,   and $T\in(0,+\infty]$. Suppose that $(a_{0},v_{0})\in(B_{p,1}^{s})^{2}$ and $(F,G)\in(\widetilde{L}_{T}^{\rho_{1}}(B_{p,1}^{s-2+\frac{2}{\rho_{1}}}))^{2}$. Then,   system \eqref{gama3-A-2-1}
has a unique solution $(a,v)$ in $(\tilde{C}_{T}(B_{p,1}^{s})\cap\widetilde{L}_{T}^{\rho_{1}}(B_{p,1}^{s-2+\frac{2}{\rho_{1}}}))^{2}$. Moreover,  there exists $C>0$ depending only on $\nu,\rho,\rho_{1}$ such that
\begin{equation}\label{eq:3.47}
\|(a,v)\|_{\tilde{L}_{T}^{\rho_1}(B_{p,1}^{s+\frac{2}{\rho_1}})}^h\leq C\Big(
\|(a_{0},v_{0})\|_{B_{p,1}^{s}}^h+\|(F,G)\|_{\widetilde{L}_{T}^{\rho_{2}}(B_{p,1}^{s-2+\frac{2}{\rho_{2}}})}^h\Big).
\end{equation}
\end{Lemma}

\noindent{\bf Proof.}\, In terms of  Green's matrix and Duhamel's principle, the solution of   system \eqref{gama3-A-2-1} can be expressed as
\begin{align}\label{equ:4.5}
\left(\begin{array}{ll} \!a \!\!\vspace{.15cm}\\ \!v \!\end{array}\right)=
\mathcal{G}(x,t)\ast\left(\begin{array}{ll} \!a_0 \!\vspace{.15cm}\\ \!v_0 \!\end{array}\right)
+\int_0^t\mathcal{G}(x,t-\tau)\ast\left(\begin{array}{ll}\!F\!\vspace{.15cm}\\
\!G\!\end{array}\right)d\tau.
\end{align}
Applying  homogeneous frequency localization operators
$\dot{\Delta}_j$ on both sides of \eqref{equ:4.5}, we get
\begin{align*}
\left(\begin{array}{ll} \!\dot{\Delta}_ja \!\!\vspace{.15cm}\\ \!\dot{\Delta}_jv \!\end{array}\right)=
\mathcal{G}(x,t)\ast\left(\begin{array}{ll} \!\dot{\Delta}_ja_0 \!\vspace{.15cm}\\ \!\dot{\Delta}_jv_0 \!\end{array}\right)
+\int_0^t\mathcal{G}(x,t-\tau)\ast\left(\begin{array}{ll}\!\dot{\Delta}_jF\!\vspace{.15cm}\\
\!\dot{\Delta}_jG\!\end{array}\right)d\tau.
\end{align*}
From  \eqref{3.41} and Young's inequality, we infer that
\begin{align*}
\|\dot{\Delta}_ja(t)\|_{L^p}&+\|\dot{\Delta}_jv(t)\|_{L^p}\le Ce^{-c2^{2j} t}\bigl(\|\dot{\Delta}_ja_0\|_{L^p}+\|\dot{\Delta}_jv_0\|_{L^p}\bigr)\\
&+C\int_0^te^{-c2^{2j} (t-\tau)}
\Big(\|\dot{\Delta}_jF(\tau)\|_{L^p}+\|\dot{\Delta}_jG(\tau)\|_{L^p}\Big)d\tau.
\end{align*}
Taking $L^{\rho_1}$ norm with respect to $t$, and using convolution inequality with $1+\frac{1}{\rho_1}=\frac{1}{\rho_2}+\frac{1}{\rho}$  yield  that
\begin{align}\label{equ:4.6}\begin{split}
&\|\dot{\Delta}_ja\|_{L^{\rho_1}_tL^p}+\|\dot{\Delta}_jv\|_{L^{\rho_1}_tL^p}\\&\le C2^{-\f {2j}{\rho_1}}
\Big(\|\dot{\Delta}_ja_0\|_{L^p}+\|\dot{\Delta}_jv_0\|_{L^p}
+2^{-2j+\f {2j}{\rho_2}}\bigl(\|\dot{\Delta}_jF\|_{L^{\rho_2}_tL^p}+\|\dot{\Delta}_jG\|_{L^{\rho_2}_tL^p}\bigr)\Big).
\end{split}\end{align}
Multiplying $2^{js}$ on both sides of \eqref{equ:4.6}, and summing up over  $j\geq j_0$, where $j_0\in \mathbb{Z}$,  we conclude that  \eqref{eq:3.47}. The proof of Lemma \ref{lemma3.4} is complete.\endproof

\noindent{\bf Proof of Proposition \ref{high part}.}\, Applying  Lemma \ref{lemma3.4} to system \eqref{gama3-A-1}, we have
\begin{equation}
\begin{split}
\label{high part1-A}
&\|(\Delta c^{+},\,\textrm{div}\cQ u^{+},\,\Delta c^{-},\,\textrm{div}\cQ u^{-})\|_{\wt L^\infty_t(\dot B^{\frac{N}{p}-2}_{p,1})}^h
+\|(\Delta c^{+},\,\textrm{div}\cQ u^{+},\,\Delta c^{-},\,\textrm{div}\cQ u^{-})\|_{L^1_t(\dot B^{\frac{N}{p}}_{p,1})}^h
\\&\lesssim \|(\Delta c^{+}_0,\,\textrm{div}\cQ u^{+}_0,\,\Delta c^{-}_0,\,\textrm{div}\cQ u^{-}_0)\|_{\dot B^{\frac{N}{p}-2}_{p,1}}^h
+\|(\textrm{div}\cQ H_{2},\,\textrm{div}\cQ H_{4})\|_{ L^1_t(\dot B^{\frac{N}{p}-2}_{p,1})}^h\\&\quad+\|(\Delta c^{+},\,\Delta c^{-})\|_{ L^1_t(\dot B^{\frac{N}{p}-2}_{p,1})}^h+\|(\Delta H_{1},\,\Delta H_{3})\|_{ L^1_t(\dot B^{\frac{N}{p}-2}_{p,1})}^h.
\end{split}
\end{equation}
Thanks to the high-frequency cutoff, we get
$$\|(\Delta c^{+},\,\Delta c^{-})\|_{ L^1_t(\dot B^{\frac{N}{p}-2}_{p,1})}^h\lesssim2^{-2j}\|(\Delta c^{+},\,\Delta c^{-})\|^{h}_{L_{t}^{1}(\dot{B}_{p,1}^{\frac{N}{p}})}\lesssim2^{-2j_0}\|(\Delta c^{+},\,\Delta c^{-})\|^{h}_{L_{t}^{1}(\dot{B}_{p,1}^{\frac{N}{p}})},\quad \hbox{for}\quad j\geq j_0.$$
Hence, taking $j_0$ large enough, it follows from \eqref{high part1-A}  that
\begin{equation}
\begin{split}
\label{high part11-A}
&\|(\Delta c^{+},\,\textrm{div}\cQ u^{+},\,\Delta c^{-},\,\textrm{div}\cQ u^{-})\|_{\wt L^\infty_t(\dot B^{\frac{N}{p}-2}_{p,1})}^h
+\|(\Delta c^{+},\,\textrm{div}\cQ u^{+},\,\Delta c^{-},\,\textrm{div}\cQ u^{-})\|_{L^1_t(\dot B^{\frac{N}{p}}_{p,1})}^h
\\&\lesssim \|(\Delta c^{+}_0,\,\textrm{div}\cQ u^{+}_0,\,\Delta c^{-}_0,\,\textrm{div}\cQ u^{-}_0)\|_{\dot B^{\frac{N}{p}-2}_{p,1}}^h
+\|(\cQ H_{2},\,\cQ H_{4})\|_{ L^1_t(\dot B^{\frac{N}{p}-1}_{p,1})}^h\\&\quad+\|( H_{1},\, H_{3})\|_{ L^1_t(\dot B^{\frac{N}{p}}_{p,1})}^h.
\end{split}
\end{equation}
Noticing that $u^{\pm}=\cP u^{\pm}+\cQ u^{\pm}$,   and  then combining  with \eqref{eq:HFo}  and \eqref{high part11-A} yields that
\begin{equation}
\begin{split}
\label{high part111-A}
\|(u^{+},\,u^{-})\|_{\wt L^\infty_t(\dot B^{\frac{N}{p}-1}_{p,1})}^h
&+\|(u^{+},\,u^{-})\|_{L^1_t(\dot B^{\frac{N}{p}+1}_{p,1})}^h
+\|(c^{+},\,c^{-})\|_{\wt L^\infty_t(\dot B^{\frac{N}{p}}_{p,1})}^h
+\|(c^{+},\,c^{-})\|_{L^1_t(\dot B^{\frac{N}{p}}_{p,1})}^h
\\&\lesssim X(0)+\|( H_{1},\, H_{3})\|_{ L^1_t(\dot B^{\frac{N}{p}}_{p,1})}^h+\|(H_{2},\,H_{4})\|_{ L^1_t(\dot B^{\frac{N}{p}-1}_{p,1})}^h,
\end{split}
\end{equation}where we have used  $\cP$ and $\cQ$ are  continuous on  $\dot B^s_{p,1}$ (being
$0$ order multipliers).

In what follows,  we shall bound the nonlinear terms on the right-side of \eqref{high part111-A}. First, due to
 $H_{1}=-\textrm{div}(c^{+}u^{+})=-c^{+}\textrm{div}u^{+}-u^{+}\cdot\nabla c^{+},\, H_{3}=-\textrm{div}(c^{-}u^{-})=-c^{-}\textrm{div}u^{-}-u^{-}\cdot\nabla c^{-},$
employing Proposition \ref{p26} yields that
\begin{equation*}
\begin{split}\label{eq:all}
 &\|( H_{1},\, H_{3})\|_{ L^1_t(\dot B^{\frac{N}{p}}_{p,1})}^h
 \\&\lesssim\int_{0}^{t}\|\nabla u^{+}\|_{\dot B^{\frac{N}{p}}_{p,1}}
\|c^{+}\|_{\dot B^{\frac{N}{p}}_{p,1}}d\tau+\int_{0}^{t}\|\nabla u^{-}\|_{\dot B^{\frac{N}{p}}_{p,1}}
\|c^{-}\|_{\dot B^{\frac{N}{p}}_{p,1}}d\tau\\&\quad+\int_{0}^{t}\|\nabla c^{+}\|_{\dot B^{\frac{N}{p}}_{p,1}}
\|u^{+}\|_{\dot B^{\frac{N}{p}}_{p,1}}d\tau+\int_{0}^{t}\|\nabla c^{-}\|_{\dot B^{\frac{N}{p}}_{p,1}}
\|u^{-}\|_{\dot B^{\frac{N}{p}}_{p,1}}d\tau
\\&\lesssim\|c^{+}\|_{L^\infty_t(\dot B^{\frac Np}_{p,1})}\|u^{+}\|_{L^1_t(\dot B^{\frac Np+1}_{p,1})}+\|c^{-}\|_{L^\infty_t(\dot B^{\frac Np}_{p,1})}\|u^{-}\|_{L^1_t(\dot B^{\frac Np+1}_{p,1})}\\&\quad+\|c^{+}\|_{L^2_t(\dot B^{\frac Np+1}_{p,1})}\|u^{+}\|_{L^2_t(\dot B^{\frac Np}_{p,1})}+ \|c^{-}\|_{L^2_t(\dot B^{\frac Np+1}_{p,1})}\|u^{-}\|_{L^2_t(\dot B^{\frac Np}_{p,1})}
 \\&\lesssim X^{2}(t).
\end{split}
\end{equation*}
For the term $\|(H_2,H_4)\|_{L^1_t(\dot B^{\frac Np-1}_{p,1})}^h$, omitting some positive constants,  it suffices to bound $\|(H_2,H_4)\|_{L^1_t(\dot B^{\frac Np-1}_{p,1})}^h$.
For the  term $H_{2}^{i}$,  thanks to  Proposition \ref{p26} and Proposition \ref{p27}(i),  we get
\begin{equation*}
\begin{split}
\|g_{+}(c^{+},c^{-})\partial_{i} c^{+}\|_{L^1_t(\dot B^{\frac Np-1}_{p,1})}^h
&\lesssim\|g_{+}(c^{+},c^{-})\|_{L^2_t(\dot{B}^{\f{N}{p}}_{p,1})}\|\nabla c^{+}\|_{L^2_{t}
(\dot{B}^{\f{N}{p}-1}_{p,1})}
\\&\lesssim\Big(1+X^{2}(t)\Big)^{[\frac{N}{p}]+1}\| (c^{+},c^{-})\|_{\tilde{L}^2_{t}
(\dot{B}^{\f{N}{p}}_{p,1})}^{2}
\\&\lesssim \Big(1+X^{2}(t)\Big)^{[\frac{N}{p}]+1}X^{2}(t),
\end{split}
\end{equation*}
\begin{equation*}
\begin{split}
\|u^{+}\cdot\nabla u_{i}^{+}\|_{L^1_t(\dot B^{\frac Np-1}_{p,1})}^h&\lesssim\|u_{+}\|_{L^2_t(\dot{B}^{\f{N}{p}}_{p,1})}\|\nabla u^{+}\|_{L^2_{t}
(\dot{B}^{\f{N}{p}-1}_{p,1})}
\\&\lesssim\| u^{+}\|_{L^2_{t}
(\dot{B}^{\f{N}{p}}_{p,1})}^{2}
\\&\lesssim X^{2}(t),
\end{split}
\end{equation*}
\begin{equation*}
\begin{split}
&\|h_{+}(c^{+},c^{-})\partial_{j}c^{+}\partial_{j}u^{+}_{i}\|_{L^1_t(\dot B^{\frac Np-1}_{p,1})}^h
\\&\lesssim\|h_{+}(c^{+},c^{-})\partial_{j}c^{+}\|_{L^\infty_t(\dot{B}^{\f{N}{p}-1}_{p,1})}\|\partial_{j}u^{+}_{i}\|_{L^1_{t}
(\dot{B}^{\f{N}{p}}_{p,1})}
\\&\lesssim\|h_{+}(c^{+},c^{-})\|_{L^\infty_{t}
(\dot{B}^{\f{N}{p}}_{p,1})}
\|\nabla c^{+}\|_{L^\infty_{t}
(\dot{B}^{\f{N}{p}-1}_{p,1})}
\|\nabla u^{+}\|_{L^1_t(\dot{B}^{\f{N}{p}}_{p,1})}
\\&\lesssim\Big(1+ \big(1+\|(c^{+},c^{-})\|^{2}_{\tilde{L}^\infty_t(\dot{B}^{\f{N}{p}}_{p,1})}\big)^{[\frac{N}{p}+1]}
\|(c^{+},c^{-})\|_{\tilde{L}^\infty_t(\dot{B}^{\f{N}{p}}_{p,1})}\Big)
\| c^{+}\|_{L^\infty_{t}
(\dot{B}^{\f{N}{p}}_{p,1})}
\|u^{+}\|_{L^1_{t}
(\dot{B}^{\f{N}{p}+1}_{p,1})}
\\&\lesssim X^{2}(t)+\Big(1+X^{2}(t)\Big)^{[\frac{N}{p}]+1}X^{3}(t),
\end{split}
\end{equation*}
and
\begin{equation*}
\begin{split}
\|l_{+}(c^{+},c^{-})\partial_{j}^{2}u_{i}^{+}\|_{L^1_t(\dot B^{\frac Np-1}_{p,1})}^h
&\lesssim\|l_{+}(c^{+},c^{-})\|_{L^\infty_{t}
(\dot{B}^{\f{N}{p}}_{p,1})}
\|\partial_{j}^{2}u_{i}^{+}\|_{L^1_{t}
(\dot{B}^{\f{N}{p}-1}_{p,1})}
\\&\lesssim\|(c^{+},c^{-})\|_{L^\infty_{t}
(\dot{B}^{\f{N}{p}}_{p,1})}
\|u^{+}\|_{L^1_{t}
(\dot{B}^{\f{N}{p}+1}_{p,1})}
\\&\lesssim \Big(1+X^{2}(t)\Big)^{[\frac{N}{p}]+1}X^{2}(t).
\end{split}
\end{equation*}
Similarly,  we also obtain  the corresponding estimates of other terms
$(\mu^{+}k_{+}(c^{+},c^{-})\partial_{j}c^{-}\partial_{j}u^{+}_{i})^{h},$\\
$(\mu^{+}h_{+}(c^{+},c^{-})\partial_{j}c^{+}\partial_{i}u^{+}_{j})^{h},$
$(\mu^{+}k_{+}(c^{+},c^{-})\partial_{j}c^{-}\partial_{i}u^{+}_{j})^{h},
(\lambda^{+}h_{+}(c^{+},c^{-})\partial_{i}c^{+}\partial_{j}u^{+}_{j})^{h},$\\
$(\lambda^{+}k_{+}(c^{+},c^{-})\partial_{i}c^{-}\partial_{j}u^{+}_{j})^{h}$, $((\mu^{+}+\lambda^{+})l_{+}
(c^{+},c^{-})\partial_{i}\partial_{j}u^{+}_{j})^{h}$. Here, we omit the details.\\
Hence
\begin{equation*}
\|H_{2}\|_{L^{1}([0,t];\dot{B}^{\f{N}{p}-1}_{2,1})}^{h}
\lesssim X^{2}(t)+\Big(1+X^{2}(t)\Big)^{[\frac{N}{p}]+1}\Big(X^{2}(t)+X^{3}(t)\Big).
\end{equation*}
The term  $H_{4}$ may be treated along the same lines,  and we omit it.

Plugging  all the above nonlinear estimates of the terms  $H_{1}$--$H_{4}$ into  \eqref{high part111-A}, we conclude that \eqref{high part1}. This completes the proof of Proposition \ref{high part}.\endproof
\subsection{The Unique Global Solvability}\ \ \ \ \
Combining  with Propositions \ref{low part}  and \ref{high part},  we have the following \emph{a priori}  estimates in all frequencies.
\begin{Proposition}\label{ab} Let $T\geq0$,
$N\geq2$, $p$ satisfy $2\leq p\leq\min(4,\frac{2N}{N-2})$ and, additionally, $ p\neq 4$ if $ N=2$.
Assume that $(c^{+},\, u^{+},\, c^{-},\, u^{-})$ is a solution to the  Cauchy problem \eqref{equ:CTFS2}-\eqref{equ:CTFS3} on $[0,T]\times \mathbb{R}^{N}$,  then
\begin{equation}\label{eq:5.1}X(t)\lesssim\Big(X(0)+\big(1+X^{2}(t)\big)^{[\frac{N}{p}]+1}\big(X^{2}(t)+X^{3}(t)\big)\Big)  \quad \text{for all} \quad t\in[0,T].\end{equation}
\end{Proposition}

From  \eqref{eq:5.1}, it is not difficult to work out a fixed point
argument as in  \cite{XC},  we  finally conclude  the  global well-posedness  of the Cauchy problem \eqref{equ:CTFS2}-\eqref{equ:CTFS3}. This completes the proof of Theorem \ref{th:main1}.

\section{The proof of Theorem \ref{th:decay}}
\ \ \ \ \ In this section,  our central task is to prove Theorem \ref{th:decay} taking for granted the global-in-time existence result of Theorem \ref{th:main1}.
The proof is divided into two parts, according to the two terms of the time-weighted functional $D(t)$ (see \eqref{1.9}). In what follows, we shall use frequently
elementary  fact  that the global solution $(c^{+},u^{+},c^{-},u^{-})$ provided by Theorem \ref{th:main1} fulfills
\begin{equation}\label{99111B}
\|(c^{+},c^{-})\|_{\wt L^\infty_t(\dot B^{\frac{N}{p}}_{p,1})}\leq c\ll1 \quad \text{for all}\quad t\geq0.
\end{equation}
\subsection{In the low frequencies}\ \ \ \ \
Denoting by $A(D)$ the semi-group associated to system \eqref{equ:CTFS2} with $H_1\equiv H_2\equiv H_3\equiv H_4\equiv0$ for $U=(c^{+},u^{+},c^{-},u^{-})$,      and using Parseval's equality, \eqref{gama8} and the definition of $\ddq$,   we get for all $q\leq q_{0}$
\begin{align*}\|e^{tA(D)}\ddq U\|_{L^2}&\lesssim  e^{-c_{0}2^{2q}t}\|\ddq U\|_{L^2}.
\end{align*}
Hence,   multiplying by $t^{\frac{s+s_0}{2}}2^{qs}$ and summing up on $q$, we readily have
\begin{equation}
\begin{split}
\label{low.2}
t^{\frac{s+s_0}{2}}\sum_{q\leq q_0}2^{qs}\|e^{tA(D)}\ddq U\|_{L^2}
&\lesssim
\sum_{q\leq q_0}2^{qs}e^{-c_{0}2^{2q}t}\|\ddq U\|_{L^2}t^{\frac{s+s_0}{2}}\\
&\lesssim
\sum_{q\leq q_0}2^{q(s+s_0)}e^{-c_{0}2^{2q}t}\|\ddq U\|_{L^2}2^{q(-s_0)}t^{\frac{s+s_0}{2}}\\
&\lesssim
\|U\|_{\dot B^{-s_0}_{2,\infty}}^\ell\sum_{q\leq q_0}2^{q(s+s_0)}e^{-c_{0}2^{2q}t}t^{\frac{s+s_0}{2}}\\
&\lesssim
\|U\|_{\dot B^{-s_0}_{2,\infty}}^\ell\sum_{q\leq q_0}2^{q(s+s_0)}e^{-c_{0}2^{2q}t}
t^{\frac{s+s_0}{2}}.
\end{split}
\end{equation}
As for any $\sigma>0$ there  exists a constant $C_\sigma$ so that
\begin{equation}\label{low.3}
\sup_{t\geq0}\sum_{q\in\mathbb{Z}}t^{\frac\sigma2}2^{q\sigma}e^{-c_{0}2^{2q}t}\leq C_\sigma.
\end{equation}
We get from \eqref{low.2} and \eqref{low.3} that for $s+s_0>0,$
$$
\sup_{t\geq0}\, t^{\frac{s+s_0}{2}}\|e^{tA(D)}U\|_{\dot B^s_{2,1}}^\ell
\lesssim\|U\|_{\dot B^{-\frac{N}{2}}_{2,\infty}}^\ell.
$$
Furthermore, it is  obvious that  for $s+s_0>0,$
$$
\|e^{tA(D)}U\|_{\dot B^s_{2,1}}^\ell
\lesssim \|U\|_{\dot B^{-s_0}_{2,\infty}}^\ell\sum_{q\leq q_0}2^{q(s+s_0)}\lesssim\|U\|_{\dot B^{-s_0}_{2,\infty}}^\ell.
$$
Hence, setting $\langle t\rangle\eqdefa(1+t)$, we get
\begin{equation}
\label{U}
\sup_{t\geq0}\, \langle t\rangle^{\frac{s+s_0}{2}}\|e^{tA(D)}U\|_{\dot B^s_{2,1}}^\ell
\lesssim\|U\|_{\dot B^{-s_0}_{2,\infty}}^\ell.
\end{equation}
Thus, from \eqref{U} and Duhamel's formula, we have
\begin{equation}\label{eq:5.4}
\begin{split}
&\big\|\big(c^{+},u^{+},c^{-},u^{-}\big)\big\|_{\dot B^s_{2,1}}^\ell\\
&\quad\lesssim \sup_{t\geq0}\, \langle t\rangle^{-\frac{s+s_0}{2}}\big\|\big(c^{+}_0,u^{+}_0,c^{-}_0,u^{-}_0\big)\big\|_{\dot B^{-s_0}_{2,\infty}}^\ell\\&\qquad+\int_{0}^{t}\langle t-\tau\rangle^{-\frac{s+s_0}{2}}\big\|(H_{1},H_{2},H_{3},H_{4})(\tau)
\big\|_{\dot B^{-s_0}_{2,\infty}}^\ell d\tau.
\end{split}
\end{equation}
We claim that for all $s\in[\varepsilon-s_0,\frac{N}{2}+1]$ and $t\geq0$,  then
\begin{equation}\label{s1234low1}\begin{split}
\int_{0}^{t}\langle t-\tau\rangle^{-\frac{s+s_0}{2}}\big\|\big(H_{1},H_{2},H_{3},H_{4}\big)(\tau)
\big\|_{\dot B^{-s_0}_{2,\infty}}^\ell d\tau
\lesssim\langle t\rangle^{-\frac{s+s_0}{2}}
\Big(X^2(t)+D^2(t)\Big),
\end{split}
\end{equation}
where $X(t)$ and $D(t)$ have been defined in \eqref{1.6} and \eqref{1.9}, respectively.\medbreak

In  order to prove our claim,  we first present the following some important inequalities which will be frequently used in our process later.
\begin{equation}
\begin{split}
\label{6.6}
\big\|\langle\tau\rangle^{\alpha}\big(\nabla c^{+},u^{+},\nabla c^{-},u^{-}\big)\big\|_{\wt L^\infty_t(\dot B^{\frac Np-1}_{p,1})}^h
&\lesssim\big\|\big(\nabla c^{+},u^{+},\nabla c^{-},u^{-}\big)\big\|_{\wt L^\infty_t(\dot B^{\frac{N}{p}-1}_{p,1})}^h.
\\&\qquad+\big\|\tau^{\alpha}\big(\nabla c^{+},u^{+},\nabla c^{-},u^{-}\big)\big\|_{\wt L^\infty_t(\dot B^{\frac Np-1}_{p,1})}^h
\\&\quad\lesssim X(t)+D(t),
\end{split}
\end{equation}
\begin{equation}\label{low.8}
\begin{split}
&\|(c^{+},u^{+},c^{-},u^{-})\|_{\dot B^{1-\frac{N}{p}}_{p,1}}^{\ell}
\lesssim\|(c^{+},u^{+},c^{-},u^{-})\|_{\dot B^{1-s_0}_{2,1}}^{\ell}
\lesssim\langle \tau\rangle^{-\frac 12}D(\tau),
\end{split}
\end{equation}
\begin{equation}\label{low.10}
\begin{split}
&\|(\nabla c^{+},\nabla u^{+},\nabla c^{-},\nabla u^{-})\|_{\dot B^{\frac{N}{2}-1}_{2,1}}^{\ell}
\lesssim\|(c^{+},u^{+},c^{-},u^{-})\|_{\dot B^{\frac{N}{2}}_{2,1}}^{\ell}
\lesssim\langle \tau\rangle^{-\frac Np}D(\tau).
\end{split}
\end{equation}
\begin{equation}\label{low.11}
\begin{split}
\|(c^{+},c^{-})\|_{\dot B^{1-\frac{N}{p}}_{p,1}}^{h}
&\lesssim\|(c^{+},c^{-})\|_{\dot B^{\frac{N}{p}}_{p,1}}^{h}
\lesssim\langle \tau\rangle^{-\alpha}\Big(D(\tau)+X(\tau)\Big),
\end{split}
\end{equation}
\begin{equation}\label{low.12}
\begin{split}
\|(c^{+},u^{+})\|_{\dot B^{\frac{N}{p}}_{2,1}}^{\ell}
\lesssim\langle \tau\rangle^{-(\frac{3N}{2p}-\frac{N}{4})}D(\tau).
\end{split}
\end{equation}
To bound the term $H_{1}$, using low-high frequency  decomposition    we see that
$$ H_{1}=u^{+}\cdot\nabla c^{+\ell}+u^{+}\cdot\nabla c^{+h}+ c^{+}\, \div (u^{+})^\ell + c^{+}\,\div (u^{+})^h.$$
We now bound  the above  items one by one.  For the term $u^{+}\cdot\nabla c^{+\ell}$, employing \eqref{low.7} yields that
\begin{equation}
\begin{split}\label{s11}
&\int_0^t\langle t-\tau\rangle^{-\frac{s+s_0}{2}}\|(u^{+}\cdot\nabla c^{+\ell})(\tau)\|_{\dot B^{-s_0}_{2,\infty}}^{\ell}\,d\tau\\
&\quad\lesssim
\int_0^t\langle t-\tau\rangle^{-\frac{s+s_0}{2}}
\|u^{+}\|_{\dot B^{1-\frac{N}{p}}_{p,1}}
\|\nabla c^{+\ell}\|_{\dot B^{\frac{N}{2}-1}_{2,1}}\\
&\quad\lesssim\int_0^t\langle t-\tau\rangle^{-\frac{s+s_0}{2}}\Big(\|u^{+\ell}\|_{\dot B^{1-\frac{N}{p}}_{p,1}}+\|u^{+h}\|_{\dot B^{1-\frac{N}{p}}_{p,1}}\Big)\|\nabla c^{+\ell}\|_{\dot B^{\frac{N}{2}-1}_{2,1}}\,d\tau.
\end{split}
\end{equation}
We are going to bound  $\|u^{+h}\|_{\dot B^{1-\frac{N}{p}}_{p,1}}$ in \eqref{s11}. Here, we shall proceed differently depending
on whether  $2\leq p\leq N$ or $N<p<2N$. 
If $2\leq p\leq N$ then $1-\frac{N}{p}\leq\frac{N}{p}-1$, then we conclude, from \eqref{6.6}, that
\begin{equation}\label{low.91}
\begin{split}
\|u^{+h}\|_{\dot B^{1-\frac{N}{p}}_{p,1}}\lesssim\|u^{+h}\|_{\dot B^{\frac{N}{p}-1}_{p,1}}\lesssim\langle \tau\rangle^{-\alpha}\langle \tau\rangle^{\alpha}\|u^{+h}\|_{\dot B^{\frac{N}{p}-1}_{p,1}}\lesssim\langle \tau\rangle^{-\alpha}\Big(D(\tau)+X(\tau)\Big).
\end{split}
\end{equation}
Plugging \eqref{low.91} and \eqref{low.10} into \eqref{s11} implies that
\begin{equation}
\begin{split}\label{s11-A}
&\int_0^t\langle t-\tau\rangle^{-\frac{s+s_0}{2}}\|(u^{+}\cdot\nabla c^{+\ell})(\tau)\|_{\dot B^{-s_0}_{2,\infty}}^{\ell}\,d\tau\\
&\quad\lesssim \Big(D^{2}(t)+X^{2}(t)\Big)\int_0^t\langle t-\tau\rangle^{-\frac{s+s_0}{2}}
\Big(\langle \tau\rangle^{-{(\frac 12+\frac{N}{p})}}
+\langle \tau\rangle^{-{(\frac \alpha2+\frac{N}{p})}}\Big)\,d\tau
\\&\quad\lesssim \langle t\rangle^{-\frac{s+s_0}{2}} \Big(D^{2}(t)+X^{2}(t)\Big),
\end{split}
\end{equation}
where we have  used  Lemma \ref{lemma2.13} and the fact $\frac{s_{0}+s}{2}\leq \min\{\frac{1}{2}+\frac{N}{p},\frac{\alpha}{2}+\frac{N}{p}\}$ with $\alpha=\frac{N}{p}+\frac{1}{2}-\varepsilon$ for all $s\leq\frac{N}{2}+1$.
If $N<p<2N$ then $1-\frac{N}{p}\leq\frac{N}{p}$. Employing interpolation inequality and \eqref{6.6} yields that
\begin{equation}\label{low.92}
\begin{split}
\|u^{+h}\|_{\dot B^{1-\frac{N}{p}}_{p,1}}
&\lesssim\|u^{+h}\|_{\dot B^{\frac{N}{p}}_{p,1}}
\lesssim\Big(\|u^{+h}\|_{\dot B^{\frac{N}{p}-1}_{p,1}}
\|u^{+h}\|_{\dot B^{\frac{N}{p}+1}_{p,1}}\Big)^{\frac{1}{2}}
\\&\lesssim\Big( \langle\tau\rangle^{-\alpha}\langle \tau\rangle^{\alpha}\|u^{+}\|_{\dot B^{\frac{N}{p}-1}_{p,1}}^{h}
\tau^{-\alpha}\tau^{\alpha}\|u^{+}\|_{\dot B^{\frac{N}{p}+1}_{p,1}}^{h}\Big)^{\frac{1}{2}}
\\&\lesssim\langle \tau\rangle^{-\frac \alpha2}\tau^{-\frac{\alpha}{2}}\Big(X(\tau)+D(\tau)\Big).
\end{split}
\end{equation}
Putting \eqref{low.92} and \eqref{low.10} into \eqref{s11}, noticing $\frac{\alpha}{2}=\frac{N}{2p}+\frac{1}{4}-\frac{\varepsilon}{2}<1$, $\frac{s+s_0}{2}\leq\alpha+\frac{N}{p}$ for all $s\leq\frac{N}{2}+1$, and $\alpha+\frac{N}{p}>1$,   and then  using  Lemmas \ref{lemma2.13}-\ref{lemma2.14} give rise to
\begin{equation*}
\begin{split}\label{s11-B}
&\int_0^t\langle t-\tau\rangle^{-\frac{s+s_0}{2}}\|(u^{+}\cdot\nabla c^{+\ell})(\tau)\|_{\dot B^{-s_0}_{2,\infty}}^{\ell}\,d\tau\\
&\quad\lesssim \Big(D^{2}(t)+X^{2}(t)\Big)\int_0^t\langle t-\tau\rangle^{-\frac{s+s_0}{2}}
(\langle \tau\rangle^{-{(\frac 12+\frac{N}{p})}}
+\langle \tau\rangle^{-{(\frac \alpha2+\frac{N}{p})}}\tau^{-\frac{\alpha}{2}})\,d\tau
\\&\quad\lesssim \langle t\rangle^{-\frac{s+s_0}{2}} \Big(D^{2}(t)+X^{2}(t)\Big).
\end{split}
\end{equation*}
For the term  $c^{+}\,\div (u^{+})^\ell$,  using \eqref{low.7}, \eqref{low.8}, \eqref{low.10} and \eqref{low.11} implies that
\begin{equation}
\begin{split}
\label{s12}
&\int_0^t\langle t-\tau\rangle^{-\frac{s+s_0}{2}}\|c^{+}\div (u^{+})^{\ell}\|_{\dot B^{-s_0}_{2,\infty}}^{\ell}\,d\tau\\
&\quad\lesssim
\int_0^t\langle t-\tau\rangle^{-\frac{s+s_0}{2}}
\|c^{+}\|_{\dot B^{1-\frac{N}{p}}_{p,1}}
\|\div (u^{+})^{\ell}\|_{\dot B^{\frac{N}{2}-1}_{2,1}}\\
&\quad\lesssim\int_0^t\langle t-\tau\rangle^{-\frac{s+s_0}{2}}\big(\|c^{+}\|_{\dot B^{1-\frac{N}{p}}_{p,1}}^{\ell}+\|c^{+}\|_{\dot B^{1-\frac{N}{p}}_{p,1}}^h\big)\|\nabla u^{+}\|_{\dot B^{\frac{N}{2}-1}_{2,1}}^{\ell}\,d\tau\\
&\quad\lesssim\Big(D^{2}(t)+X^{2}(t)\Big)\int_0^t\langle t-\tau\rangle^{-\frac{s+s_0}{2}}\langle \tau\rangle^{-\min({\frac 12+\frac Np,\alpha+\frac{N}{p}})}\,d\tau\\
&\quad\lesssim \langle t\rangle^{-\frac{s+s_0}{2}} \Big(D^{2}(t)+X^{2}(t)\Big),
\end{split}
\end{equation}
where we have  used  Lemma \ref{lemma2.13} and the fact $\frac{s_{0}+s}{2}\leq \min\{\frac{1}{2}+\frac{N}{p},\alpha+\frac{N}{p}\}$ with $\alpha=\frac{N}{p}+\frac{1}{2}-\varepsilon$ for all $s\leq\frac{N}{2}+1$.\\
For the term $u^{+}\cdot\nabla c^{+h},$ we shall also proceed differently depending
on whether  $2\leq p\leq N$ or $N<p<2N$. If $2\leq p\leq N$,  we observe that applying  \eqref{prolow2} with $\sigma=\frac{N}{p}-1$ yields
\begin{equation}\label{low.14}
\|fg^{h}\|_{\dot B^{-s_0}_{2,\infty}}^{\ell}
\lesssim\|f\|_{\dot B^{1-\frac{N}{p}}_{p,1}}\Big(\|g^{h}\|_{\dot B^{\frac{N}{p}-1}_{p,1}}+\|\dot{S}_{q_0+N_0}g^{h}\|_{L^{p^{*}}}\Big)
\lesssim\|f\|_{\dot B^{1-\frac{N}{p}}_{p,1}}\|g^{h}\|_{\dot B^{\frac{N}{p}-1}_{p,1}},
\end{equation}
where we used Bernstein's inequality (recall that $p^{*}=\frac{2p}{p-2}\geq p$) and the fact that only middle frequencies of $g$ are involved in $\dot{S}_{q_0+N_0}g^{h}.$ Employing \eqref{low.8}, \eqref{low.91},  \eqref{low.11},  noticing $\frac{s_{0}+s}{2}\leq \min\{\frac{1}{2}+\alpha,2\alpha\}$ with $\alpha=\frac{N}{p}+\frac{1}{2}-\varepsilon$ for all $s\leq\frac{N}{2}+1$,  and then using  Lemma \ref{lemma2.13},  we  get
\begin{equation*} \label{s13i}
\begin{split}
&\int_0^t\langle t-\tau\rangle^{-\frac{s+s_0}{2}}\|u^{+}\cdot\nabla c^{+h}(\tau)\|_{\dot B^{-s_0}_{2,\infty}}\,d\tau\\
&\quad\lesssim\int_0^t\langle t-\tau\rangle^{-\frac{s+s_0}{2}} \|u^{+}\|_{\dot B^{1-\frac{N}{p}}_{p,1}}
\|\nabla c^{+h}\|_{\dot B^{\frac{N}{p}-1}_{p,1}}\,d\tau\\
&\quad\lesssim\int_0^t\langle t-\tau\rangle^{-\frac{s+s_0}{2}} \Big(\|u^{+\ell}\|_{\dot B^{1-\frac{N}{p}}_{p,1}}+\|u^{+h}\|_{\dot B^{1-\frac{N}{p}}_{p,1}}\Big)
\|\nabla c^{+h}\|_{\dot B^{\frac{N}{p}-1}_{p,1}}\,d\tau\\
&\quad\lesssim\Big(D^{2}(t)+X^{2}(t)\Big)\int_0^t\langle t-\tau\rangle^{-\frac{s+s_0}{2}}
\Big(\langle \tau\rangle^{-{(\frac 12+\alpha)}}
+\langle \tau\rangle^{-2\alpha}\Big) \,d\tau\\
&\quad\lesssim \langle t\rangle^{-\frac{s+s_0}{2}}
\Big(D^{2}(t)+X^{2}(t)\Big).
\end{split}
\end{equation*}
If $N<p<2N$, applying  \eqref{prolow1} with $\sigma=1-\frac{N}{p}$ yields
\begin{equation}\label{prolow11}
\|fg^{h}\|_{\dot B^{-s_0}_{2,\infty}}^{\ell}
\lesssim\Big(\|f\|_{\dot B^{1-\frac{N}{p}}_{p,1}}+\|\dot{S}_{q_0+N_0}f\|_{L^{p^{*}}}\Big)
\|g^{h}\|_{\dot B^{\frac{N}{p}-1}_{p,\infty}}
\lesssim\Big(\|f\|_{\dot B^{1-\frac{N}{p}}_{p,1}}+\|f\|_{\dot B^{\frac{N}{p}}_{2,1}}^{\ell}\Big)
\|g^{h}\|_{\dot B^{\frac{N}{p}-1}_{p,1}},
\end{equation}
where $p^{*}=\frac{2p}{p-2}.$ Then  it follows from  \eqref{low.8}, \eqref{low.92}, \eqref{low.12} and \eqref{prolow11} that
\begin{equation*} \label{s13ii}
\begin{split}
&\int_0^t\langle t-\tau\rangle^{-\frac{s+s_0}{2}}\|u^{+}\cdot\nabla c^{+h}(\tau)\|_{\dot B^{-s_0}_{2,\infty}}\,d\tau\\
&\quad\lesssim\int_0^t\langle t-\tau\rangle^{-\frac{s+s_0}{2}}
(\|u^{+}\|_{\dot B^{1-\frac{N}{p}}_{p,1}}+\|u^+\|_{\dot B^{\frac{N}{p}}_{2,1}}^{\ell})
\|\nabla c^{+h}\|_{\dot B^{\frac{N}{p}-1}_{p,1}}\,d\tau\\
&\quad\lesssim\int_0^t\langle t-\tau\rangle^{-\frac{s+s_0}{2}}
\Big(\|u^{+\ell}\|_{\dot B^{1-\frac{N}{p}}_{p,1}}+\|u^{+h}\|_{\dot B^{1-\frac{N}{p}}_{p,1}}+\|u^+\|_{\dot B^{\frac{N}{p}}_{2,1}}^{\ell}\Big)
\|\nabla c^{+h}\|_{\dot B^{\frac{N}{p}-1}_{p,1}}\,d\tau\\
&\quad\lesssim\Big(D^{2}(t)+X^{2}(t)\Big)\int_0^t\langle t-\tau\rangle^{-\frac{s+s_0}{2}}
\Big(\langle \tau\rangle^{-{\min\{\frac 12+\alpha,\frac{3N}{2p}-\frac{N}{4}+\alpha\}}}
+\langle \tau\rangle^{-\frac{3\alpha}{2}}\tau^{-\frac{\alpha}{2}}\Big) \,d\tau\\
&\quad\lesssim \langle t\rangle^{-\frac{s+s_0}{2}}
\Big(D^{2}(t)+X^{2}(t)\Big),
\end{split}
\end{equation*}
 where we have  used  Lemmas \ref{lemma2.13}-\ref{lemma2.14} and the fact $\frac{\alpha}{2}=\frac{N}{2p}+\frac{1}{4}-\frac{\varepsilon}{2}<1$, $\frac{s+s_0}{2}\leq\min\{\frac 12+\alpha,\frac{3N}{2p}-\frac{N}{4}+\alpha,2\alpha\}$ for all $s\leq\frac{N}{2}+1$, and $2\alpha>1$.\\
For the term $c^{+}\,\div (u^{+})^h,$  we shall also proceed differently depending
on whether  $2\leq p\leq N$ or $N<p<2N$. If $2\leq p\leq N$,  applying  \eqref{low.6.23}
yields that
\begin{equation*} \label{s14ii-A}
\begin{split}
&\int_0^t\langle t-\tau\rangle^{-\frac{s+s_0}{2}}\|c^{+}\,\div (u^{+})^h(\tau)\|_{\dot B^{-s_0}_{2,\infty}}\,d\tau\\
&\quad\lesssim\int_0^t\langle t-\tau\rangle^{-\frac{s+s_0}{2}}
\|c^{+}\|_{\dot B^{\frac{N}{p}-1}_{p,1}}
\|\div u^{+h}\|_{\dot B^{\frac{N}{p}-1}_{p,1}}\,d\tau.
\end{split}
\end{equation*}
In what follows, we split the integral on $[0,t]$ into  integrals
on  $[0,2]$ and $[2,t],$ respectively.
The case $t\leq2$ is obvious as $\langle t\rangle\approx1$ and $\langle t-\tau\rangle\approx1$ for $0\leq\tau\leq t\leq 2$, and  we infer that
\begin{equation*} \label{s14ii-A}
\begin{split}
&\int_0^t\langle t-\tau\rangle^{-\frac{s+s_0}{2}}\|c^{+}\,\div (u^{+})^h(\tau)\|_{\dot B^{-s_0}_{2,\infty}}\,d\tau\\
&\quad\lesssim\int_0^t\langle t-\tau\rangle^{-\frac{s+s_0}{2}}
\Big(\|c^{+}\|_{\dot B^{\frac{N}{p}-1}_{p,1}}^{\ell}+\|c^{+}\|_{\dot B^{\frac{N}{p}-1}_{p,1}}^{h}\Big)
\|u^{+}\|_{\dot B^{\frac{N}{p}}_{p,1}}^{h}\,d\tau\\
&\quad\lesssim\int_0^t\langle t-\tau\rangle^{-\frac{s+s_0}{2}}
\Big(\|c^{+}\|_{\dot B^{\frac{N}{2}-1}_{2,1}}^{\ell}+\|c^{+}\|_{\dot B^{\frac{N}{p}}_{p,1}}^{h}\Big)
\| u^{+}\|_{\dot B^{\frac{N}{p}+1}_{p,1}}^{h}\,d\tau\\
&\quad\lesssim \langle t\rangle^{-\frac{s+s_0}{2}}
X^{2}(t).
\end{split}
\end{equation*}
On the other hand, if $t\geq2,$  we infer that
\begin{equation*} \label{s14ii-B}
\begin{split}
&\int_0^t\langle t-\tau\rangle^{-\frac{s+s_0}{2}}\|c^{+}\,\div (u^{+})^h(\tau)\|_{\dot B^{-s_0}_{2,\infty}}\,d\tau\\
&\quad\lesssim\int_0^1\langle t-\tau\rangle^{-\frac{s+s_0}{2}}
\|c^{+}\|_{\dot B^{\frac{N}{p}-1}_{p,1}}
\|\div u^{+h}\|_{\dot B^{\frac{N}{p}-1}_{p,1}}\,d\tau\\&\qquad+\int_1^t\langle t-\tau\rangle^{-\frac{s+s_0}{2}}
\|c^{+}\|_{\dot B^{\frac{N}{p}-1}_{p,1}}
\|\div u^{+h}\|_{\dot B^{\frac{N}{p}-1}_{p,1}}\,d\tau
\\&\eqdefa I_{1}+I_{2}.
\end{split}
\end{equation*}
Remembering the definition of $X(t)$ and $D(t)$,  we  obtain
\begin{equation*}\label{4.15}
I_{1}\lesssim\langle t\rangle^{-\frac{s_{0}}{2}-\frac s2}X^{2}(1).
\end{equation*}
To bound the term $II_{2}$, according to  \eqref{low.11},
and using the fact that $\langle \tau\rangle\approx\tau$ when $\tau\geq1$, we conclude that
\begin{equation*}\label{4.17}
\begin{split}
I_{2}&\lesssim
\int_1^t\langle t-\tau\rangle^{-\frac{s+s_0}{2}}
\|c^{+}\|_{\dot B^{\frac{N}{p}-1}_{p,1}}
\|\div u^{+h}\|_{\dot B^{\frac{N}{p}-1}_{p,1}}\,d\tau\\
&\lesssim\int_1^t\langle t-\tau\rangle^{-\frac{s+s_0}{2}}
\Big(\|c^{+}\|_{\dot B^{\frac{N}{p}-1}_{p,1}}^{\ell}+\|c^{+}\|_{\dot B^{\frac{N}{p}-1}_{p,1}}^{h}\Big)
\| u^{+}\|_{\dot B^{\frac{N}{p}}_{p,1}}^{h}\,d\tau\\
&\lesssim\int_1^t\langle t-\tau\rangle^{-\frac{s+s_0}{2}}
\Big(\|c^{+}\|_{\dot B^{\frac{N}{2}-1}_{2,1}}^{\ell}+\|c^{+}\|_{\dot B^{\frac{N}{p}}_{p,1}}^{h}\Big)
\| u^{+}\|_{\dot B^{\frac{N}{p}}_{p,1}}^{h}\,d\tau\\
&\lesssim\int_1^t\langle t-\tau\rangle^{-\frac{s+s_0}{2}}
\Big(\|c^{+}\|_{\dot B^{\frac{N}{2}-1}_{2,1}}^{\ell}+\|c^{+}\|_{\dot B^{\frac{N}{p}}_{p,1}}^{h}\Big)
\| u^{+}\|_{\dot B^{\frac{N}{p}}_{p,1}}^{h}\,d\tau\\
&\lesssim \Big(D^{2}(t)+X^{2}(t)\Big)\int_1^t\langle t-\tau\rangle^{-\frac{s+s_0}{2}}\Big(\langle \tau\rangle^{-(\frac{s_0}{2}+\frac{N}{4}-\frac{1}{2}+\alpha)}+\langle \tau\rangle^{-2\alpha}\Big)\,d\tau\\
&\lesssim\langle t\rangle^{-\frac{s_{0}}{2}-\frac{s}{2}}\Big(D^{2}(t)+X^{2}(t)\Big),
\end{split}
\end{equation*}
where we have  used  Lemma \ref{lemma2.13} and the fact $\frac{s_{0}}{2}+\frac{s}{2}\leq \min\{\frac{s_0}{2}+\frac{N}{4}-\frac{1}{2}+\alpha, 2\alpha\}$ with $\alpha=\frac{N}{p}+\frac{1}{2}-\varepsilon$ for all $s\leq\frac{N}{2}+1$.
If $N<p<2N$,  thanks to \eqref{low.8},  \eqref{low.92}, \eqref{low.11}, \eqref{low.12} and \eqref{prolow11}, we deduce that
\begin{equation*} \label{s14ii}
\begin{split}
&\int_0^t\langle t-\tau\rangle^{-\frac{s+s_0}{2}}\|c^{+}\,\div (u^{+})^h(\tau)\|_{\dot B^{-s_0}_{2,\infty}}\,d\tau\\
&\quad\lesssim\int_0^t\langle t-\tau\rangle^{-\frac{s+s_0}{2}}
\Big(\|c^{+}\|_{\dot B^{1-\frac{N}{p}}_{p,1}}+\|c^{+}\|_{\dot B^{\frac{N}{p}}_{2,1}}^{\ell}\Big)
\|\div u^{+h}\|_{\dot B^{\frac{N}{p}-1}_{p,1}}\,d\tau\\
&\quad\lesssim\int_0^t\langle t-\tau\rangle^{-\frac{s+s_0}{2}}
\Big(\|c^{+}\|_{\dot B^{1-\frac{N}{p}}_{p,1}}^{\ell}+\|c^{+}\|_{\dot B^{1-\frac{N}{p}}_{p,1}}^{h}+\|c^{+}\|_{\dot B^{\frac{N}{p}}_{2,1}}^{\ell}\Big)
\| u^{+h}\|_{\dot B^{\frac{N}{p}}_{p,1}}\,d\tau\\
&\quad\lesssim\Big(D^{2}(t)+X^{2}(t)\Big)\int_0^t\langle t-\tau\rangle^{-\frac{s+s_0}{2}}
\langle \tau\rangle^{-\min\{\frac{1+\alpha}{2},\frac {3\alpha}{2},\frac{3N}{2p}-\frac{N}{4}+\frac{\alpha}{2}\}}\tau^{-\frac{\alpha}{2}}\,d\tau\\
&\quad\lesssim \langle t\rangle^{-\frac{s+s_0}{2}}
\Big(D^{2}(t)+X^{2}(t)\Big),
\end{split}
\end{equation*}
where we have  used  Lemma \ref{lemma2.14} and the fact $\frac{\alpha}{2}<1$, $\min\{\frac 12+\alpha,2\alpha,\frac{3N}{2p}-\frac{N}{4}+\alpha\}>1$ for $N\leq p\leq\min(4,\frac{2N}{N-2})$, $\frac{s_{0}}{2}+\frac{s}{2}\leq \min\{\frac 12+\alpha,2\alpha,\frac{3N}{2p}-\frac{N}{4}+\alpha\}$ with $\alpha=\frac{N}{p}+\frac{1}{2}-\varepsilon$ for all $s\leq\frac{N}{2}+1$.\\
Then we infer that
\begin{equation}\label{s1}
\int_0^t\langle t-\tau\rangle^{-\frac{s+s_0}{2}} \big\|H_{1}(\tau)\big\|_{\dot B^{-\frac{N}{2}}_{2,\infty}}^\ell d\tau
\lesssim\langle t\rangle^{-{\frac{s+s_0}{2}}}
\Big(D^{2}(t)+X^{2}(t)\Big).
\end{equation}
The term  $H_{3}$ may be treated along the same lines,  and we obtain
\begin{equation}\label{s3}
\int_0^t\langle t-\tau\rangle^{-\frac{s+s_0}{2}} \big\|H_{3}(\tau)\big\|_{\dot B^{-\frac{N}{2}}_{2,\infty}}^\ell d\tau
\lesssim\langle t\rangle^{-\frac{s+s_0}{2}}
\Big(D^{2}(t)+X^{2}(t)\Big).
\end{equation}

Next, we bound the term $H_{2}^{i}$. To handle the first term  $g_{+}(c^{+},c^{-})\partial_{i}c^{+}$ in $H_{2}^{i}$,  we decompose  it into
$$g_{+}(c^{+},c^{-})\partial_{i}c^{+}=g_{+}(c^{+},c^{-})\partial_{i}c^{+\ell}
+g_{+}(c^{+},c^{-})\partial_{i}c^{+h},
$$
where $g_{+}$ stands for some smooth function vanishing at $0$.\\
For the term $g_{+}(c^{+},c^{-})\partial_{i}c^{+\ell}$, employing  \eqref{low.7}, \eqref{low.8}, \eqref{low.10}, \eqref{low.11}, \eqref{99111B},  Proposition \ref{p28} and Lemma \ref{lemma2.13}, we conclude that
\begin{equation*}
\begin{split}\label{s21l}
&\int_0^t\langle t-\tau\rangle^{-\frac{s+s_0}{2}}\|g_{+}(c^{+},c^{-})\partial_{i}c^{+\ell}(\tau)\|_{\dot B^{-s_0}_{2,\infty}}\,d\tau\\
&\quad\lesssim
\int_0^t\langle t-\tau\rangle^{-\frac{s+s_0}{2}}
\|g_{+}(c^{+},c^{-})\|_{\dot B^{1-\frac{N}{p}}_{p,1}}
\|\nabla c^{+\ell}\|_{\dot B^{\frac{N}{2}-1}_{2,1}}\\
&\quad\lesssim\int_0^t\langle t-\tau\rangle^{-\frac{s+s_0}{2}}\Big(\|(c^{+},c^{-})^{\ell}\|_{\dot B^{1-\frac{N}{p}}_{p,1}}+\|(c^{+},c^{-})^{h}\|_{\dot B^{1-\frac{N}{p}}_{p,1}}\Big)\|\nabla c^{+\ell}\|_{\dot B^{\frac{N}{2}-1}_{2,1}}\,d\tau\\
&\quad\lesssim D^{2}(t)\int_0^t\langle t-\tau\rangle^{-\frac{s+s_0}{2}}
(\langle \tau\rangle^{-\min\{\frac 12+\frac{N}{p},\alpha+\frac{N}{p}\}}\,d\tau
\\&\quad\lesssim \langle t\rangle^{-\frac{s+s_0}{2}} D^{2}(t).
\end{split}
\end{equation*}
For the term $g_{+}(c^{+},c^{-})\partial_{i}c^{+h}$, if $2\leq p\leq N,$  it follows from  \eqref{low.8}, \eqref{low.11}, \eqref{low.14}, \eqref{99111B},  Proposition \ref{p28}  and Lemma \ref{lemma2.13},  that
\begin{equation*} \label{s21hi}
\begin{split}
&\int_0^t\langle t-\tau\rangle^{-\frac{s+s_0}{2}}\|g_{+}(c^{+},c^{-})\partial_{i}c^{+h}(\tau)\|_{\dot B^{-s_0}_{2,\infty}}\,d\tau\\
&\quad\lesssim\int_0^t\langle t-\tau\rangle^{-\frac{s+s_0}{2}} \|(c^{+},c^{-})\|_{\dot B^{1-\frac{N}{p}}_{p,1}}
\|\nabla c^{+h}\|_{\dot B^{\frac{N}{p}-1}_{p,1}}\,d\tau\\
&\quad\lesssim\int_0^t\langle t-\tau\rangle^{-\frac{s+s_0}{2}}
(\langle \tau\rangle^{-{(\frac 12+\alpha)}}
+\langle \tau\rangle^{2\alpha}) \,d\tau\\
&\quad\lesssim \langle t\rangle^{-\frac{s+s_0}{2}}
D^{2}(\tau).
\end{split}
\end{equation*}
If $N<p\leq 2N,$  employing \eqref{low.8}, \eqref{low.11}, \eqref{low.12}, \eqref{prolow11}, \eqref{99111B}, Proposition \ref{p27} and Lemma \ref{lemma2.13} yields that
\begin{equation*} \label{s21hii}
\begin{split}
&\int_0^t\langle t-\tau\rangle^{-\frac{s+s_0}{2}}\|g_{+}(c^{+},c^{-})\partial_{i}c^{+h}(\tau)\|_{\dot B^{-s_0}_{2,\infty}}\,d\tau\\
&\quad\lesssim\int_0^t\langle t-\tau\rangle^{-\frac{s+s_0}{2}}
\Big(\|(c^{+},c^{-})\|_{\dot B^{1-\frac{N}{p}}_{p,1}}+\|(c^{+},c^{-})\|_{\dot B^{\frac{N}{p}}_{2,1}}^{\ell}\Big)
\|\nabla c^{+h}\|_{\dot B^{\frac{N}{p}-1}_{p,1}}\,d\tau\\
&\quad\lesssim\int_0^t\langle t-\tau\rangle^{-\frac{s+s_0}{2}}
\langle \tau\rangle^{-\min\{\frac 12+\alpha,2\alpha,\frac{3N}{2p}-\frac{N}{4}+\alpha\}} \,d\tau\\
&\quad\lesssim \langle t\rangle^{-\frac{s+s_0}{2}}
D^{2}(\tau). \\
\end{split}
\end{equation*}
Thus
\begin{equation*} \label{ds21}
\begin{split}
\int_0^t\langle t-\tau\rangle^{-\frac{s+s_0}{2}}\|g_{+}(c^{+},c^{-})\partial_{i}c^{+}(\tau)\|_{\dot B^{-s_0}_{2,\infty}}\,d\tau
\lesssim \langle t\rangle^{-\frac{s+s_0}{2}}
D^{2}(\tau). \\
\end{split}
\end{equation*}
Similarly,
\begin{equation*}
\begin{split}\label{s22}
&\int_0^t\langle t-\tau\rangle^{-\frac{s+s_0}{2}}\|\tilde{g}_{+}(c^{+},c^{-})\partial_{i}c^{-}\|_{\dot B^{-s_0}_{2,\infty}}\,d\tau
\lesssim \langle t\rangle^{-\frac{s+s_0}{2}} D^{2}(t).
\end{split}
\end{equation*}
To bound the term with $(u^{+}\cdot\nabla)u_{i}^{+}$ in $H_{2}^{i}$, we employ the following decomposition:
$$(u^{+}\cdot\nabla)u_{i}^{+}=(u^{+}\cdot\nabla)(u_{i}^{+})^{\ell}
+(u^{+}\cdot\nabla)(u_{i}^{+})^{h}.
$$
For the term  $(u^{+}\cdot\nabla)(u_{i}^{+})^{\ell}$, if $2\leq p\leq N$, it follows from  \eqref{low.7},  \eqref{low.8}, \eqref{low.10}, and \eqref{low.91},  that
\begin{equation*}
\begin{split}
&\int_0^t\langle t-\tau\rangle^{-\frac{s+s_0}{2}}\|(u^{+}\cdot\nabla)(u_{i}^{+})^{\ell}(\tau)\|_{\dot B^{-s_0}_{2,\infty}}\,d\tau\\
&\quad\lesssim
\int_0^t\langle t-\tau\rangle^{-\frac{s+s_0}{2}}
\|u^{+}(\tau)\|_{\dot B^{1-\frac{N}{p}}_{p,1}}
\|\nabla u_{i}^{+\ell}(\tau)\|_{\dot B^{\frac{N}{2}-1}_{2,1}}\\
&\quad\lesssim\int_0^t\langle t-\tau\rangle^{-\frac{s+s_0}{2}}\Big(\|u^{+}(\tau)\|_{\dot B^{1-\frac{N}{p}}_{p,1}}^{\ell}+\|u^{+}(\tau)\|_{\dot B^{1-\frac{N}{p}}_{p,1}}^h\Big)\|\nabla u^{+}\|_{\dot B^{\frac{N}{2}-1}_{2,1}}^{\ell}\,d\tau\\
&\quad\lesssim \Big(D^{2}(t)+X^{2}(t)\Big)\int_0^t\langle t-\tau\rangle^{-\frac{s+s_0}{2}}\Big(\langle \tau\rangle^{-{(\frac Np+\frac 12)}}+\langle \tau\rangle^{-{(\frac Np+\alpha)}} \Big)\,d\tau\\
&\quad\lesssim \langle t\rangle^{-\frac{s+s_0}{2}} \Big(D^{2}(t)+X^{2}(t)\Big).
\end{split}
\end{equation*}
If $N< p\leq 2N,$  using \eqref{low.7}, \eqref{low.8}, \eqref{low.10}, and \eqref{low.92}  yields that
\begin{equation*}
\begin{split}
&\int_0^t\langle t-\tau\rangle^{-\frac{s+s_0}{2}}\|(u^{+}\cdot\nabla)(u_{i}^{+})^{\ell}(\tau)\|_{\dot B^{-s_0}_{2,\infty}}\,d\tau\\
&\quad\lesssim
\int_0^t\langle t-\tau\rangle^{-\frac{s+s_0}{2}}
\|u^{+}(\tau)\|_{\dot B^{1-\frac{N}{p}}_{p,1}}
\|\nabla u_{i}^{+\ell}(\tau)\|_{\dot B^{\frac{N}{2}-1}_{2,1}}\\
&\quad\lesssim\int_0^t\langle t-\tau\rangle^{-\frac{s+s_0}{2}}\Big(\|u^{+}(\tau)\|_{\dot B^{1-\frac{N}{p}}_{p,1}}^{\ell}+\|u^{+}(\tau)\|_{\dot B^{1-\frac{N}{p}}_{p,1}}^h\Big)\|\nabla u^{+}(\tau)\|_{\dot B^{\frac{N}{2}-1}_{2,1}}^{\ell}\,d\tau\\
&\quad\lesssim \Big(D^{2}(t)+X^{2}(t)\Big)\int_0^t\langle t-\tau\rangle^{-\frac{s+s_0}{2}}(\langle \tau\rangle^{-{(\frac Np+\frac 12)}}+\langle \tau\rangle^{-{(\frac Np+\frac \alpha2)}} \tau^{-{\frac \alpha2}})\,d\tau\\
&\quad\lesssim \langle t\rangle^{-\frac{s+s_0}{2}} \Big(D^{2}(t)+X^{2}(t)\Big).
\end{split}
\end{equation*}
To deal with the term $(u^{+}\cdot\nabla)(u_{i}^{+})^{h},$ let us first consider the case $2\leq p\leq N.$
Applying \eqref{low.6.23} implies  that
\begin{equation*} \label{s14ii-A-C}
\begin{split}
&\int_0^t\langle t-\tau\rangle^{-\frac{s+s_0}{2}}\|(u^{+}\cdot\nabla)(u_{i}^{+})^{h}(\tau)\|_{\dot B^{-s_0}_{2,\infty}}\,d\tau\\
&\quad\lesssim\int_0^t\langle t-\tau\rangle^{-\frac{s+s_0}{2}}
\|u^{+}(\tau)\|_{\dot B^{\frac{N}{p}-1}_{p,1}}
\|\nabla u^{+}(\tau)\|_{\dot B^{\frac{N}{p}-1}_{p,1}}^{h}\,d\tau.
\end{split}
\end{equation*}
In what follows, we divide the integral on $[0,t]$ into  integrals
on  $[0,2]$ and $[2,t],$ respectively.
When $t\leq2$, thus  $\langle t\rangle\approx1$ and $\langle t-\tau\rangle\approx1$ for $0\leq\tau\leq t\leq 2$.  We have
\begin{equation*} \label{s14ii-A}
\begin{split}
&\int_0^t\langle t-\tau\rangle^{-\frac{s+s_0}{2}}\|(u^{+}\cdot\nabla)(u_{i}^{+})^{h}(\tau)\|_{\dot B^{-s_0}_{2,\infty}}\,d\tau\\
&\quad\lesssim\int_0^t\langle t-\tau\rangle^{-\frac{s+s_0}{2}}
\Big(\|u^{+}(\tau)\|_{\dot B^{\frac{N}{p}-1}_{p,1}}^{\ell}+\|u^{+}(\tau)\|_{\dot B^{\frac{N}{p}-1}_{p,1}}^{h}\Big)
\|u^{+}(\tau)\|_{\dot B^{\frac{N}{p}}_{p,1}}^{h}\,d\tau\\
&\quad\lesssim\int_0^t\langle t-\tau\rangle^{-\frac{s+s_0}{2}}
\Big(\|u^{+}(\tau)\|_{\dot B^{\frac{N}{2}-1}_{2,1}}^{\ell}+\|u^{+}(\tau)\|_{\dot B^{\frac{N}{p}}_{p,1}}^{h}\Big)
\| u^{+}(\tau)\|_{\dot B^{\frac{N}{p}+1}_{p,1}}^{h}\,d\tau\\
&\quad\lesssim \langle t\rangle^{-\frac{s+s_0}{2}}
X^{2}(t).
\end{split}
\end{equation*}
When $t\geq2,$  we deduce  that
\begin{equation*} \label{s14ii-B}
\begin{split}
&\int_0^t\langle t-\tau\rangle^{-\frac{s+s_0}{2}}\|(u^{+}\cdot\nabla)(u_{i}^{+})^{h}(\tau)\|_{\dot B^{-s_0}_{2,\infty}}\,d\tau\\
&\quad\lesssim\int_0^t\langle t-\tau\rangle^{-\frac{s+s_0}{2}}
\|u^{+}(\tau)\|_{\dot B^{\frac{N}{p}-1}_{p,1}}
\|\nabla u^{+}(\tau)\|_{\dot B^{\frac{N}{p}-1}_{p,1}}^{h}\,d\tau\\
&\quad=\int_0^1\langle t-\tau\rangle^{-\frac{s+s_0}{2}}
\|u^{+}(\tau)\|_{\dot B^{\frac{N}{p}-1}_{p,1}}
\|\nabla u^{+}(\tau)\|_{\dot B^{\frac{N}{p}-1}_{p,1}}^{h}\,d\tau\\&\qquad+\int_1^t\langle t-\tau\rangle^{-\frac{s+s_0}{2}}
\|u^{+}(\tau)\|_{\dot B^{\frac{N}{p}-1}_{p,1}}
\|\nabla u^{+}(\tau)\|_{\dot B^{\frac{N}{p}-1}_{p,1}}^{h}\,d\tau
\\&\eqdefa II_{1}+II_{2}.
\end{split}
\end{equation*}
Using the definition of $X(t)$ and $D(t)$,  we  get
\begin{equation*}\label{4.15}
II_{1}\lesssim\langle t\rangle^{-\frac{s_{0}}{2}-\frac s2}X^{2}(1).
\end{equation*}
For the term $II_{2}$, according to \eqref{low.11},
and using the fact that $\langle \tau\rangle\approx\tau$ when $\tau\geq1$, we conclude,  from Lemma \ref{lemma2.13},  that
\begin{equation*}\label{4.17}
\begin{split}
II_{2}&\lesssim
\int_1^t\langle t-\tau\rangle^{-\frac{s+s_0}{2}}
\|u^{+}(\tau)\|_{\dot B^{\frac{N}{p}-1}_{p,1}}
\|\nabla u^{+}(\tau)\|_{\dot B^{\frac{N}{p}-1}_{p,1}}^{h}\,d\tau\\
&\lesssim\int_1^t\langle t-\tau\rangle^{-\frac{s+s_0}{2}}
\Big(\|u^{+}\|_{\dot B^{\frac{N}{p}-1}_{p,1}}^{\ell}+\|u^{+}\|_{\dot B^{\frac{N}{p}-1}_{p,1}}^{h}\Big)
\| u^{+}\|_{\dot B^{\frac{N}{p}}_{p,1}}^{h}\,d\tau\\
&\lesssim\int_1^t\langle t-\tau\rangle^{-\frac{s+s_0}{2}}
\Big(\|u^{+}\|_{\dot B^{\frac{N}{2}-1}_{2,1}}^{\ell}+\|u^{+}\|_{\dot B^{\frac{N}{p}}_{p,1}}^{h}\Big)
\| u^{+}\|_{\dot B^{\frac{N}{p}}_{p,1}}^{h}\,d\tau\\
&\lesssim\int_1^t\langle t-\tau\rangle^{-\frac{s+s_0}{2}}
\Big(\|u^{+}\|_{\dot B^{\frac{N}{2}-1}_{2,1}}^{\ell}+\|u^{+}\|_{\dot B^{\frac{N}{p}}_{p,1}}^{h}\Big)
\| u^{+}\|_{\dot B^{\frac{N}{p}}_{p,1}}^{h}\,d\tau\\
&\lesssim \Big(D^{2}(t)+X^{2}(t)\Big)\int_1^t\langle t-\tau\rangle^{-\frac{s+s_0}{2}}\Big(\langle \tau\rangle^{-(\frac{s_0}{2}+\frac{N}{4}-\frac{1}{2}+\alpha)}+\langle \tau\rangle^{-2\alpha}\Big)\,d\tau\\
&\lesssim\langle t\rangle^{-\frac{s_{0}}{2}-\frac{s}{2}}\Big(D^{2}(t)+X^{2}(t)\Big).
\end{split}
\end{equation*}
If $N<p\leq 2N, $ according  to \eqref{low.8},  \eqref{low.92}, \eqref{low.11}, \eqref{low.12} and \eqref{prolow11}, we deduce, from Lemma \ref{lemma2.13}, that
\begin{equation*} \label{s14ii-F}
\begin{split}
&\int_0^t\langle t-\tau\rangle^{-\frac{s+s_0}{2}}\|(u^{+}\cdot\nabla)(u_{i}^{+})^{h}(\tau)\|_{\dot B^{-s_0}_{2,\infty}}\,d\tau\\
&\quad\lesssim\int_0^t\langle t-\tau\rangle^{-\frac{s+s_0}{2}}
\Big(\|u^{+}\|_{\dot B^{1-\frac{N}{p}}_{p,1}}+\|u^{+}\|_{\dot B^{\frac{N}{p}}_{2,1}}^{\ell}\Big)
\|\div u^{+h}\|_{\dot B^{\frac{N}{p}-1}_{p,1}}\,d\tau\\
&\quad\lesssim\int_0^t\langle t-\tau\rangle^{-\frac{s+s_0}{2}}
\Big(\|u^{+}\|_{\dot B^{1-\frac{N}{p}}_{p,1}}^{\ell}+\|u^{+}\|_{\dot B^{1-\frac{N}{p}}_{p,1}}^{h}+\|u^{+}\|_{\dot B^{\frac{N}{p}}_{2,1}}^{\ell}\Big)
\| u^{+h}\|_{\dot B^{\frac{N}{p}}_{p,1}}\,d\tau\\
&\quad\lesssim\Big(D^{2}(t)+X^{2}(t)\Big)\int_0^t\langle t-\tau\rangle^{-\frac{s+s_0}{2}}
\langle \tau\rangle^{-\min\{\frac{1+\alpha}{2},\frac {3\alpha}{2},\frac{3N}{2p}-\frac{N}{4}+\frac{\alpha}{2}\}}\tau^{-\frac{\alpha}{2}}\,d\tau\\
&\quad\lesssim \langle t\rangle^{-\frac{s+s_0}{2}}
\Big(D^{2}(t)+X^{2}(t)\Big).
\end{split}
\end{equation*}
Therefore
\begin{equation*} \label{ds22}
\begin{split}
\int_0^t\langle t-\tau\rangle^{-\frac{s+s_0}{2}}\|(u^{+}\cdot\nabla)u_{i}^{+}(\tau)\|_{\dot B^{-s_0}_{2,\infty}}\,d\tau
\lesssim \langle t\rangle^{-\frac{s+s_0}{2}}
D^{2}(\tau). \\
\end{split}
\end{equation*}
To bound the term with $\mu^{+}h_{+}(c^{+},c^{-})\partial_{j}c^{+}\partial_{j}u^{+}_{i}$ in $H_{2}^{i}$,  employing the following low-high frequency decomposition yields that
$$\mu^{+}h_{+}(c^{+},c^{-})\partial_{j}c^{+}\partial_{j}u^{+}_{i}
=\mu^{+}h_{+}(c^{+},c^{-})\partial_{j}c^{+}\partial_{j}u^{+\ell}_{i}
+\mu^{+}h_{+}(c^{+},c^{-})\partial_{j}c^{+}\partial_{j}u^{+h}_{i}.
$$
For the term $\mu^{+}h_{+}(c^{+},c^{-})\partial_{j}c^{+}\partial_{j}u^{+\ell}_{i}$, it follows from  $\alpha\geq\frac{N}{p}$, \eqref{low.10}, \eqref{low.11}, \eqref{99111B}, Proposition \ref{p28} and Lemma \ref{lemma2.13}, that
\begin{equation*}
\begin{split}\label{s21l}
&\int_0^t\langle t-\tau\rangle^{-\frac{s+s_0}{2}}\|\mu^{+}h_{+}(c^{+},c^{-})\partial_{j}c^{+}\partial_{j}u^{+\ell}_{i}
(\tau)\|_{\dot B^{-s_0}_{2,\infty}}^{\ell}\,d\tau\\
&\quad\lesssim
\int_0^t\langle t-\tau\rangle^{-\frac{s+s_0}{2}}
\|\mu^{+}h_{+}(c^{+},c^{-})\partial_{j}c^{+}\partial_{j}u^{+\ell}_{i}
(\tau)\|_{\dot B^{-\frac{N}{p}}_{2,\infty}}^{\ell}\\
&\quad\lesssim
\int_0^t\langle t-\tau\rangle^{-\frac{s+s_0}{2}}
\|h_{+}(c^{+},c^{-})\nabla c^{+}\|_{\dot B^{\frac{N}{p}-1}_{p,1}}
\|\nabla u^{+\ell}\|_{\dot B^{1-\frac{N}{p}}_{2,1}}\\
&\quad\lesssim\int_0^t\langle t-\tau\rangle^{-\frac{s+s_0}{2}}\| c^{+}\|_{\dot B^{\frac{N}{p}}_{p,1}}\| u^{+}\|_{\dot B^{2-\frac{N}{p}}_{2,1}}^{\ell}\,d\tau\\
&\quad\lesssim\int_0^t\langle t-\tau\rangle^{-\frac{s+s_0}{2}}\Big(\|c^{+}\|_{\dot B^{\frac{N}{2}}_{2,1}}^{\ell}+\|c^{+}\|_{\dot B^{\frac{N}{p}}_{p,1}}^{h}\Big)\| u^{+}\|_{\dot B^{2-\frac{N}{p}}_{2,1}}^{\ell}\,d\tau\\
&\quad\lesssim \Big(X^{2}(t)+D^{2}(t)\Big)\int_0^t\langle t-\tau\rangle^{-\frac{s+s_0}{2}}
\Big(\langle \tau\rangle^{-(1-\frac N4+\frac{3N}{2p})}+\langle\tau\rangle^{-(1+\alpha+\frac{N}{p}-\frac N4)}\Big)\,d\tau
\\&\quad\lesssim \langle t\rangle^{-\frac{s+s_0}{2}}\Big(X^{2}(t)+D^{2}(t)\Big),
\end{split}
\end{equation*}
where we have  used  the embedding  $\dot B^{-\frac{N}{p}}_{2,\infty} \hookrightarrow \dot B^{-s_0}_{2,\infty}$ as $s_0\leq\frac{N}{p}$ (for $p\geq2$),       $\|fg\|_{\dot B^{-\frac{N}{p}}_{2,\infty}}
\lesssim\|f\|_{\dot B^{\frac{N}{p}-1}_{p,1}}
\|g\|_{\dot B^{1-\frac{N}{p}}_{2,1}}$ and the fact $\frac{s_{0}}{2}+\frac{s}{2}\leq \min\{1+\alpha+\frac{N}{p}-\frac N4, 1-\frac N4+\frac{3N}{2p}\}$ for all $s\leq\frac{N}{2}+1, p\leq\frac{2N}{N-2}$.\\
For the term $\mu^{+}h_{+}(c^{+},c^{-})\partial_{j}c^{+}\partial_{j}u^{+h}_{i}$,  let us first consider the case $2\leq p\leq N.$  It follows from  \eqref{99111B},  \eqref{low.6.23} and Proposition \ref{p27},   that
\begin{equation*}
\begin{split}
&\int_0^t\langle t-\tau\rangle^{-\frac{s+s_0}{2}}\|\mu^{+}h_{+}(c^{+},c^{-})\partial_{j}c^{+}\partial_{j}u^{+h}_{i}
(\tau)\|_{\dot B^{-s_0}_{2,\infty}}\,d\tau\\
&\quad\lesssim\int_0^t\langle t-\tau\rangle^{-\frac{s+s_0}{2}}\|\nabla c^{+}
(\tau)\|_{\dot B^{\frac{N}{p}-1}_{p,1}}
\|\nabla u^{+h}(\tau)\|_{\dot B^{\frac{N}{p}-1}_{p,1}}\,d\tau\\
&\quad\lesssim\int_0^t\langle t-\tau\rangle^{-\frac{s+s_0}{2}}\| c^{+}
(\tau)\|_{\dot B^{\frac{N}{p}}_{p,1}}
\|\nabla u^{+}(\tau)\|_{\dot B^{\frac{N}{p}-1}_{p,1}}^{h}\,d\tau.
\end{split}
\end{equation*}
Due to  $\langle t\rangle\approx1$ and $\langle t-\tau\rangle\approx1$ for $0\leq\tau\leq t\leq 2$, we  deduce,  if $t\leq2$,  that
\begin{equation*} \label{s14ii-A-D}
\begin{split}
&\int_0^t\langle t-\tau\rangle^{-\frac{s+s_0}{2}}\|\mu^{+}h_{+}(c^{+},c^{-})\partial_{j}c^{+}\partial_{j}u^{+h}_{i}
(\tau)\|_{\dot B^{-s_0}_{2,\infty}}\,d\tau\\
&\quad\lesssim\int_0^t\langle t-\tau\rangle^{-\frac{s+s_0}{2}}
\Big(\|c^{+}\|_{\dot B^{\frac{N}{p}-1}_{p,1}}^{\ell}+\|c^{+}\|_{\dot B^{\frac{N}{p}}_{p,1}}^{h}\Big)
\|u^{+}\|_{\dot B^{\frac{N}{p}+1}_{p,1}}^{h}\,d\tau\\
&\quad\lesssim\int_0^t\langle t-\tau\rangle^{-\frac{s+s_0}{2}}
\Big(\|c^{+}\|_{\dot B^{\frac{N}{2}-1}_{2,1}}^{\ell}+\|c^{+}\|_{\dot B^{\frac{N}{p}}_{p,1}}^{h}\Big)
\| u^{+}\|_{\dot B^{\frac{N}{p}+1}_{p,1}}^{h}\,d\tau\\
&\quad\lesssim \langle t\rangle^{-\frac{s+s_0}{2}}
X^{2}(t).
\end{split}
\end{equation*}
If $t\geq2,$  we conclude  that
\begin{equation*} \label{s14ii-B-D}
\begin{split}
&\int_0^t\langle t-\tau\rangle^{-\frac{s+s_0}{2}}\|\mu^{+}h_{+}(c^{+},c^{-})\partial_{j}c^{+}\partial_{j}u^{+h}_{i}
(\tau)\|_{\dot B^{-s_0}_{2,\infty}}\,d\tau
\\&\quad\lesssim\int_0^1\langle t-\tau\rangle^{-\frac{s+s_0}{2}}
\|c^{+}\|_{\dot B^{\frac{N}{p}}_{p,1}}
\|\nabla u^{+}\|_{\dot B^{\frac{N}{p}-1}_{p,1}}^{h}\,d\tau\\&\qquad+\int_1^t\langle t-\tau\rangle^{-\frac{s+s_0}{2}}
\|c^{+}\|_{\dot B^{\frac{N}{p}}_{p,1}}
\|\nabla u^{+}\|_{\dot B^{\frac{N}{p}-1}_{p,1}}^{h}\,d\tau
\\&\eqdefa III_{1}+III_{2}.
\end{split}
\end{equation*}
From the definition of $X(t)$ and $D(t)$,  we  get
\begin{equation*}\label{4.15-DD}
III_{1}\lesssim\langle t\rangle^{-\frac{s_{0}}{2}-\frac s2}X^{2}(1).
\end{equation*}
To handle the term $II_{2}$, according to \eqref{low.11}, Lemma \ref{lemma2.13},
and using the fact that $\langle \tau\rangle\approx\tau$ when $\tau\geq1$, we infer that
\begin{equation*}\label{4.17-DD}
\begin{split}
II_{2}&\lesssim
\int_1^t\langle t-\tau\rangle^{-\frac{s+s_0}{2}}
\|c^{+}\|_{\dot B^{\frac{N}{p}}_{p,1}}
\|\nabla u^{+}\|_{\dot B^{\frac{N}{p}-1}_{p,1}}^{h}\,d\tau\\
&\lesssim\int_1^t\langle t-\tau\rangle^{-\frac{s+s_0}{2}}
\Big(\|c^{+}\|_{\dot B^{\frac{N}{p}-1}_{p,1}}^{\ell}+\|c^{+}\|_{\dot B^{\frac{N}{p}}_{p,1}}^{h}\Big)
\|\nabla u^{+}\|_{\dot B^{\frac{N}{p}}_{p,1}}^{h}\,d\tau\\
&\lesssim\int_1^t\langle t-\tau\rangle^{-\frac{s+s_0}{2}}
\Big(\|c^{+}\|_{\dot B^{\frac{N}{2}-1}_{2,1}}^{\ell}+\|c^{+}\|_{\dot B^{\frac{N}{p}}_{p,1}}^{h}\Big)
\|\nabla u^{+}\|_{\dot B^{\frac{N}{p}}_{p,1}}^{h}\,d\tau\\
&\lesssim \Big(D^{2}(t)+X^{2}(t)\Big)\int_1^t\langle t-\tau\rangle^{-\frac{s+s_0}{2}}\Big(\langle \tau\rangle^{-(\frac{s_0}{2}+\frac{N}{4}-\frac{1}{2}+\alpha)}+\langle \tau\rangle^{-2\alpha}\Big)\,d\tau\\
&\lesssim\langle t\rangle^{-\frac{s_{0}}{2}-\frac{s}{2}}\Big(D^{2}(t)+X^{2}(t)\Big).
\end{split}
\end{equation*}
Let us consider the case $N\leq p\leq 2N.$  Applying \eqref {prolow2} with $\sigma=1-\frac{N}{p}$ yields that
\begin{equation*}\label{prolow2}\begin{split}
\|h_{+}(c^{+},c^{-})\partial_{j}c^{+}\partial_{j}u^{+h}_{i}\|_{\dot B^{-s_0}_{2,\infty}}^{\ell}
&\lesssim\Big(\|\nabla u^{+h}\|_{\dot B^{1-\frac{N}{p}}_{p,1}}+\sum_{k=k_0}^{k_0+N_0-1}\|\dot{\Delta}_{k}\nabla u^{+h}\|_{L^{p^{*}}}\Big)
\|h_{+}(c^{+},c^{-})\nabla c^{+}\|_{\dot B^{\frac{N}{p}-1}_{p,\infty}}
\\&\lesssim \|\nabla u^{+h}\|_{\dot B^{1-\frac{N}{p}}_{p,1}}
\|h_{+}(c^{+},c^{-})\nabla c^{+}\|_{\dot B^{\frac{N}{p}-1}_{p,\infty}},
\end{split}\end{equation*}
where we have used $1-\frac{N}{p}>0$ and  Bernstein inequality $\|\dot{\Delta}_{k}\nabla u^{+h}\|_{L^{p^{*}}}\lesssim \|\dot{\Delta}_{k}\nabla u^{+h}\|_{L^{p}}$ for $k_0\leq k\leq k_0+N_0$ as $p^{*}\geq p$. Further, due to the smooth function $h_{+}(c^{+},c^{-})$ vanishing at $(0,0)$ and $1-\frac{N}{p}\leq \frac{N}{p}$, employing \eqref{99111B} and Proposition \ref{p28}, we thus get
\begin{equation*}
\begin{split}
&\int_0^t\langle t-\tau\rangle^{-\frac{s+s_0}{2}}\|\mu^{+}h_{+}(c^{+},c^{-})\partial_{j}c^{+}\partial_{j}u^{+h}_{i}
(\tau)\|_{\dot B^{-s_0}_{2,\infty}}\,d\tau\\
&\quad\lesssim\int_0^t\langle t-\tau\rangle^{-\frac{s+s_0}{2}}
\|\nabla c^{+}\|_{\dot B^{\frac{N}{p}-1}_{p,1}}
\|\nabla u^{+}\|_{\dot B^{1-\frac{N}{p}}_{p,1}}^{h}\,d\tau\\
&\quad\lesssim\int_0^t\langle t-\tau\rangle^{-\frac{s+s_0}{2}}
\|c^{+}\|_{\dot B^{\frac{N}{p}}_{p,1}}
\|\nabla u^{+}\|_{\dot B^{\frac{N}{p}}_{p,1}}^{h}\,d\tau.
\end{split}
\end{equation*}
As $\langle t\rangle\approx1$ and $\langle t-\tau\rangle\approx1$ for $0\leq\tau\leq t\leq 2$,  we deduce, for the case $t\leq2$, that
\begin{equation*} \label{s14ii-A-E}
\begin{split}
&\int_0^t\langle t-\tau\rangle^{-\frac{s+s_0}{2}}\|\mu^{+}h_{+}(c^{+},c^{-})\partial_{j}c^{+}\partial_{j}u^{+h}_{i}
(\tau)\|_{\dot B^{-s_0}_{2,\infty}}\,d\tau\\
&\quad\lesssim\int_0^t\langle t-\tau\rangle^{-\frac{s+s_0}{2}}
\Big(\|c^{+}\|_{\dot B^{\frac{N}{p}-1}_{p,1}}^{\ell}+\|c^{+}\|_{\dot B^{\frac{N}{p}}_{p,1}}^{h}\Big)
\|u^{+}\|_{\dot B^{\frac{N}{p}+1}_{p,1}}^{h}\,d\tau\\
&\quad\lesssim\int_0^t\langle t-\tau\rangle^{-\frac{s+s_0}{2}}
\Big(\|c^{+}\|_{\dot B^{\frac{N}{2}-1}_{2,1}}^{\ell}+\|c^{+}\|_{\dot B^{\frac{N}{p}}_{p,1}}^{h}\Big)
\| u^{+}\|_{\dot B^{\frac{N}{p}+1}_{p,1}}^{h}\,d\tau\\
&\quad\lesssim \langle t\rangle^{-\frac{s+s_0}{2}}
X^{2}(t).
\end{split}
\end{equation*}
On the other hand,  when $t\geq2,$  we  have
\begin{equation*} \label{s14ii-B-E}
\begin{split}
&\int_0^t\langle t-\tau\rangle^{-\frac{s+s_0}{2}}\|\mu^{+}h_{+}(c^{+},c^{-})\partial_{j}c^{+}\partial_{j}u^{+h}_{i}
(\tau)\|_{\dot B^{-s_0}_{2,\infty}}\,d\tau\\
&\quad\lesssim\int_0^1\langle t-\tau\rangle^{-\frac{s+s_0}{2}}
\|c^{+}\|_{\dot B^{\frac{N}{p}}_{p,1}}
\|\nabla u^{+}\|_{\dot B^{\frac{N}{p}}_{p,1}}^{h}\,d\tau\\&\qquad+\int_1^t\langle t-\tau\rangle^{-\frac{s+s_0}{2}}
\|c^{+}\|_{\dot B^{\frac{N}{p}}_{p,1}}
\|\nabla u^{+}\|_{\dot B^{\frac{N}{p}}_{p,1}}^{h}\,d\tau
\\&\eqdefa IV_{1}+IV_{2}.
\end{split}
\end{equation*}
Remembering the definition of $X(t)$ and $D(t)$ implies that
\begin{equation*}\label{4.15-EE}
IV_{1}\lesssim\langle t\rangle^{-\frac{s_{0}}{2}-\frac s2}X^{2}(1).
\end{equation*}
To handle the term $IV_{2}$, thanks to \eqref{low.11} and  the fact that $\langle \tau\rangle\approx\tau$ when $\tau\geq1$, we deduce  that
\begin{equation*}\label{4.17-EE}
\begin{split}
IV_{2}&\lesssim\int_1^t\langle t-\tau\rangle^{-\frac{s+s_0}{2}}
\Big(\|c^{+}\|_{\dot B^{\frac{N}{p}-1}_{p,1}}^{\ell}+\|c^{+}\|_{\dot B^{\frac{N}{p}}_{p,1}}^{h}\Big)
\|\nabla u^{+}\|_{\dot B^{\frac{N}{p}}_{p,1}}^{h}\,d\tau\\
&\lesssim\int_1^t\langle t-\tau\rangle^{-\frac{s+s_0}{2}}
\Big(\|c^{+}\|_{\dot B^{\frac{N}{2}-1}_{2,1}}^{\ell}+\|c^{+}\|_{\dot B^{\frac{N}{p}}_{p,1}}^{h}\Big)
\|\nabla u^{+}\|_{\dot B^{\frac{N}{p}}_{p,1}}^{h}\,d\tau\\
&\lesssim \Big(D^{2}(t)+X^{2}(t)\Big)\int_1^t\langle t-\tau\rangle^{-\frac{s+s_0}{2}}\Big(\langle \tau\rangle^{-(\frac{s_0}{2}+\frac{N}{4}-\frac{1}{2}+\alpha)}+\langle \tau\rangle^{-2\alpha}\Big)\,d\tau\\
&\lesssim\langle t\rangle^{-\frac{s_{0}}{2}-\frac{s}{2}}\Big(D^{2}(t)+X^{2}(t)\Big).
\end{split}
\end{equation*}
Thus
\begin{equation*}
\begin{split}
\label{s24}
\int_0^t\|\mu^{+}h_{+}(c^{+},c^{-})\partial_{j}c^{+}\partial_{j}u^{+}_{i}
(\tau)\|_{\dot B^{-s_0}_{2,\infty}}\,d\tau\lesssim \langle t\rangle^{-({\frac{s+s_0}{2}})} \Big(D^{2}(t)+X^{2}(t)\Big).
\end{split}
\end{equation*}
Similarly,  we also obtain  the 	corresponding estimates of other terms
$\mu^{+}k_{+}(c^{+},c^{-})\partial_{j}c^{-}\partial_{j}u^{+}_{i},$\\
$\mu^{+}h_{+}(c^{+},c^{-})\partial_{j}c^{+}\partial_{i}u^{+}_{j},$
$\mu^{+}k_{+}(c^{+},c^{-})\partial_{j}c^{-}\partial_{i}u^{+}_{j},
\lambda^{+}h_{+}(c^{+},c^{-})\partial_{i}c^{+}\partial_{j}u^{+}_{j}~\text{and}~
\lambda^{+}k_{+}(c^{+},c^{-})\partial_{i}c^{-}\partial_{j}u^{+}_{j}$. Here, we omit them.\\
To  bound  the term $\mu^{+}l_{+}(c^{+},c^{-})\partial_{j}^{2}u_{i}^{+}$ in $H_{2}^{i}$. We decompose it into
\begin{equation*}
\mu^{+}l_{+}(c^{+},c^{-})\partial_{j}^{2}u_{i}^{+}
=\mu^{+}l_{+}(c^{+},c^{-})\partial_{j}^{2}(u_{i}^{+})^\ell
+\mu^{+}l_{+}(c^{+},c^{-})\partial_{j}^{2}(u_{i}^{+})^h,
\end{equation*}
where $l_{+}$ stands for some smooth function vanishing at $0$. The term  $\mu^{+}l_{+}(c^{+},c^{-})\partial_{j}^{2}(u_{i}^{+})^\ell$ may be treated  as $c^{+}\,\div (u^{+})^\ell$ in \eqref{s12}, that is,
\begin{equation*}
\begin{split}
&\int_0^t\langle t-\tau\rangle^{-\frac{s+s_0}{2}}\|\mu^{+}l_{+}(c^{+},c^{-})\partial_{j}^{2}(u_{i}^{+})^\ell\|_{\dot B^{-s_0}_{2,\infty}}\,d\tau\lesssim \langle t\rangle^{-\frac{s+s_0}{2}} \Big(D^{2}(t)+X^{2}(t)\Big).
\end{split}
\end{equation*}
To handle the term $\mu^{+}l_{+}(c^{+},c^{-})\partial_{j}^{2}(u_{i}^{+})^h$, if $2\leq p\leq N,$  it follows from \eqref{99111B},  \eqref{low.6.23} and Proposition \ref{p28}, that
\begin{equation*}
\begin{split}
&\int_0^t\langle t-\tau\rangle^{-\frac{s+s_0}{2}}\|\mu^{+}l_{+}(c^{+},c^{-})\partial_{j}^{2}(u_{i}^{+})^h
(\tau)\|_{\dot B^{-s_0}_{2,\infty}}\,d\tau\\
&\quad\lesssim\int_0^t\langle t-\tau\rangle^{-\frac{s+s_0}{2}}\|(c^{+},c^{-})
(\tau)\|_{\dot B^{\frac{N}{p}-1}_{p,1}}
\|\nabla^{2} u^{+h}(\tau)\|_{\dot B^{\frac{N}{p}-1}_{p,1}}\,d\tau.
\end{split}
\end{equation*}
When  $t\leq2$,  then $\langle t\rangle\approx1$ and $\langle t-\tau\rangle\approx1$ for $0\leq\tau\leq t\leq 2$. We have
\begin{equation*} \label{s14ii-A-D}
\begin{split}
&\int_0^t\langle t-\tau\rangle^{-\frac{s+s_0}{2}}\|\mu^{+}l_{+}(c^{+},c^{-})\partial_{j}^{2}(u_{i}^{+})^h
(\tau)\|_{\dot B^{-s_0}_{2,\infty}}\,d\tau\\
&\quad\lesssim\int_0^t\langle t-\tau\rangle^{-\frac{s+s_0}{2}}
\Big(\|(c^{+},c^{-})
(\tau)\|_{\dot B^{\frac{N}{p}-1}_{p,1}}^{\ell}+\|(c^{+},c^{-})
(\tau)\|_{\dot B^{\frac{N}{p}}_{p,1}}^{h}\Big)
\|u^{+}(\tau)\|_{\dot B^{\frac{N}{p}+1}_{p,1}}^{h}\,d\tau\\
&\quad\lesssim\int_0^t\langle t-\tau\rangle^{-\frac{s+s_0}{2}}
\Big(\|(c^{+},c^{-})
(\tau)\|_{\dot B^{\frac{N}{2}-1}_{2,1}}^{\ell}+\|(c^{+},c^{-})
(\tau)\|_{\dot B^{\frac{N}{p}}_{p,1}}^{h}\Big)
\| u^{+}(\tau)\|_{\dot B^{\frac{N}{p}+1}_{p,1}}^{h}\,d\tau\\
&\quad\lesssim \langle t\rangle^{-\frac{s+s_0}{2}}
X^{2}(t).
\end{split}
\end{equation*}
When $t\geq2,$
\begin{equation*} \label{s14ii-B-D}
\begin{split}
&\int_0^t\langle t-\tau\rangle^{-\frac{s+s_0}{2}}\|\mu^{+}l_{+}(c^{+},c^{-})\partial_{j}^{2}(u_{i}^{+})^h
(\tau)\|_{\dot B^{-s_0}_{2,\infty}}\,d\tau\\
&\quad\lesssim\int_0^1\langle t-\tau\rangle^{-\frac{s+s_0}{2}}
\|(c^{+},c^{-})
(\tau)\|_{\dot B^{\frac{N}{p}-1}_{p,1}}
\|\nabla^{2} u^{+h}(\tau)\|_{\dot B^{\frac{N}{p}-1}_{p,1}}\,d\tau\\&\qquad+\int_1^t\langle t-\tau\rangle^{-\frac{s+s_0}{2}}
\|(c^{+},c^{-})
(\tau)\|_{\dot B^{\frac{N}{p}-1}_{p,1}}
\|\nabla^{2} u^{+h}(\tau)\|_{\dot B^{\frac{N}{p}-1}_{p,1}}\,d\tau
\\&\eqdefa V_{1}+V_{2}.
\end{split}
\end{equation*}
Obviously,
\begin{equation*}\label{4.15-DD}
V_{1}\lesssim\langle t\rangle^{-\frac{s_{0}}{2}-\frac s2}X^{2}(1).
\end{equation*}
To deal with the term $V_{2}$,
based on the fact that $\langle \tau\rangle\approx\tau$ when $\tau\geq1$, we conclude,  according to \eqref{low.11}, that
\begin{equation*}\label{4.17-DD}
\begin{split}
V_{2}&\lesssim
\int_1^t\langle t-\tau\rangle^{-\frac{s+s_0}{2}}
\|(c^{+},c^{-})
(\tau)\|_{\dot B^{\frac{N}{p}-1}_{p,1}}
\|\nabla^{2} u^{+h}(\tau)\|_{\dot B^{\frac{N}{p}-1}_{p,1}}\,d\tau\\
&\lesssim\int_1^t\langle t-\tau\rangle^{-\frac{s+s_0}{2}}
\Big(\|(c^{+},c^{-})
(\tau)\|_{\dot B^{\frac{N}{p}-1}_{p,1}}^{\ell}+\|(c^{+},c^{-})
(\tau)\|_{\dot B^{\frac{N}{p}}_{p,1}}^{h}\Big)
\|\nabla u^{+}(\tau)\|_{\dot B^{\frac{N}{p}}_{p,1}}^{h}\,d\tau\\
&\lesssim\int_1^t\langle t-\tau\rangle^{-\frac{s+s_0}{2}}
\Big(\|(c^{+},c^{-})
(\tau)\|_{\dot B^{\frac{N}{2}-1}_{2,1}}^{\ell}+\|(c^{+},c^{-})
(\tau)\|_{\dot B^{\frac{N}{p}}_{p,1}}^{h}\Big)
\|\nabla u^{+}(\tau)\|_{\dot B^{\frac{N}{p}}_{p,1}}^{h}\,d\tau\\
&\lesssim \Big(D^{2}(t)+X^{2}(t)\Big)\int_1^t\langle t-\tau\rangle^{-\frac{s+s_0}{2}}\Big(\langle \tau\rangle^{-(\frac{s_0}{2}+\frac{N}{4}-\frac{1}{2}+\alpha)}+\langle \tau\rangle^{-2\alpha}\Big)\,d\tau\\
&\lesssim\langle t\rangle^{-\frac{s_{0}}{2}-\frac{s}{2}}\Big(D^{2}(t)+X^{2}(t)\Big),
\end{split}
\end{equation*}
If $N<p\leq 2N,$  employing \eqref{99111B}, \eqref{prolow11} and  Proposition \ref{p28} yields that
\begin{equation*}
\begin{split}
&\int_0^t\langle t-\tau\rangle^{-\frac{s+s_0}{2}}\|\mu^{+}l_{+}(c^{+},c^{-})\partial_{j}^{2}(u_{i}^{+})^h(\tau)\|_{\dot B^{-s_0}_{2,\infty}}\,d\tau\\
&\quad\lesssim\int_0^t\langle t-\tau\rangle^{-\frac{s+s_0}{2}}
\Big(\|(c^{+},c^{-})\|_{\dot B^{1-\frac{N}{p}}_{p,1}}+\|(c^{+},c^{-})(\tau)\|_{\dot B^{\frac{N}{p}}_{2,1}}^{\ell}\Big)
\|\nabla^{2} u^{+h}(\tau)\|_{\dot B^{\frac{N}{p}-1}_{p,1}}\,d\tau.
\end{split}
\end{equation*}
When $t\leq2$,  then $\langle t\rangle\approx\langle t\rangle\approx\langle t-\tau\rangle\approx1$ for $0\leq\tau\leq t\leq 2$. We get
\begin{equation*} \label{s14ii-A-F}
\begin{split}
&\int_0^t\langle t-\tau\rangle^{-\frac{s+s_0}{2}}\|\mu^{+}l_{+}(c^{+},c^{-})\partial_{j}^{2}(u_{i}^{+})^h(\tau)\|_{\dot B^{-s_0}_{2,\infty}}\,d\tau\\
&\quad\lesssim\int_0^t\langle t-\tau\rangle^{-\frac{s+s_0}{2}}
\Big(\|(c^{+},c^{-})\|_{\dot B^{1-\frac{N}{p}}_{p,1}}^{\ell}+\|(c^{+},c^{-})\|_{\dot B^{1-\frac{N}{p}}_{p,1}}^{h}+\|(c^{+},c^{-})(\tau)\|_{\dot B^{\frac{N}{p}}_{2,1}}^{\ell}\Big)
\| u^{+}\|_{\dot B^{\frac{N}{p}+1}_{p,1}}^{h}\,d\tau\\
&\quad\lesssim\int_0^t\langle t-\tau\rangle^{-\frac{s+s_0}{2}}
\Big(\|(c^{+},c^{-})\|_{\dot B^{\frac{N}{p}-1}_{p,1}}^{\ell}+\|(c^{+},c^{-})\|_{\dot B^{\frac{N}{p}}_{p,1}}^{h}+\|(c^{+},c^{-})(\tau)\|_{\dot B^{\frac{N}{p}}_{2,1}}^{\ell}\Big)
\| u^{+}\|_{\dot B^{\frac{N}{p}+1}_{p,1}}^{h}\,d\tau\\
&\quad\lesssim\int_0^t\langle t-\tau\rangle^{-\frac{s+s_0}{2}}
\Big(\|(c^{+},c^{-})\|_{\dot B^{\frac{N}{2}-1}_{2,1}}^{\ell}+\|(c^{+},c^{-})\|_{\dot B^{\frac{N}{p}}_{p,1}}^{h}+\|(c^{+},c^{-})(\tau)\|_{\dot B^{\frac{N}{p}}_{2,1}}^{\ell}\Big)
\| u^{+}\|_{\dot B^{\frac{N}{p}+1}_{p,1}}^{h}\,d\tau\\
&\quad\lesssim \langle t\rangle^{-\frac{s+s_0}{2}}
\Big(D^{2}(t)+X^{2}(t)\Big),
\end{split}
\end{equation*}
where we have used $\frac{N}{p}-1<1-\frac{N}{p}\leq\frac{N}{p}$ for $N<p\leq 2N$. When $t\geq2,$  we have
\begin{equation*} \label{s14ii-B-F}
\begin{split}
&\int_0^t\langle t-\tau\rangle^{-\frac{s+s_0}{2}}\|\mu^{+}l_{+}(c^{+},c^{-})\partial_{j}^{2}(u_{i}^{+})^h(\tau)\|_{\dot B^{-s_0}_{2,\infty}}\,d\tau\\
&\quad\lesssim\int_0^1\langle t-\tau\rangle^{-\frac{s+s_0}{2}}
\Big(\|(c^{+},c^{-})\|_{\dot B^{1-\frac{N}{p}}_{p,1}}+\|(c^{+},c^{-})(\tau)\|_{\dot B^{\frac{N}{p}}_{2,1}}^{\ell}\Big)
\|\nabla^{2} u^{+h}(\tau)\|_{\dot B^{\frac{N}{p}-1}_{p,1}}\,d\tau\\&\qquad+\int_1^t\langle t-\tau\rangle^{-\frac{s+s_0}{2}}
\Big(\|(c^{+},c^{-})\|_{\dot B^{1-\frac{N}{p}}_{p,1}}+\|(c^{+},c^{-})(\tau)\|_{\dot B^{\frac{N}{p}}_{2,1}}^{\ell}\Big)
\|\nabla^{2} u^{+h}(\tau)\|_{\dot B^{\frac{N}{p}-1}_{p,1}}\,d\tau
\\&\eqdefa VI_{1}+VI_{2}.
\end{split}
\end{equation*}
It follows from the definition of $X(t)$ and $D(t)$,  that
\begin{equation*}\label{4.15-FF}
VI_{1}\lesssim\langle t\rangle^{-\frac{s_{0}}{2}-\frac s2}\Big(D^{2}(1)+X^{2}(1)\Big).
\end{equation*}
To bound the term $VI_{2}$, according to \eqref{low.8}, \eqref{low.11} and \eqref{low.12},
and using the fact that $\langle \tau\rangle\approx\tau$ when $\tau\geq1$, we conclude  that
\begin{equation*}\label{4.17-EF}
\begin{split}
VI_{2}&\lesssim\int_0^t\langle t-\tau\rangle^{-\frac{s+s_0}{2}}
\Big(\|(c^{+},c^{-})\|_{\dot B^{1-\frac{N}{p}}_{p,1}}^{\ell}+\|(c^{+},c^{-})\|_{\dot B^{1-\frac{N}{p}}_{p,1}}^{h}+\|(c^{+},c^{-})(\tau)\|_{\dot B^{\frac{N}{p}}_{2,1}}^{\ell}\Big)
\| u^{+}\|_{\dot B^{\frac{N}{p}+1}_{p,1}}^{h}\,d\tau\\
&\quad\lesssim \Big(D^{2}(t)+X^{2}(t)\Big)\int_1^t\langle t-\tau\rangle^{-\frac{s+s_0}{2}}\langle \tau\rangle^{-\min\{\frac 12+\alpha,2\alpha,\frac{3N}{2p}-\frac{N}{4}+\alpha\}}\,d\tau\\
&\quad\lesssim\langle t\rangle^{-\frac{s_{0}}{2}-\frac{s}{2}}\Big(D^{2}(t)+X^{2}(t)\Big).
\end{split}
\end{equation*}
Hence
\begin{equation*} \label{ds2l}
\begin{split}
\int_0^t\langle t-\tau\rangle^{-\frac{s+s_0}{2}}\|\mu^{+}l_{+}(c^{+},c^{-})\partial_{j}^{2}(u_{i}^{+})\|_{\dot B^{-s_0}_{2,\infty}}\,d\tau
\lesssim \langle t\rangle^{-\frac{s+s_0}{2}}
\Big(X^{2}(t)+D^{2}(t)\Big).
\end{split}
\end{equation*}
Similarly,
\begin{equation*}
\label{s1last}
\begin{split}
\int_0^t\langle t-\tau\rangle^{-\frac{s+s_0}{2}}\|(\mu^{+}+\lambda^{+})l_{+}(c^{+},c^{-})\partial_{i}\partial_{j}
u^{+}_{j}\|_{\dot B^{-s_0}_{2,\infty}}\,d\tau\lesssim \langle t\rangle^{-\frac{s+s_0}{2}} \Big(X^{2}(t)+D^{2}(t)\Big).
\end{split}
\end{equation*}
We finally conclude that
\begin{equation}\label{s2}
\int_0^t\langle t-\tau\rangle^{-\frac{s+s_0}{2}} \big\|H_{2}(\tau)\big\|_{\dot B^{-\frac{N}{2}}_{2,\infty}}^\ell d\tau
\lesssim\langle t\rangle^{-\frac{s+s_0}{2}}
\Big(X^{2}(t)+D^{2}(t)\Big).
\end{equation}
The term  $H_{4}$ may be treated along the same lines,  and we have
\begin{equation}\label{s4}
\int_0^t\langle t-\tau\rangle^{-\frac{s+s_0}{2}} \big\|H_{4}(\tau)\big\|_{\dot B^{-\frac{N}{2}}_{2,\infty}}^\ell d\tau
\lesssim\langle t\rangle^{-\frac{s+s_0}{2}}
\Big(X^{2}(t)+D^{2}(t)\Big).
\end{equation}

Thus, putting together inequalities \eqref{s1}, \eqref{s2}, \eqref{s3} and \eqref{s4}, we complete the proof of \eqref{s1234low1}. Further, combining  with \eqref{U} and \eqref{s1234low1}, we conclude that for all $t\geq0$ and $s\in(\varepsilon-\frac{N}{2},\frac{N}{2}+1],$
\begin{equation}
\label{low}
\langle t\rangle^{\frac{s+s_0}{2}}
 \|(c^{+},\,u^{+},\,c^{-},\,u^{-}\|^{\ell}_{\dot  B^{s}_{2,1}}\lesssim D_{0}+X^{2}(t)+D^{2}(t).
\end{equation}
\subsection{In the high frequencies}\ \ \ \ \
This part is devoted to bounding the last   term of $D(t)$. We first introduce the following system in terms of the weighted unknowns the term $(t^{\alpha}c^{+},t^{\alpha}u^{+},t^{\alpha}c^{-},t^{\alpha}u^{-} )$
\begin{equation}\label{equ:CTFS5}
\left\{
\begin{aligned}{}
&\p_t(t^{\alpha}c^{+})+\textrm{div}(t^{\alpha}u^{+})=\alpha t^{\alpha-1}c^{+}+  t^{\alpha}H_{1},\\
&\p_t(t^{\alpha}u^{+})+\beta_{1}\nabla (t^{\alpha}c^{+})
+\beta_{2}\nabla (t^{\alpha}c^{-})-\nu_{1}^{+}\Delta (t^{\alpha}u^{+})
-\nu_{2}^{+}\nabla\textrm{div}(t^{\alpha}u^{+})-\nabla \Delta (t^{\alpha}c^{+})\\&\qquad\qquad\qquad\qquad\quad=\alpha t^{\alpha-1}u^{+}+t^{\alpha}H_{2},
\\&\p_t(t^{\alpha}c^{-})+\textrm{div}(t^{\alpha}u^{-})=\alpha t^{\alpha-1}c^{-}+t^{\alpha}H_{3},
\\&\p_t(t^{\alpha}u^{-})+\beta_{3}\nabla(t^{\alpha} c^{+})
+\beta_{4}\nabla(t^{\alpha} c^{-})-\nu_{1}^{-}\Delta(t^{\alpha} u^{-})
-\nu_{2}^{-}\nabla\textrm{div}(t^{\alpha}u^{-})-\nabla \Delta (t^{\alpha}c^{-})\\&\qquad\qquad\qquad\qquad\quad=\alpha t^{\alpha-1}u^{-}+t^{\alpha}H_{4},\\
&(t^{\alpha}c^{+},t^{\alpha}u^{+},t^{\alpha}c^{-},t^{\alpha}u^{-})|_{t=0}=(0,0,0,0).
\end{aligned}
\right.
\end{equation}
Applying Lemma \ref{lemma3.4} and Proposition \ref{Pro:3}  to  system \eqref{equ:CTFS5}, by a similar derivation process of \eqref{high part111-A},  we  also have
\begin{equation}\label{regular estimate1}\begin{split}
\big\|\tau^{\alpha}\big(\nabla c^{+},u^{+},\nabla c^{-},u^{-}\big)\big\|
_{\wt L^\infty_t(\dot B^{\frac{N}{p}+1}_{p,1})}^{h}
&\lesssim\big\|\alpha\tau^{\alpha-1}(c^{+},u^{+}, c^{-},u^{-})\big\|
_{\wt L^\infty_t(\dot B^{\frac{N}{p}-1}_{p,1})}^{h}\\
&\qquad+\big\|\tau^{\alpha}(H_{1},H_{2},H_{3},H_{4})
\big\|_{\tilde{L}_{t}^{\infty}(\dot{B}^{\frac{N}{p}-1}_{p,1})}^{h}.
\end{split}\end{equation}
We now handle the lower order linear terms on the right hand-side of the above inequality. When $v\in \{c^{+},u^{+},c^{-},u^{-}\}$, for $0\leq \tau \leq t\leq2$,
we have $$\|\alpha\tau^{\alpha-1}v
\|_{\tilde{L}_{t}^{\infty}(\dot{B}^{\frac{N}{p}-1}_{p,1})}^{h}\lesssim \|v
\|_{\tilde{L}_{t}^{\infty}(\dot{B}^{\frac{N}{p}-1}_{p,1})}^{h}\lesssim X(t).$$
 When $t\geq2$, for  $0\leq \tau \leq 1$, we get
 $$\|\alpha\tau^{\alpha-1}v
\|_{\tilde{L}^{\infty}([0,1];\dot{B}^{\frac{N}{p}-1}_{p,1})}^{h}\lesssim \|v
\|_{\tilde{L}_{t}^{\infty}(\dot{B}^{\frac{N}{p}-1}_{p,1})}^{h}\lesssim X(t).$$
 When $t\geq2$, for  $1\leq \tau \leq t$, we have
\begin{equation}\label{regular estimate11}\begin{split}
\|\alpha\tau^{\alpha-1}v
\|_{\tilde{L}^{\infty}([1,t];\dot{B}^{\frac{N}{p}-1}_{p,1})}^{h}&=\alpha\sum_{j\geq j_0}2^{j(\frac{N}{p}+1)}2^{-2j}\|\tau^{\alpha-1}\dot{\Delta}_{j} v\|_{L^{\infty}([1,t],L^{p})}
\\&\lesssim \alpha 2^{-2j_0}\sum_{j\geq j_0}2^{j(\frac{N}{p}+1)}\|\tau^{\alpha}\dot{\Delta}_{j} v\|_{L^{\infty}([1,t],L^{p})}
\\&\lesssim \alpha 2^{-2j_0}\|\tau^{\alpha}v
\|_{\tilde{L}^{\infty}([1,t];\dot{B}^{\frac{N}{p}+1}_{p,1})}^{h}.
\end{split}\end{equation}
Choosing $j_0$ large enough such that  $$C\alpha 2^{-2j_0}\leq\frac{1}{4},$$ which implies that \eqref{regular estimate11} may be absorbed by  the  left hand-side of \eqref{regular estimate1}.
Thus
\begin{equation}\label{regular estimate111}\begin{split}
\big\|\tau^{\alpha}\big(\nabla c^{+},u^{+},\nabla c^{-},u^{-}\big)\big\|
_{\wt L^\infty_t(\dot B^{\frac{N}{p}+1}_{p,1})}^{h}\lesssim  X(t)+\big\|\tau^{\alpha}(H_{1},H_{2},H_{3},H_{4})
\big\|_{\tilde{L}_{t}^{\infty}(\dot{B}^{\frac{N}{p}-1}_{p,1})}^{h}.
\end{split}\end{equation}
It now comes down to estimating the above nonlinear terms.
We first show the following inequalities which are repeatedly used later.
\begin{equation}\label{regular estimate1111}\begin{split}
\|\tau^{\alpha}\nabla(c^{+},c^{-})\|
_{\wt L^\infty_t(\dot B^{\frac{N}{p}}_{p,1})}&\lesssim  \|\tau^{\alpha}\nabla(c^{+},c^{-})\|
_{\wt L^\infty_t(\dot B^{\frac{N}{p}}_{p,1})}^\ell+\|\tau^{\alpha}\nabla(c^{+},c^{-})\|
_{\wt L^\infty_t(\dot B^{\frac{N}{p}}_{p,1})}^{h}\\
&\lesssim  \|\tau^{\alpha}(c^{+},c^{-})\|
_{L^\infty_t(\dot B^{\frac{N}{2}+1-\varepsilon}_{2,1})}^\ell+\|\tau^{\alpha}\nabla(c^{+},c^{-})\|
_{\wt L^\infty_t(\dot B^{\frac{N}{p}}_{p,1})}^{h}
\\&\lesssim D(t),
\end{split}\end{equation}
\begin{equation}\label{regular estimate11111}\begin{split}
\|\tau^{\alpha}(u^{+},u^{-})\|
_{\wt L^\infty_t(\dot B^{\frac{N}{p}+1}_{p,1})}&\lesssim  \|\tau^{\alpha}(u^{+},u^{-})\|
_{\wt L^\infty_t(\dot B^{\frac{N}{p}+1}_{p,1})}^\ell+\|\tau^{\alpha}(u^{+},u^{-})\|
_{\wt L^\infty_t(\dot B^{\frac{N}{p}+1}_{p,1})}^{h}\\
&\lesssim  \|\tau^{\alpha}(u^{+},u^{-})\|
_{L^\infty_t(\dot B^{\frac{N}{2}+1-\varepsilon}_{2,1})}^\ell+\|\tau^{\alpha}(u^{+},u^{-})\|
_{\wt L^\infty_t(\dot B^{\frac{N}{p}+1}_{p,1})}^{h}
\\&\lesssim D(t).
\end{split}\end{equation}
For $\|\tau^{\alpha}H_{1}
\|_{\tilde{L}_{t}^{\infty}(\dot{B}^{\frac{N}{2}-1}_{2,1})}^{h}$, from \eqref{regular estimate1111}, \eqref{regular estimate11111} and Proposition \ref{p26},  we have
\begin{equation}\label{regular estimate6.54}\begin{split}
\|\tau^{\alpha}H_{1}
\|_{\tilde{L}_{t}^{\infty}(\dot{B}^{\frac{N}{p}-1}_{p,1})}^{h}&\lesssim  \|\tau^{\alpha}c^{+}\textrm{div}u^{+}\|
_{\wt L^\infty_t(\dot B^{\frac{N}{p}-1}_{p,1})}^h+\|\tau^{\alpha}u^{+}\cdot \nabla c^{+}\|
_{\wt L^\infty_t(\dot B^{\frac{N}{p}-1}_{p,1})}^{h}\\
&\lesssim \|\tau^{\alpha}c^{+}\textrm{div}u^{+}\|
_{\wt L^\infty_t(\dot B^{\frac{N}{p}}_{p,1})}^h+\|\tau^{\alpha}u^{+}\cdot \nabla c^{+}\|
_{\wt L^\infty_t(\dot B^{\frac{N}{p}-1}_{p,1})}^{h}\\
&\lesssim \|c^{+}\|_{\wt L^\infty_t(\dot B^{\frac{N}{p}}_{p,1})} \|\tau^{\alpha}u^{+}\|
_{\wt L^\infty_t(\dot B^{\frac{N}{p}+1}_{p,1})} +\|u^{+}\|
_{\wt L^\infty_t(\dot B^{\frac{N}{p}-1}_{p,1})}\|\tau^{\alpha}\nabla c^{+}\|
_{\wt L^\infty_t(\dot B^{\frac{N}{p}}_{p,1})}
\\&\lesssim X(t)D(t).
\end{split}\end{equation}
Similarly,
\begin{equation}\label{regular estimate6.55}\begin{split}
\|\tau^{\alpha}H_{3}
\|_{\tilde{L}_{t}^{\infty}(\dot{B}^{\frac{N}{p}-1}_{p,1})}^{h}\lesssim X(t)D(t).
\end{split}\end{equation}
In what follows, we bound the term $\|\tau^{\alpha}H_{2}
\|_{\tilde{L}_{t}^{\infty}(\dot{B}^{\frac{N}{2}-1}_{2,1})}^{h}$. To bound the first part of $H_{2}^{i}$, employing \eqref{regular estimate1111}, Proposition \ref{p26}, Proposition \ref{p28}  and \eqref{99111B},  we infer that
\begin{equation}
\begin{split}\label{regular estimate6.56}
&\|\tau^{\alpha}g_{+}(c^{+},c^{-})\partial_{i}c^{+}
\|_{\tilde{L}_{t}^{\infty}(\dot{B}^{\frac{N}{p}-1}_{p,1})}^{h}\\
&\quad\lesssim\|g_{+}(c^{+},c^{-})\|_{\tilde{L}_{t}^{\infty}(\dot{B}^{\frac{N}{p}-1}_{p,1})}\|\tau^{\alpha}\partial_{i}c^{+}
\|_{\tilde{L}_{t}^{\infty}(\dot{B}^{\frac{N}{p}}_{p,1})}\\
&\quad\lesssim \|(c^{+},c^{-})\|_{\tilde{L}_{t}^{\infty}(\dot{B}^{\frac{N}{p}-1}_{p,1})}
\|\tau^{\alpha}\nabla c^{+}
\|_{\tilde{L}_{t}^{\infty}(\dot{B}^{\frac{N}{p}}_{p,1})}\\
&\quad\lesssim X(t)D(t).
\end{split}
\end{equation}
Similarly,
\begin{equation}
\begin{split}\label{regular estimate6.57}
\|\tau^{\alpha}\tilde{g}_{+}(c^{+},c^{-})\partial_{i}c^{-}
\|_{\tilde{L}_{t}^{\infty}(\dot{B}^{\frac{N}{p}-1}_{p,1})}^{h}\lesssim X(t)D(t).
\end{split}
\end{equation}
To bound the term  $(u^{+}\cdot\nabla)u_{i}^{+}$, from Proposition \ref{p26} and \eqref{regular estimate11111}, we get
\begin{equation}
\begin{split}\label{regular estimate6.58}
\|\tau^{\alpha}(u^{+}\cdot\nabla)u_{i}^{+}
\|_{\tilde{L}_{t}^{\infty}(\dot{B}^{\frac{N}{p}-1}_{p,1})}^{h}&\lesssim \|u^{+}
\|_{\tilde{L}_{t}^{\infty}(\dot{B}^{\frac{N}{p}-1}_{p,1})}\|\tau^{\alpha}\nabla u^{+}
\|_{\tilde{L}_{t}^{\infty}(\dot{B}^{\frac{N}{p}}_{p,1})}\\&\lesssim X(t)D(t).
\end{split}
\end{equation}
Using \eqref{regular estimate11111}, Proposition \ref{p26}, Proposition \ref{p28}  and \eqref{99111B} yields that
 \begin{equation}
\begin{split}\label{regular estimate6.59}
&\|\tau^{\alpha}\mu^{+}h_{+}(c^{+},c^{-})\partial_{j}c^{+}\partial_{j}u^{+}_{i}
\|_{\tilde{L}_{t}^{\infty}(\dot{B}^{\frac{N}{p}-1}_{p,1})}^{h}\\
&\quad\lesssim\|\nabla c^{+}
\|_{\tilde{L}_{t}^{\infty}(\dot{B}^{\frac{N}{p}-1}_{p,1})}\|\tau^{\alpha}\nabla u^{+}
\|_{\tilde{L}_{t}^{\infty}(\dot{B}^{\frac{N}{p}}_{p,1})}\\
&\quad\lesssim X(t)D(t).
\end{split}
\end{equation}
Similarly,  we also obtain  the 	corresponding estimates of other terms
$\mu^{+}k_{+}(c^{+},c^{-})\partial_{j}c^{-}\partial_{j}u^{+}_{i},$\\
$\mu^{+}h_{+}(c^{+},c^{-})\partial_{j}c^{+}\partial_{i}u^{+}_{j},$
$\mu^{+}l_{+}(c^{+},c^{-})\partial_{j}^{2}u_{i}^{+},$
$(\mu^{+}+\lambda^{+})l_{+}(c^{+},c^{-})\partial_{i}\partial_{j}
u^{+}_{j},$
$\mu^{+}k_{+}(c^{+},c^{-})\partial_{j}c^{-}\partial_{i}u^{+}_{j},$
$\lambda^{+}h_{+}(c^{+},c^{-})\partial_{i}c^{+}\partial_{j}u^{+}_{j}~\text{and}~
\lambda^{+}k_{+}(c^{+},c^{-})\partial_{i}c^{-}\partial_{j}u^{+}_{j}$. Here, we omit the details.

Then
\begin{equation}\label{regular estimate6.614}\begin{split}
\big\|\tau^{\alpha}\big(\nabla c^{+},u^{+},\nabla c^{-},u^{-}\big)\big\|
_{\wt L^\infty_t(\dot B^{\frac{N}{p}+1}_{p,1})}^{h}
\lesssim  X(t)+X(t)D(t),
\end{split}\end{equation}
which together  with \eqref{low}   for all $t\geq0,$  yields that
\begin{align*}D(t)\lesssim  D_{0}+X(t)+X^2(t)+D^2(t).
\end{align*}
As Theorem \ref{th:main1} ensures that $X(t)\lesssim  X(0)$  with $X(0) $ being small, and $X(0)^{\ell}=\|( c_{0}^{+}, u_{0}^{+}, c_{0}^{-}, u_{0}^{-})\|^{\ell}_{\dot B^{\frac N2-1}_{2,1}}\lesssim \|( c_{0}^{+}, u_{0}^{+}, c_{0}^{-}, u_{0}^{-})\|^{\ell}_{\dot B^{-s_0}_{2,\infty}}$,   one can conclude that \eqref{1.8} is fulfilled for all time if $\big\|\big(\nabla R^{+}_{0},\,u^{+}_{0},\,\nabla R^{-}_{0},\,u^{-}_{0}\big)\big\|^h_{\dot B^{\frac Np-1}_{p,1}}$ and  $D_{0}$ are small enough. This completes the proof of Theorem \ref{th:decay}.
\begin{center}

\end{center}
\end{document}